\documentclass[a4paper]{amsart}  

\usepackage{amssymb} \usepackage{amscd} \usepackage{times}

\usepackage{amsmath}

\usepackage[french]{babel}

\def\zbb{\mathbb{Z}}  
  
  \def\phi{\varphi}
 \def\p1{{\mathbb{P}^1_\zbb}}

\newcommand{\be} {\begin{equation}}
\newcommand{\ee} {\end{equation}}

\title{ Cas d'existence de solutions d'EDP.}

\author{Samy Skander Bahoura}

\address{Equipe d'Analyse Complexe et Geometrie, Universite Pierre et Marie Curie, 75005, Paris, France.}

\email{samybahoura@gmail.com}

\begin{document}

\maketitle

\section{Cas d'existence de solutions d'EDP}

\begin{abstract}
 Estimations uniformes pour l'equation de Yamabe, l'equation du type Yamabe, l'equation de la courbure scalaire prescrite, l'equation du type courbure scalaire prescrite, l'equation de Gauss, l'equation de type Gauss. On se pose la question de savoir s'il y a des solutions.
\end{abstract}

\bigskip

{\bf En dimension 2 :}

\smallskip
 
1) Les solutions de Liouville en fonction des fonctions holomorphes dont les d\'erivees ne s'annulent pas, voir le livre de C. Bandle.

2) Sur la boule unite, il y a des solutions radiales, nulle au bord, cela revient a resoudre une equation differentielle.

3) Si le domaine n'est plus une boule, par la methode des sur et sous solutions, il y a des solutions au probleme de Gelfand, $ u=0 $ est sous solution et la  fonction distance  est une sur solution.(Voir Dupaigne).

4) Par la methode variationelle, il y a toujours des solutions, c'est li\'e \`a la compacit\'e de l'injection de Moser-Trudinger dans $ L^1$.

Pour la dimension 2, sur des surfaces compactes sans bord, il y a des resultats d'existence dans le livre d'Aubin, par la methode variationelle si la courbure scalaire $ S=-1 $ on peut minorer les solutions, lorsque $ S=0 $ il se peut qu'on ne puisse pas minorer les solutions.

Pour une courbure scalaire negative, on peut resoudre le probleme de la courbure scalaire prescrite avec une courbure prescrite $ f \leq 0 $ par la methode des sur et sous solutions, puis on rajoute une constante ou fonction petite $ f+c $ ou $ f+g $, $ c,g >0 $. Dans ce cas, il se peut qu'on puisse minorer les solutions. C'est fait dans le livre d'Aubin.

Dans le cas nul (courbure scalaire nulle, par la methode variationelle), on peut resoudre le probleme de la courbure prescrite, pour une courbure prescrite $ f $ changeant de signe et $ \int f <0 $.Dans ce cas il se peut que les solutions ne soient pas minorables. C'est le probleme de Kazdan-Warner. C'est fait dans le livre d'Aubin. 

En ce qui nous concerne, on se placera sur la partie positive de $ f+c= c >0 $ et $ f+g = g >0 $.

On peut resoudre le Probleme de la courbure scalaire sur une surface compacte non orientable, c'est fait dans le livre d'Aubin  pour le Projectif r\'eel de dimension 2, $ {\mathbb P}_2 $.

\bigskip

{\bf En dimension $ \geq  3 $.}

\smallskip

1) Pour l'\'equation de Yamabe (5 et 6) et l'\'equation de la courbure prescrite en dimension 4 :
Quand la vari\'et\'e est compacte sans bord, les r\'esultats existent par Aubin-Schoen en dimension 5 et 6.

Si on considere une fonction $ f $ changeant de signe et que la courbure scalaire $ S_g=0 $ ( et de metrique $ g $), par la methode variationelle $ f $ est courbure scalaire d'une metrique conforme, voir le livre de Hebey (pour les dimensions 3, 4, le theoreme d'Escobar Schoen) et le livre d'Aubin (le Theoreme d'Escobar-Schoen dans le cas conformement plat et le Theoreme de Bismuth). Voir aussi dans le livre d'Aubin d'autres conditions sur  $ f $ pour qu'elle soit courbure scalaire d'une metrique conforme. Dans le cas, changeant de signe on se placera sur  l'ensemble $ \{f>0\} $. Ceci est valable en toute dimension $ \geq 3 $ (conformement plat, Escobar-Schoen et Aubin et invariant de Yamabe positif).

Pour la dimension 6, on r\'esout le probl\`eme de la courbure scalaire prescrite sur une vari\'et\'e de courbure scalaire $ S_g=-1 $ (de metrique $ g $), pour une courbure prescrite $ f \leq 0 $, par exemple une cutoff avec $ {f=0} $ une boule, puis par un th\'eor\`eme dans le livre d'Aubin, il existe un voisinage  de $ f $  dans $ C^{\alpha} $, dans lequel on peut r\'esoudre le probl\`eme de la courbure prescrite.

Voir dans le livre d'Aubin la condition de Rauzy, Ouyang, Tang, Vazquez-Veron, la valeur propre du Laplacien doit etre plus grande que la courbure scalaire pour resoudre dans le cas n\'egatif le probleme de la courbure scalaire prescrite avec une courbure prescrite $ f $ cutoff. Or en prenant des petites boules, on a grace au Theoreme de Levy-Brhul, le comportement asymptotique des valeurs propres du Laplacien des petites boules geodesiques (elle le prouve dans le cas de varietes Riemanniennes analytiques), voir l'article et le monograph de Leon Karp-Mark Pinsky (ils le prouvent dans le cas de vari\'et\'es Riemanniennes $ C^{\infty} $).(C'est la condition de Rellich sur les valeurs propres, developement asymptotique).

 On r\'esout le probleme de la coubrure prescrite pour des fonctions cutoff $ f \leq 0 $, par exemple (de boules petites), puis on ajoute une constante positive petite pour avoir $ f+c $ courbure prescrite d'une m\'etrique conforme telle que $ {f+c=c>0} =$ boule et $ f+c $ est n\'egative dans un autre ensemble. Alors on remarque que la ou la solution est minimum, alors elle est minor\'ee par une constante strictement positive. On vient de construire un exemple de solutions minor\'ees par une constante positive, qui  est un exemple de solutions en dimension 6 pour la quelle le minimum est minor\'e par une constante positive  et est solution de Yamabe dans la boule ou $ f+c = c >0 $. (le minorant est $ m=1/(1-c) $, si on suppose, $ 0\geq f \geq -1 $, par exemple,( on peut remplacer $-1 $ par $ - \bar M <0 $ quelconque dans $ 0\geq f\geq - \bar M $),  utiliser l'\'equation et le fait qu'en $ P $ point minimum de la solution $ u $, $ \Delta u(P)\leq 0 $).(Journal of Functional Analysis). Sur les ensembles o\`u $ f<0 $ voir travail de Th\`ese dans le cas n\'egatif (courbure scalaire prescrite negative et unifor.holderienne).

Pour la dimension 4 et l'\'equation de la courbure prescrite, il suffit de suivre le meme chemin que pour la dimension 6, en remplacant $ c >0 $ par une fonction positive voisine de $ c $, $ c >0 $ petite, alors les solutions sont minor\'ees par une constante $ m >0 $. (Voir Travail de These dans le cas plat (JMPA) et Journal of Mathematical Analysis and Applications).

Quand la vari\'et\'e a un bord :

En dimension 5, il y a des solutions voir l'article de Z-C. Han et YY.Li, l'operateur doit etre de type positif, ce qui est possible si on consid\`ere des boules g\'eodesiques, la courbure moyenne est alors presque constante et positive (proche de de la sph\`ere), on choisit la courbure scalaire et de Ricci non constantes pour avoir un point non ombilique.(Journal of Functional Analysis).

2) Pour l'\'equation du type Yamabe : si  $ R=\dfrac{(n-2)}{4(n-1)} S_g $ avec $ S_g $ est la courbure scalaire et $ g $ la metrique.

Quand la vari\'et\'e est compacte sans bord, \`a partir de la dimension 4, il y a les r\'esultats d'Aubin pour un op\'erateur du type ($ \Delta +a $, avec $ a < R $ en un point).(Bulletin des Sciences Math\'ematiques).

Quand la vari\'et\'e a un bord, \`a partir de la dimension 5, il y a les r\'esultats de Holcman pour un op\'erateur du type ($ \Delta +a $, avec $ a < R $ en un point interieur).(Bulletin des Sciences Math\'ematiques).

Il y a l'article de Hebey-Vaugon (nodales solutions, 1994). (JMPA (dimension 4 cas plat avec 2 contraintes), Journal of Functional Analysis, Bulletin des Sciences Mathematiques et Mathematica Aeterna, Journal of Mathematical Analysis and Applications (dimension 4), ici on a des exemples avec ou sans minimum tendant vers zero).

3) Pour l'\'equation du type courbure scalaire prescrite :

Dans le cas sous-critique tendant vers le critique, il y a des solutions par la m\'ethode variationnelle car l'injection de Sobolev est compacte dans ce cas.(Travail de Th\`ese).

4) Dans le cas radial, il y a des exemples dans l'article de C.C.Chen- C.S. Lin, propositions 4.2 et 4.3 (de: Comm. Pure. Applied. Math. 1997). avec des courbures prescrites de la forme $ V=1-Kr^{\rho} $, $ 1< \rho < n-2 $. (Voir travail de These, pour les fonctions radiales).

Il y a les r\'esultats de Brezis-Nirenberg dans le cas d'un ouvert de l'espace euclidien en dimension $ \geq 4 $, pour une perturbation non lin\'eaire sous critique de l'\'equation de Yamabe dans le cas d'un ouvert plat. (Voir travail de Th\`ese).

Il y a le r\'esultat de Druet, sur une vari\'et\'e compacte sans bord avec potentiel $ f >0 $ et une perturbation sous-critique $ hu^{q-1} $, $ h >0, q < \frac{2n}{n-2}, n\geq 5 $, pour des valeures de $ q $ dependant d'un $ q_c=2 $, on peut se placer localement et dans le cas plat ou non plat et dans le cas plat se ramener a un ouvert de l'espace euclidien. (Voir travail de Th\`ese et equations avec perturbation sous-critique dans le cas Riemannien).

Pour l'equation de la courbure scalaire prescrite en dimension 4, on en a parl\'e ci-dessus, il suffit de considerer une courbure prescrite changeant de signe de telle mani\`ere que les solutions soient positives et minor\'ees et se placer sur la partie positive de la courbure prescrite. (Voir Travail de These dans le cas plat (JMPA) et Journal of Mathematical Analysis and Applications).

Pour la dimension 2, les estimations a priori ont pour but de prouver l'existence de solutions par la m\'ethode du degr\'e topologique de Leray-Schauder, comme dans De Figueiredo-Lions-Nussbaum, une des conditions est que l'operateur compact $ M $ est nul en 0,  ici cette condition n'est pas satisfaite car il y a une exponentielle. Cette m\'ethode peut s'appliquer \`a l'exponentielle (Equations en $\sinh $, plasmas et 2d-turbulences, avec poids), et aux problemes des valeurs propres non-lineaires.(Crandall-Rabinowitz prouvent l'existence d'une solution $ u_{\lambda} $ stable et la condition de stabilit\'e implique que $(u_{\lambda}) $ est born\'ee dans $ C^0 $. De Figueiredo-Lions-Nussbaum prouvent en utilisant le degr\'e topologique qu'il y a une deuxieme solution $ v_{\lambda} $ en utilisant le fait que $ (u_{\lambda})$ est born\'ee dans $ C^0 $, l'estimation a priori etait deja donn\'e dans la condition de stablit\'e de Crandall-Rabinowitz, par contre on ne sait pas si $ (v_{\lambda})$ est born\'ee dans $ C^0 $, c'est en partie l'objet de ce qu'on fait ici, borner les solutions et en particulier $ (v_{\lambda}) $).

D'autre part ce cas particulier (l'exponentielle) peut induire l'existence de solutions par la m\'ethode des sur et sous solutions (voir  Dupaigne), en fait ici on s'interesse aux estimations a priori des solutions d'\'equations du type Gauss. On a une compacit\'e \`a la Gromov des suites de m\'etriques ou bien un theoreme de compacit\'e de fonctions (Ascoli ou d'autres theoremes, voir les livres de Brezis et Aubin, dans $ L^p $, fonctions born\'ees dans une certaine topologie, $ L^p, H^1 $ par exemple et donc compacit\'e faible) ou si on revient aux EDP, des estimations a priori elliptiques.

Pour les dimensions $ \geq 3 $, on a des in\'egalit\'es du type Harnack, une particularit\'e de ces 
\'equations \`a exposant critique de Sobolev comme pour les fonctions harmoniques. C'est justifi\'e par l'absence de conditions au bord.

Dans le cas d'une vari\'et\'e compacte sans bord, un des but des estimations du type $ \sup \times \inf $ est de prouver qu'il n'y a que des blow-up isol\'es simples pour les solutions de ces \'equations et donc elles convergent vers uns fonction semblable \`a la fonction de Green, le Th\'eor\`eme de la masse positive s'il est vrai tout le temps, permet de dire que le terme constant est positif, alors que par l'identit\'e de Pohozaev l'integrale d'une certaine quantit\'e serait \`a la fois positive ou nulle (blow-up isol\'es simples) et strictement n\'egative (masse positive), ce qui n'est pas possible. Donc il n'y a pas de blow-up, c'est-\`a-dire que les solutions sont born\'ees dans $ C^2 $ et par un argument d'homotopie et degr\'e topologique le degr\'e serait non nul ce qui veut dire qu'il y a des solutions.(Un r\'esultat d'existence par le degr\'e topologique et par consequent l'estimation a priori).

L'existence de solutions peut etre connu par la m\'ethode variationnelle, on se borne \`a consid\'erer que ce travail (en partie) a pour but de chercher des estimations a priori et en particulier des inegalit\'es du type Harnack, propri\'et\'es de ces solutions \`a exposant critique de Sobolev.

On suppose que les solutions existent, ce qui est possible dans les diff\'erents cas consid\'er\'es et on prouve certaines propri\'et\'es qualitatives de ces \'equations. (Convergence de m\'etriques, \`a la Cheeger-Gromov, estimations a priori, compacit\'e de solutions, blow-up, in\'egalit\'es du type Harnack, in\'egalit\'es de Sobolev, r\'esultat d'unicit\'e et de rigidit\'e (Marc Herzlich a parl\'e de th\'eoreme de rigidit\'e, condition sur la courbure scalaire, comme pour les theoremes de masse positive, voir son expos\'e dans Numdam.), on a aussi, equilibre et stabilit\'e d'un systeme, dynamique d'un systeme, notions de blow-up et de diviseurs ou sans diviseurs (dimension 2 pour des domaines analytiques (Voir l'article de Marc Troyanov, le probleme est pos\'e sur des surfaces de Riemann, 'environnement' cas des domaines complexes) et equation de Yamabe et de courbure scalaire prescrite (Yamabe avec diviseurs, un exemple, voir la these de Farid Madani)), intervalles de precisions de nombres, branche du maximum en fonction du minimum, sup et inf li\'es ou majorants et minorants li\'es, il y a une notion d'enroulement car on part du minimum ou du minorant et on obtient le maximum ou le majorant). Du point de vue de la physique, de la chimie, de l'astronomie, de la biologie, amplitude de "fonction d'onde" et des emissions, dynamique de particules ou des astres, cordes (fermees ou ouvertes), D-branes en dimension 2 (cordes ouvertes avec conditions de Dirichlet), non-agr\'egation de cellules (non-concentration de cellules) (amibes). Ondes, gazs, combustion, equilibre gravitationnel d'\'etoiles, non-concentration d'amibes. Voir les articles de Crandall et Rabinowitz, De Figueiredo-Lions-Nussbaum, Chen-Li, Lin-Ni-Takagi.

Concernant les vari\'et\'es qu'on considere, elles peuvent etre non conform\'ement plates, il y a quelques exemples de ces vari\'et\'es.

1)	Les vari\'et\'es de dimension 1 et 2 sont localement conform\'ement plates.

2)	Les sph\`eres, sommes connexes de sph\`eres, produit de sph\`eres et d'un cerle, l'espace hyperbolique et le produit d'une variete de courbure sectionelle constante par un cercle, sont conformement plates.

Par le revetement, les projectifs r\'eels sont conformement plats car de courbure sectionelle constante egale a 1 (Voir le livre d'Hebey, \'egalit\'e des courbures sectionnelles, quand on a un revetement riemannien). Les projectifs de dimension paire sont non orientables, les projectifs de dimension impaire sont orientables (voir le livre de Lafontaine).

3)	Il existe une vari\'et\'e compacte de dimension 3, non conform\'ement plate, plus precisement, elle ne possede pas de structure localement conformement plate. (William Goldman).

Si on considere une variete de courbure sectionelle non constante et son produit par un cercle, alors la variete obtenue est non localement conformement plate.

4)	Le projectif complexe de dimension 2 complexe donc de dimension 4 r\'eelle est non conformement plat (m\'etrique de Fubini-Study), il est d'Einstein et de courbure scalaire constante positive.

5)	Les surfaces K3, non conform\'ement plat de dimension 4 r\'eelle. Il est d'Einstein de constante 0.

6)	Le produit $ S2 \times S2 $ est de dimension 4 r\'eelle et non conformement plat, d'Einstein, de courbure scalaire constante positive.

7)	Les projectifs complexes de dimensions $ > 2 $, non conform\'ement plat.(M\'etrique de Fubini-Study).

8)	Maintenant, si on veut que la courbure scalaire soit, -1, 0, 1, on utilise les proc\'ed\'es suivant :

-On construit des vari\'et\'es produit \`a partir des vari\'et\'es pr\'ec\'edentes, elles sont alors soit plates soit non plate (il suffit qu'une ne soit pas plate, voir le Hebey).

-On utilise les sommes connexes, en se r\'ef\'erant \`a un article de Dominic Joyce, sur les diff\'erentes combinaisons possibles des sommes connexes donnant des vari\'et\'es  de courbures scalaire -1, 0, 1. On construit de telles vari\'et\'es et elles peuvent etre plates ou non plates.

\smallskip

Exemples dans le cas negatif: quand la courbure scalaire prescrite est negative et uniformement holderienne. Voir le livre d'Aubin, cas de Yamabe negatif ou nul, probleme de Kazdan-Warner dans le cas nul et se placer sur la partie ou la courbure scalaire prescrite est negative et reguliere. Voir travail de these, cas negatif.

\smallskip

{\it Exemple en dimension 4, cas non plat :}

\smallskip

Soit $ S_2 $ la sphere de dimension 2. On consid\`ere une surface de genre grand, 2, la somme connexe de deux Tores $ T_2 \sharp T_2 $, on r\'esout l'equation de la courbure prescrite constante $ k =-1 $, c'est fait dans le livre d'Aubin. On peut supposer qu'il existe une surface $ S $ de courbure scalaire constante (et donc sectionelle constante $ -3/2 $) \'egale \`a -3 par exemple. On consid\`ere le produit $ S_2 \times S $, il est non conform\'ement plat par le th\'eor\`eme dans le livre d'Hebey (courbures sectionelles constantes non oppos\'ees) et de courbure scalaire -1. On peut alors appliquer le proc\'ed\'e d'existence de solutions pour l'\'equation de la courbure prescrite en dimension 4 minor\'ees par $ m>0 $ (pour la fonction cutoff nulle au voisinage du point $ P $ ou $ Weyl(P) \not = 0 $, est non nul). On obtient alors des solutions de l'\'equation de la courbure prescrite en dimension 4 minor\'ees par une constante positive et sur une vari\'et\'e non localement plate.

\smallskip

{\it Exemples en dimensions 5 et 6, cas non plat:}

\smallskip

On peut faire la meme chose pour l'\'equation de Yamabe en dimension 5 et 6. On considere $ S_1 \times S_2 \times S $ et $ S_1 \times S_1 \times S_2 \times S $ qui sont non plate et de courbure scalaire -1, ici, $ S_1 $ est le cercle unit\'e.

\smallskip

{\it Un autre exemple en dimension 5, cas non plat:}

\smallskip

Soit $ S_3 $ la sphere de dimension 3 (courbure scalaire 6 et courbure sectionelle 1) et $ S $ une surface compacte de courbure scalaire -6 (donc sectionelle -3). Alors le produit $ S_3 \times S $ est non localement plat (voir theoreme du livre d'Hebey, courbures sectionelles constantes non opposees) et de courbure scalaire $ R=0 $. En considerant une fonction changeant de signe et cutoff $ f=\epsilon >0 $ dans une petite boule centree en $ P $ tel que $ Weyl(P) \not = 0 $ et $ f <0 $ ailleurs telle que $ \int f <0 $ ($\epsilon $ petit). On peut appliquer le theoreme de Bismuth en dimension 5 d'existence de solutions de l'equation de la coubure prescrite $ f $ changeant de signe (avec une courbure scalaire nulle sur une variete compacte). Ici on n'a pas forcement un minorant du minimum des solutions.

On a un exemple de solutions de l'equation de Yamabe sur $ \{f=\epsilon>0\} $, sur une variete non conformement plate de dimension 5 ou on n'a pas forcement un minorant des solutions.

\smallskip

{\it Un autre exemple en dimension 5:}

Un theoreme dans le livre d'Aubin dit que si on considere une vari\'et\'e $ M $ de dimension $ n\geq 3 $ non conformement diffeomorphe a la sphere (en particulier non localment conformement plate), alors il existe $ k >1 $, ne dependant que de la vari\'et\'e $ M $, tel que toute fonction $ f >0 $ verifiant:

$$ \sup_M f \leq k \inf_M f, $$

est courbure scalaire d'une metrique conforme. On applique cela \`a $ M = S_1 \times P_2({\mathbb C}) $ (le cercle fois le projectif complexe de dimension 2, $ M $ est non localement conformement plate par le theoreme dans le livre d'Hebey, de courbure scalaire constante $ R >0 $, $ dim (M)=5 $), et, $ f $ une fonction cutoff telle que $ \dfrac{1+2\epsilon}{k}\leq f \leq 1+\epsilon $ et $ f\equiv 1+\epsilon $ dans une boule  centree en $ P $ telle que $ Weyl(P) \not =0 $, ici, $ S_1 $ est le cercle unit\'e.

Ici, on ne sait pas minorer les solutions.

\smallskip

{\it Un autre exemple en dimension 6, cas non plat:}

 Soit $ S_1 $ le cercle unit\'e et $ S_3 $ la sphere de dimension 3 (courbure scalaire 6 et courbure sectionelle 1) et $ S $ une surface compacte de courbure scalaire -6 (donc sectionelle -3). Alors le produit $ S_3 \times S $ est non localement conform\'ement plat (voir le livre d'Hebey). On prend alors $ S_1 \times S_3 \times S $ qui est une vari\'et\'e de dimension 6 non localement conform\'ement plate (voir le livre d'Hebey, dans un produit, des qu'une des vari\'et\'es est non plate, la vari\'et\'e produit est non plate) et de courbure scalaire 0. On applique le r\'esultat d'Aubin-Hebey et Bismuth, pour les vari\'et\'es compactes non localement conform\'ement plate de dimension 6 ayant une courbure scalaire nulle. Comme l'exemple en dimension 5. En considerant une fonction changeant de signe et cutoff $ f=\epsilon >0 $ dans une petite boule centree en $ P $ tel que $ Weyl(P) \not = 0 $ et $ f <0 $ ailleurs telle que $ \int f <0 $ ($\epsilon $ petit). On peut appliquer le theoreme d'Aubin-Hebey  et Bismuth en dimension 6 d'existence de solutions de l'equation de la coubure prescrite $ f $ changeant de signe (avec une courbure scalaire nulle $ S_g=0 $ sur une variete compacte $ (M, g) $). Ici on n'a pas forcement un minorant du minimum des solutions.
 Ici, on ne peut pas forc\'ement minorer le minimum et sur une vari\'et\'e non-localement conform\'ement plate. (voir le livre d'Aubin).

\smallskip

{\it Exemple en dimension 4, cas plat et equation de la courbure scalaire prescrite:}

\smallskip

Soit $ S_4 $ la sphere de dimension 4. On considere une variete de courbure scalaire -1 si elle est conformement plat, c'est fini, on utilise le meme procede que precedement dans le cas non plat. Si cette variete est non plate, par exemple $ S_2 \times S $ du cas non plat, on utilise la somme connexe avec la sphere $ S_4 $ de Dominic Joyce ; on obtient une variete de courbure scalaire -1 dans le quelle (la partie de $ S_4$ est plate, car on a soit la metrique de $ S_4$ soit une fonction fois la metrique de $S_4$ et le tenseur de Weyl reste nul car on a un invariant conforme). On resout le probleme de la courbure prescrite pour une fonction cutoff nulle au voisinage du point $P$ ou la metrique est plate et negative ailleurs puis on rajoute une fonction $ C^{0,1} $. Alors  on a des solutions minorees et comme la metrique est plate au voisinage de $ P $ et que c'est l'equation de la courbure prescrite, en faisant un changement de metrique conforme, on obtient l'equation sur un ouvert de $ {\mathbb R}^4 $ avec des solutions minorees.

\smallskip

Davis, M, W, (1985),(Closed orientable hyperbolic 4-manifold), montre l'existence d'une vari\'et\'e compacte sans bord orientable hyperbolique de dimension 4, c'est-a-dire de courbure sectionelle constante -1, (en fait, comme pour le Tore, ces varietes sont obtenues par un revetement riemannien du au quotient de l'espace hyperbolique par un groupe discret d'isometries). Donc elle est localement conformement plate et de courbure scalaire constante strictement negative. On peut utiliser le meme procede d'existence de solutions minor\'ees par une constante positive $ m >0 $ et solutions de l'equation de la courbure prescrite. Par un changement de metrique conforme on se ramene a un ouvert de l'espace $ {\mathbb R}^4 $.

\smallskip

{\it D'autres exemples en dimension 4: (voir la section 5 pour la formulation en dimension 4)}

\smallskip

Ici on regroupe les deux exemples "plat" et "non plat" pr\'ec\'edents en un seul exemple.

\smallskip

La vari\'et\'e (par exemple) $ M =S\times S_2 $ de metrique $ g $ et courbure scalaire  $ S_g=-6 $ dans le cas non-plat. ($ P $ tel que $ Weyl(P)\not =0 $).

La vari\'et\'e (par exemple) $ M=(S\times S_2)\sharp S_4 $ ou $ M = Davis \, Manifold $, dans le cas plat.(de metrique $ g $ et courbure scalaire $ S_g=-6 $ (on multiplie la metrique par un scalaire adequat)). ($ P $ tel que $ Weyl(Q)\equiv 0 $ dans une boule $ B_r(P) $).

\bigskip

a) Pour un $ m >0 $ il existe dans le livre d'Aubin des exemples de $ (u,V) $ solutions de l'\'equation de la courbure scalaire prescrite en dimension 4 pour un $ V $ (Lipschitzien) verifiant: $||\nabla V||_{\infty} = o(1) \times m \dfrac{c(\bar M)\sqrt {c(\bar M)}}{32e^2\sqrt 2} \leq o(1) \times \dfrac{c(\bar M)\sqrt {c(\bar M)}}{32e^2\sqrt 2} (\inf_M u) $. (On suppose ici $\min_M u \geq m >0 $ et $ M $ devient la boule $ B_r(P) $). (la fonction cutoff $ -\bar M \leq f\leq 0 $ devient  $ -\bar M <-\bar M+c\leq f + c \leq c $, $ c=c(\bar M)>0 $ assez petit ($c<1 $ et $ \bar M\to +\infty $)). (On aura $ c(\bar M)/2 \leq V \leq c(\bar M) $ et  $ \inf u \geq 1/\sqrt{\bar M} $). ( si $ c(\bar M) \to 0 $, on remplace $ V $ par $ W=V/c(\bar M) $ et $ u $  par $ v=\sqrt {c(\bar M)} u $) .

\smallskip

b) En faisant varier $ m $, qui devient $ m_i \to 0 $ ($ m_i>0 $, on peut prendre $ m_i=\sqrt{c(\bar M)}/\sqrt{\bar M} $ avec $ \bar M \to +\infty $), on exhibe deux suites $ (u_i,V_i) $ (qui deviennent $ (v_i, W_i) $ si $ c(\bar M) \to 0 $) solutions de l'equation de la courbure scalaire prescrite avec $ ||\nabla V_i ||_{\infty} =o(1) \times m_i \leq o(1) \times (\inf_M u_i) $. Sauf qu'ici on ne sait pas si $ \min_M u_i \to 0 $ ou $ \min_M u_i \geq m'>0 $.

\smallskip

c) Le point b) pr\'ec\'edent  dit qu'on a deux suites $ (u_i,V_i) $ v\'erifiant les deux contraintes sur $ u $ et $ V $ avec le fait qu'on ne sait pas minorer le minimum des solutions par une constante $ m'>0 $. C'est l'exemple voulu pour cette formulation en dimension 4 avec deux contraintes sur $ u $ et $ V $. (les deux contraintes sont: 1) l'equation de la courbure scalaire prescrite en dimenson 4, verif\'ee par $ u>0 $ relativement a un $ V $ Lipschitzien et entre deux constantes positives uniform\'ement, 2) la deuxieme contrainte est entre $ ||\nabla V||_{\infty} $ et $ \inf_M u $; $ ||\nabla V||_{\infty} \leq \dfrac{3a {\sqrt a}}{32e^2{\sqrt 2}}\inf_M u $). (Un point essentiel est d'exhiber un exemple de fonctions $ (u_i,V_i) $ verifiant les deux contraintes sans rien savoir sur le minimum des solutions $ u_i $, c'est ce qu'on a fait jusqu'ici).

\smallskip

Le terme $ o(1) $ sert a attenuer l'effet de la presence des constantes universelles et du changement de metrique conforme dans le cas plat quand on se ramene \`a un ouvert de $ {\mathbb R}^4 $.

\smallskip

{\it Exemples en dimension 3, cas non plat et plat}

\smallskip

Soit $ M_{\phi} $ le Torus bundle de William Goldman. Cette variete compacte connexe orientable ne possede pas de structure localement plate. C'est \`a dire que toute metrique $ g $ sur cette vari\'et\'e est non localement plate. D'apres Aubin, il existe une metrique $ g_1 $ de courbure sclaire $ S_{g_1}=-1 $ et a partir de cette metrique, pour toute fonction $ f $ strictement  negative quelque part, il existe une metrique $ g_2 $ telle que $ S_{g_2}=f $. Il suffit de prendre $ f $ changeant de signe. 

Comme on peut resoudre le probleme de la courbure prescrite avec $ H $ une courbure prescrite cutoff (autour de $ P $ tel que $ C_{g_2}(P) \not = 0 $, en dimension 3, le tenseur de Weyl est nul, on prend le tenseur de Cotton, $ C_{g_2} $ caracterise la courbure et c'est un invariant conforme), et $ H \leq 0 $, pour une courbure scalaire de depart egale a $ f $ relativement a la metrique $ g_2 $, puis on rajoute a $ H $ une constante $ c >0 $ ou une fonction $ k >0 $, dans ce cas il se peut qu'on ne puisse pas minorer les solutions, car le point $ Q $ ou le minimum est atteint peut etre tel que $ f(Q) \geq 0 $.

\smallskip

{\bf Remarque:} Pour trouver des solutions de l'equation de la courbure prescrite, on considere, $ f\leq 0 $ cutoff alors $ \exists \, g_2 $ telle que $ S_{g_2}=f $, puis on resout le probleme de la courbure prescrite avec $ F\equiv -1 $, c'est fait dans le livre d'Aubin, de plus la solution $ u_1 $ est unique. On considere la metrique $ g_3=u_1^{4/(n-2)} g_2 $, sa courbure scalaire est $ -1 $. On resout le Probleme de la courbure prescrite (par les sous et sur solutions, c'est fait dans le livre d'Aubin), pour une courbure prescrite $ H $ cutoff, comme on l'a dit precedemment, puis on rajoute une constante positive $ c $ ou une fonction positive $ k $. Dans ce cas les solutions $ u_c $ sont minor\'ees par $ m >0 $. Donc les solutions au final sont du type $ u_1 u_c $ sont minor\'ees. Ce qu'on vient de voir c'est qu'il existe des solutions et elles sont minor\'ees par $ m >0 $ de :

$$ -\Delta u+fu=(H+c)u^{N-1},\,\, u >0, \qquad (E_c)$$

Par contre ce qu'on ne sait pas est que si toutes les solutions sont minorables par $ m >0 $.

Soit, $ u $ une solution de $ (E_c) $ alors elle n'est pas forcement minorable, car le point $ Q $ ou le minimum est attient peut etre tel que $ f(Q)=0 $. 

\smallskip

Dans le cas plat, il suffit de considerer $ S_1 \times S_1 \times S_1 $ qui est plat de courbure scalaire nulle (produit d'un Tore de courbure sectionelle nulle et d'un cercle, avec la metrique produit) et utiliser le theoreme du livre de Hebey, comme ce qu'on a fait pour la dimension 5. (Theoreme de Schoen-Escobar d'existence de solutions changant de signe).

\smallskip

{\it Un exemple en dimension 3:}

\smallskip

On resout le probleme de la courbure prescrite sur la sphere $ S_2 $ de dimension 2 avec la courbure prescrite $ R=1+(<x|e>)^2, R(-x)=R(x) $, (c'est un resultat de Moser, voir le livre d'Aubin), $ R $ est non constante et donc la courbure sectionelle associ\'ee est non constante. 

On considere le produit $ M_3=S_1\times S_2 $ avec la metrique precedente sur $ S_2 $ de courbure sectionelle non constante. Donc $ M_3 $ est non localement conformement plate par un resultat d'Hebey (courbure sectionelle de $ S_2 $ non constante).

On utilise le theoreme dans le livre d'Aubin, comme dans l'exemple 3 en dimension 5. ($ \sup f \leq k \inf f $, $ k >1 $). On choisit $ f >0 $ telle $ \nabla f $ non born\'e, tendant vers une fonction Heavside, car un theoreme de YY.Li-M.Zhu dit que tout est born\'e, si $ f $ est $ C^2 $ et en particulier les solutions sont minor\'ees).

Dans ce cas on n'a pas forcement un minorant des solutions de l'equation de Yamabe dans la boule ou $ C_g(P)\not = 0 $, (tenseur de Cotton), sur une variete non conforment plate.

On utilise le resultat d'Escobar-Schoen en dimension 3, pour une fonction $ f = \epsilon >0 $ cuttof changeant de signe ou bien une $ f $ une fonction changeant de signe ayant un point critique en $ P $ point maximum de $ f $ tel que $ C_g(P)\not = 0 $. (voir le livre d'Hebey).

Dans ce cas aussi, on ne peut pas minorer les solutions sur une variete compacte non localement conformement plate de dimension 3. 

Ici aussi, on ne sait pas minorer les solutions.

\smallskip

{\it Exemples en dimension 3: equation du type courbure scalaire prescrite:}

\smallskip

Soit $ M $ la vari\'et\'e de dimension 3 de William Goldman. Elle ne possede pas de structure localement conform\'ement plate. D'apres Kazdan-Warner, voir dans le livre d'Aubin, cette vari\'et\'e possede une metrique $ g $ de courbure scalaire $ S_g=-1 $ et elle est non localement conform\'ement plate. On prend $ -\bar M \leq f_0\leq 0 $ une fonction cutoff nulle autour de $ P $  un voisinage non localment conform\'ement plat (tenseur de Cotton non nul). On r\'esout le probleme de la courbure scalaire prescrite avec $ f_0 $ courbure prescrite. Apres, on sait qu'il existe une voisinage dans $ C^{\alpha}, \alpha \in ]0,1[ $ dans lequel on peut resoudre le probleme de la courbure prescrite. On prend par exemple $ \frac{c(\bar M)}{2} \leq f \leq c(\bar M) $ et $ ||\nabla f||_{\infty} \leq k c(\bar M), k \geq 0 $, $ 0< c(\bar M) <1 $, et dans $ B_r(P) $. Apres, voir dans le livre d'Aubin, un theoreme de Kazdan-Warner ( par sur et sous solutions), l'equation $ \Delta u+au=hu^5 $ possede une solution avec $ \frac{1}{8} S_g=-1/8\leq a <0 $ et $ h\leq f $, ($ \Delta =-\nabla^i\nabla_i $). On prend alors $ -1/8\leq a \leq -1/16 $ et $ h = f- \eta c(\bar M)/4 \leq f $, avec $ \eta $ une fonction cutoff dans $ B_r(P) $, positive ou nulle, egale a $ 1 $ dans $ B_{r/2}(P) $ et nulle en dehors de $ B_r(P) $. On obtient alors, $ c(\bar M)/4 \leq h\leq 3c(\bar M)/4 $ et $ ||\nabla h||_{\infty} \leq k c(\bar M) $. Si $ c(\bar M) \to 0 $, on remplace $ h $ par $ W=h/c(\bar M)$ et $ u=u_{\bar M} $ par $ v_{\bar M} = [c(\bar M)]^{1/4} u_{\bar M} $. La solution $ v=v_{\bar M} $ est minor\'ee par $ k_0 [c(\bar M)]^{1/4} /{(1+ \bar M)}^{1/4} \to 0 $ si $ \bar M \to +\infty $ ($ k_0 >0 $). Tout cela dans une boule $ B_{r/2}(P) $ avec $ C_g(P)\not =0 $, $ C_g $ tenseur de Cotton.

\smallskip

a) On pouvait prendre $ M=M_3\sharp (S\times S_1) $ la somme connexe de Dominic Joyce (au final on obtient une vari\'et\'e de courbure scalaire $-1 $  avec une partie  non plate), de la vari\'et\'e non localment conform\'ement plate $ M_3=S_2\times S_1 $ de l'exemple precedent et de la vari\'et\'e localement confrom\'ement plate $ S\times S_1 $ avec $ S $ une surface de courbure $ -1 $.

\smallskip

b) On pouvait prendre $ M=S\times S_1 $ avec $ S $ une surface de courbure scalaire non constante et $ <0 $ (on resout l'equation de la courbure prescrite sur une surface de courbure $ -1 $ et la courbure prescrite une cutoff (et $ <0 $) par exemple, la nouvelle metrique conforme est de courbure la fonction cutoff non constante). D'apres le livre de Hebey, elle est non conform\'ement plate.

\smallskip

c) Dans le cas plat, on prend $ M=S\times S_1 $, $ S $ une surface de coubure scalaire $-1 $.

\smallskip

{\it Exemples en dimensions $ n\geq 4 $ et $ n=4 $: equation du type Yamabe:}

\smallskip

On s'inspire de l'exemple precedent. On a de la meme maniere des exemples en dimensions $ n \geq 4 $ ($ n >3 $) pour des equations du type Yamabe. En dimension 4, du type $\Delta u+au=u^3, u >0 $ et $ \frac{1}{6} S_g=-1/6 \leq a<0 $ qui peuveut etre minorables et non necessairement minorables. La vari\'et\'e est $ M=S\times S_2 $ comme dans le premier exemple en dimension 4, non localement conformement plate ou, $ M=(S\times S_2)\sharp S_4 $ (de Dominic Joyce) ou $ M=Davis Manifold $ dans le cas plat.

\smallskip

{\it Exemples en dimensions $ n\geq 3, n=3 $: equation de courbure scalaire prescrite  et du type courbure scalaire prescrite}

\smallskip

On peut prendre l'exemple de Hebey-Vaugon (nodales solutions, 1994), pour les equations : $ \Delta_g u_{\epsilon} + a(x) u_{\epsilon} = f {u_{\epsilon}}^{N-2} u_{\epsilon} $, $ f >0 $, avec la condition au bord $ u_{\epsilon}=\lambda_{\epsilon} h $, $ h \not \equiv 0 $ et $ \lambda_{\epsilon} >0 $. On peut prendre $ u_{\epsilon} \geq 0 $ dans la preuve de Hebey-Vaugon, si $ h\geq 0 $ ($ h $ une fonction cuttof autour de $ x_0 \in \partial \Omega $, par exemple). Par le principe du maximum $ u_{\epsilon} >0 $ ($ \Delta_g =-\nabla_i (\nabla^i), N=\frac{2n}{n-2} $ et $ \Delta_g + a $ coercif, si $ n=3 $, on peut prendre par exemple, $ S_g\equiv 1 $ et $ a \leq \frac{1}{8} S_g=\frac{1}{8} $).

(En ecrivant dans la fonctionnelle de Hebey-Vaugon, $ J(u)=G(u+\tilde h)-\int_{\Omega} |\nabla_g \tilde h|^2 dV_g $, cela revient \`a prendre l'inf sur $ u\in W_0^{1,2}(\Omega) $ avec $ u + \tilde h \geq 0 $. On a si $ u_{\alpha,q} $ est un minimiseur alors $ v_{\alpha,q}=|u_{\alpha,q}+\tilde h|-\tilde h $ est aussi un minimiseur avec $ v_{\alpha,q}+\tilde h \geq 0 $.).

On peut considerer sur une vari\'et\'e compacte sans bord $ M $, un point $ P $ tel que, $ C_g(P)\not =0 $ ($ n=3 $), $Weyl(P) \not = 0 $($ n\geq 4$), cas non plat, ou $ C_g(Q)\equiv 0 $ ($ n=3 $), $ Weyl(Q) \equiv 0 $($ n\geq 4$), $ Q\in B_r(P) $ avec $r\leq \frac{inj_g(M)}{2} $, cas plat, avec $ \Delta +\frac{n-2}{4(n-1)} S_g $ coercif et puis prendre, $ \Omega = \exp_P^{-1}(B_{r/2}(P)) $ et $ \tilde g=\exp_P^*(g) $ et $ a\leq \frac{n-2}{4(n-1)}S_{\tilde g} $. Par exemple pour $ n =3 $, $ M=S_2\times S_1 $ dans le cas plat et $ M=S_2\times S_1 $ avec une metrique de courbure sectionnelle $ >0 $ et non constante sur $ S_2 $ pour le cas non plat. Dans le cas $ n\geq 3 $ et $ h >0 $ partout sur $ \partial \Omega $ et non  constante, on a $ \inf_{\Omega} u_{\epsilon}>0 $ et $ \inf_{\Omega} u_{\epsilon} \to 0 $ si $ \epsilon \to 0 $ (equation de Yamabe), et en particulier en dimension 4, on peut choisir $ f=V_{\epsilon}\in [1/2, 3/2] $ v\'erifiant la deuxieme contrainte ci-dessus, $ ||\nabla V_{\epsilon}||_{\infty}=o(1)C(\epsilon)=o(1){\sqrt \lambda_{\epsilon}}=o(1)\inf_{\Omega} u_{\epsilon} $ (avec ici $\inf_{\Omega} u_{\epsilon} \to 0 $ si $\epsilon \to 0 $). Par exemple $ f=V_{\epsilon} = 1+\epsilon \times C(\epsilon)\times r $ avec $ \epsilon \to 0 $, ici, $ r: x \to d(x,P) $, la distance a un point et $ C(\epsilon) $ la constante de l'article de Hebey-Vaugon.

\smallskip

En dimension 4. Dans le cas de la boule unit\'e (cas plat), on peut prendre $ a(x)\equiv 0 $, $ \tilde h=h\equiv 1 $ et $ f=V_{\epsilon}=1-\epsilon \times C(\epsilon) \times r $ qui est strictement decroissante, on obtient une solution $ v_{\epsilon}\in H^1_0(\Omega) $ reguliere  et non constante (car $ \lambda_{\epsilon} \not =0$, $\lambda_{\epsilon} >0 $) et solution de $ \Delta v_{\epsilon}=\lambda_{\epsilon} V_{\epsilon} (v_{\epsilon}+1)^3 $ avec, $ v_{\epsilon}+ 1 \geq 0 $ et $ u_{\epsilon}=(\lambda_{\epsilon})^{(n-2)/4} (v_{\epsilon}+1)={\sqrt \lambda_{\epsilon}}(v_{\epsilon}+1) \geq 0 $, par le principe du maximum, $ v_{\epsilon} >0 $, et d'apres Gidas-Ni-Nirenberg, $ v_{\epsilon} $ est radiale.

Ici, en dimension 4 avec 2 contraintes, c'est a dire, $ \Delta u_{\epsilon}=V_{\epsilon} u_{\epsilon}^3 $ et $ ||\nabla V_{\epsilon}||_{\infty} =o(1)\inf_{\Omega} u_{\epsilon} $ et $ 0 < \inf_{\Omega} u_{\epsilon}={\sqrt \lambda_{\epsilon}} \to 0 $

{\bf 

Remarques: on n'a pas forc\'ement $ u_{\epsilon}(0) \geq c >0 $, car on a bien $||\nabla u_{\alpha, q}||_2 \geq C{\alpha}^{2/N} $, et, 

$$ \liminf_{q\to N} ||\nabla u_{\alpha,q}||_2 \geq ||\nabla u_{\alpha} ||_2 $$

mais en passant a la limite inf, on n'a pas $||\nabla u_{\alpha}||_2 \geq C{\alpha}^{ 2/N}$. Donc ce n'est pas vrai qu'en passant a la limite inf d'avoir : $||\nabla u_{\alpha}||_2 \geq C {\alpha }^{2/N} $.

a) On a les deux contraintes et $ \inf_{\Omega} u_{\epsilon} \to 0 $. C'est l'exemple voulu avec $ \inf \to 0 $.

b) il n'y a aucune raison que ca blow-up, dans tous les cas, il se peut que l'amplitude soit petite, moyenne ou grande, mais il n'y a pas de raison qu'elle soit grande.
}.

\smallskip

Ce procede se generalise aux dimensions $ n=3 $ et $ n\geq 5 $, on a des exemples avec $ \inf_{\Omega} u_{\epsilon}  \to 0$ .

En prenant une suite $ (\epsilon_k)_k $ avec $ \epsilon_k \to \bar \epsilon >0 $ et $ h>0 $ non constante, on a des exemples avec $ \min_{\Omega} u_{\epsilon_k} \geq m >0 $ (equation de Yamabe et equation du type Yamabe).

\smallskip

En dimension 3, il y a aussi l'exemple de YY.Li et M.Zhu en dimension 3 pour l'equation du type courbure scalaire prescrite sur une variet\'e compacte sans bord de dimension 3 non conf.diffeom. \`a la 3-sphere.

\smallskip

{\it Exemples en dimensions $ n\geq 2 $ concernant la minoration du $ \sup +\inf $ et $ \sup \times \inf $: equation de courbure scalaire prescrite:}

\smallskip

Sur la sph\`ere $ {\mathbb S}_2 $, il y a le resultat de Moser, courbure prescrite $ f $ strictement positive en un point et $ f(-x)=f(x), \forall x \in {\mathbb S}_2 $, voir le livre d'Aubin. on peut prendre $ f=V_i=\sqrt {(\frac{1}{i}+|x \cdot e_1|)} $ reguliere pour chaque $ i $ et $ 0 \leq V_i \leq b <+\infty $, $ V_i \not \equiv 0 $, $ e_1 $ est un vecteur.

\smallskip

Dans le livre d'Aubin, il y a des r\'esultats d'existence d'EDP du type Yamabe $ n\geq 4$ dans le cas compact sans bord et de courbure scalaire prescrite (Cas o\`u  la courbure prescrite $ f $ est telle que $ \sup f \leq k \inf f, k >1$) avec $ n\geq 3 $, dans le cas compact sans bord.

\smallskip

Concernant la minoration du $ \sup \times \inf $, dans le cas plat, on utilise l'exemple de C.C.Chen-C.S.Lin, dans commun.Pure.Appl.Math. 1997, avec ici la courbure prescrite $ V= 1-Kr^{\rho}, n-2 > \rho >0 $ et en dimension $ n\geq 3 $ et $ \sup \times \inf \to +\infty $, en particulier minor\'e. Il y a aussi l'article de C.C.Chen-C.S.Lin de 1999, "Blowing-up...".

\smallskip

{\it Exemples en dimension 4, cas plat et non plat}

\smallskip

On prend une surface de Riemann $ S $ de courbure scalaire $ R=-2 $ (donc sectionnelle $ -1 $) et la sphere $ S_2 $ de courbure scalaire 2 (donc sectionelle 1). Le produit $ S_2 \times S $ est conformement plat et de coubure scalaire nulle, par le theoreme du livre d'Hebey (courbure sectionnelles constantes oppos\'ees). On peut appliquer le Theoreme d'existence de solutions pour des courbures prescrites changeant de signe en dimension 4, sur une variete compacte de courbure scalaire nulle (Resultat d'Escobar-Schoen). On pouvait prendre aussi, le produit de 2 Tores, muni de la metrique produit (courbures sectionelles nulles donc oppos\'ees), $ {\mathbb T}^2 \times {\mathbb T}^2 $ de coubures scalaire nulle. Ici, on ne sait pas minorer les solutions.

Dans le cas non plat, il suffit de considerer la surface $ K3 $ munie d'une metrique non localement conformement plate, de courbure scalaire nulle. On fait la meme chose que dans le cas plat. Ici, on ne sait pas minorer les solutions.

On peut prendre $ P_2({\mathbb C}) $ de dimension reelle 4 non localement conformemnt plat et de courbure scalaire constante positive, on applique le theoreme d'Aubin comme en dimension 5, dans le cas d'invariant de Yamabe positif. Ici aussi, on ne sait pas minorer les solutions.

\smallskip

{\it Exemples en dimension $ n\geq 3 $, cas plat}

\smallskip

Sur le Tore $ {\mathbb T}^n, n \geq 3 $, il existe une metrique plate (voir, Gallot-Hulin-Lafontaine) et donc de courbure scalaire $ R=0 $. Le theoreme dans le livre d'Hebey d'existence de solutions pour des courbures scalaires prescrite changeant de signe et d'integrales negative, sur une vari\'et\'e compacte de courbure scalaire nulle, n'est pas forccement vrai, sauf en dimensions 3 et 4 (Escobar-Schoen) et en dimension 5 par Bismuth avec des conditions supplementairers, dans le cas non conformement plat. 

Par contre dans le cas conformement plat (ce qui nous interesse ici), ce theoreme est vrai en prenant des fonctions cutoff changeant de signe, c'est un theoreme d'Escobar-Schoen.

On pouvait resoudre une equation differentielle sur une boule de l'espace euclidien avec ou sans condition au bord (attention a la formule de Pohozaev). 

Par exemple, dans des couronnes, Kazdan-Warner montrent l'existence de solutions avec conditions de Dirichlet. 

Il y a aussi le resultat de Coron dans des domaines \`a trous, on peut modifier un peu le domaine pour avoir une infinit\'e de solutions pour une infinit\'e de domaines \`a trous, par exemple on enleve des disques. 

On a aussi les solutions de Caffarelli-Gidas-Spruck. Ici, on ne sait pas minorer les solutions.

\smallskip

{\it Exemple en dimension 2 pour l'equation avec singularit\'e au bord: sur la boule unit\'e.}

\smallskip

Ici, $ \Delta=-(\partial_{11}+\partial_{22}) $.

\smallskip

Soit $ \mu \geq 0 $ tel que:

$$ \mu = \inf \{ \int_{\Omega} |\nabla u |^2 dx, \,\, \int_{\Omega} \dfrac{V e^u}{|x-x_0|^{2\alpha}} dx = 1, \,\, u \in H_0^1(\Omega) \} $$

Ici, $ \Omega= B_1(0) $, la boule unit\'e.

\smallskip

1) On prouve que cet ensemble est non vide:

On considere, $ u_{\beta}(r)= \beta(1-r^2) $ et $ 0 < \epsilon \leq V \leq 2 \epsilon $. On prend d'abord la fonction nulle $ u \equiv 0 $. Alors, on choisit $ \epsilon >0 $ petit tel que:

$$ \int_{\Omega} V e^u/ |x-x_0|^{2\alpha} dx \leq 2 \epsilon \int_{\Omega} 1/|x-x_0|^{2\alpha} d x < 1. $$

On utilise $ u_{\beta} $, on utilise les coordonnes polaires:
(Pour chaque $ V $ on associe une famille de fonctions $ u_{\beta} $ donc: $ u_{\beta} $ depend de $ V$: $ u_{\beta}=u_{\beta}(V)=u_{\beta (V)} $, ou bien $ \beta \in {\mathbb R}(V)=[0,+\infty[_V$, $ V $ est une variable cach\'ee: )

$$ \int_{\Omega} Ve^u/|x-x_0|^{2\alpha} dx \geq \int_{B_{1}(0)} V e^{u_{\beta}}/|x-x_0|^{2 \alpha} dx \geq \dfrac{2\pi \epsilon}{2^{2\alpha}} \int_0^{1} r e^{\beta(1-r^2)} dr $$

Donc en prenant $ \beta $ grand, on a:

$$ \int_{\Omega} Ve^u/|x-x_0|^{2\alpha} dx \geq \dfrac{2\pi \epsilon}{2^{2\alpha + 1} \beta} e^{\beta} (1-e^{-\beta}) \to  + \infty.$$

Donc, il exsite $ \beta_0\leq \beta_{max}=\beta_{max}(\epsilon, \alpha, V)$ tel que:

$$ \int_{\Omega} Ve^{u_{\beta_0}}/|x-x_0|^{2\alpha} dx =1. $$

Donc $ \mu $ est bien definit. 

\smallskip

2) En utilisant une suite minimisante et l'injection compacte de Moser-Trudinger dans $ L^1 $. $ \mu $ est atteint et par les multiplicateurs de Lagrange on a l'existence de $ u $ et $ \lambda $ tels que:

$$ \Delta u= \lambda V e^u /|x-x_0|^{2\alpha}, \,\, u \in H_0^1, \,\,\int_{\Omega} Ve^{u}/|x-x_0|^{2\alpha} dx =1. $$

$ u $ est reguliere. Si $ \lambda \leq 0 $ par le principe du maximum $ u \leq 0 $ et donc:

$$ \int_{\Omega} Ve^{u}/|x-x_0|^{2\alpha} dx <1. $$

Ce n'est pas possible.

Donc $ \lambda >0 $ et $ u >0 $. On utilise la premiere fonction propre du Laplacien $ \phi_1 $ (dans le probleme variationel).

On vient de construire une solution du probleme variationel avec singularit\'e au bord et ces solutions( puisque $ \lambda $ est uniform\'ement born\'e) verifient des conditions du probleme du type Brezis-Merle avec bornes uniformes. Ceci dans le cas d'un $ \epsilon >0 $ fix\'e et petit.

Dans le cas o\`u on fait tendre $ \epsilon >0 $ vers $ 0 $. on a

Le terme $\lambda \times \epsilon $ et donc $ \lambda V $ sont unifrom\'ement born\'es et positifs (il suffit d'utiliser la premiere fonction propre du laplacien $ \phi_1 $ pour eliminer le terme $ \int_{\Omega} e^u/|x-x_0|^{2\alpha} \phi_1 dx $) et en utilisant la contrainte $ \int_{\Omega} e^u/|x-x_0|^{2\alpha} dx \geq 1/2\epsilon \to +\infty $ quand $ \epsilon \to 0 $.

Finalement, on a deux exemples, un avec volumes born\'es  (energies born\'ees) et un avec volumes tendant vers l'infini (energies tendant vers l'infini) et $\lambda V  $ uniformement  born\'e.
Il est clair que quand le volume tend vers l'infini (l'energie tend vers l'infini) le $ \sup_{\Omega} u $ tend vers l'infini. 
Il suffit de supposer les volumes born\'es (energies born\'ees).

Pour construire un exemple avec volume born\'es et courbures positives et uniform\'ement born\'ees avec singularit\'e au bord mais avec $ \sup $ divergeant, on part du contre exemple de Brezis-Merle (la singularit\'e est $ x_0=0 $ et $\Omega $ la boule unit\'e centr\'ee en $ (1,0) $):

$$ \Delta u_{\epsilon}=f_{\epsilon}, $$

avec condition de Dirichlet. ($\Delta =-(\partial_{11}+\partial_{22}) $).

On pose alors, $ V_{\epsilon} = |x|^{2\alpha} f_{\epsilon} e^{-u_{\epsilon}} $. Alors on a:

$$ \Delta u_{\epsilon} = V_{\epsilon} \dfrac{e^{u_{\epsilon}}}{|x|^{2\alpha}}, $$

avec condition de Dirichlet. 

On remarque que quand $ x\in B_{\epsilon}(a_{\epsilon}) $, $ |x| $ est de l'ordre de $ d_{\epsilon}=|a_{\epsilon}| $, car $ d_{\epsilon}=\epsilon^{\beta} $ avec $ \beta <1 $. Puis le support de $ f_{\epsilon} $ est dans $ B_{\epsilon}(a_{\epsilon}) $, donc:

$$ ||V_{\epsilon}||_{\infty} \leq C d_{\epsilon}^{2\alpha} \dfrac{1}{{\epsilon}^2} \left ( \dfrac{\epsilon}{d_{\epsilon}}\right)^{2A} $$
Ceci est la premiere condition.

Pour la deuxieme condition on a aussi un facteur $ d_{\epsilon}^{-2\alpha} $, dans $ w_{\epsilon} $ quand $ x\in B_{\epsilon}(a_{\epsilon}) $ et aussi, en distinguant le cas $ \{ x, |x-a_{\epsilon}|\geq d_{\epsilon}/2 \} $, dans ce cas, $ \dfrac{d_{\epsilon}}{|x-a_{\epsilon}|} \leq 2 $ et le cas, $ \{ x, \epsilon \leq |x-a_{\epsilon}|\leq d_{\epsilon}/2 \} $, dans ce cas,
$ |x|\geq d_{\epsilon}/2 $ par l'inegalit\'e triangulaire, car $ |a_{\epsilon}|=d_{\epsilon} $. On obtient alors, la condition de borne uniforme de $ \int_{\Omega} \dfrac{e^{u_{\epsilon}}}{|x|^{2\alpha}} dx $, 

$$ d_{\epsilon}^{-2\alpha}{\epsilon}^2 \left ( \dfrac{d_{\epsilon}}{\epsilon}\right)^{2A}, $$ 

Qui doit etre constant. Ce qui veut dire que l'exposant de $ \epsilon $, dans cette quantit\'e doit etre nul. Et $ d_{\epsilon}=\epsilon^{\beta} $ avec $ \beta =\dfrac{A-1}{A-\alpha} $ convient car $ \alpha \in (0,1) $.

\bigskip

Dans le cas regulier ($ \alpha = 0 $), si on suppose $ 0< a \leq V\leq b $, la premiere condition (sur la fonction nulle) est $ b <1/\pi $. Par les memes arguments que precedemment on a $ \lambda \times a \leq \lambda_1 \leq 6 $, (pour voir que $ \lambda_1 \leq 6 $, on utilise la definition de $\lambda_1 $ comme minimum d'un probleme variationel calculant le quotient pour la fonction explicite $ u=1-r\in H_0^1 $).

Comme $ \int_{\Omega} Ve^u =1 $, on a $ \int_{\Omega} \lambda V e^u = \lambda \leq \dfrac{6}{a} $ et comme $ a $ peut etre choisit arbitrairement inferieur \`a $ b< 1/\pi $, on peut choisit $ a $ tel que 
$\int_{\Omega} \lambda V e^u dx < 24\pi $.
Ainsi, on a exhib\'e dans le cas regulier des suites de volumes born\'es et tel que $\lambda V  $ born\'e et positif et $\int_{\Omega} \lambda V e^u < 24\pi $.

\smallskip

On pouvait prendre les fonctions radiales de la forme $ u_i(r)=\log \dfrac{8c_i^2}{(1+c_i^2r^2)^2}-\log \dfrac{8c_i^2}{(1+c_i^2)^2} $, ($ (c_i) $ born\'ee), nulles au bord et de volume fini  et de masse $ \leq 8\pi $.

De plus, Brezis et Merle donnent un exemple de suites de volumes born\'es et de courbure prescrite positives et born\'ees et de masses egales a $ 4\pi A $ avec $ A >1 $ quelconque (la masse est definit comme $ m =\int_{\Omega} W e^u dx $, avec $ W $ la courbure prescrite et $ u $ la solution). Donc on peut choisir $ A $ de sorte que la masse, $ m $ soit $ m < 24 \pi $, mais $ \sup_{\Omega} u \to +\infty $.

On voit que la condition de borne uniforme sur le volume est vrai dans les deux cas et est necessaire et que pour le deuxieme cas le $ \sup_{\Omega} u $ diverge.

 De meme que precedemment, comme $ b $ peut etre choisit petit, on peut exhiber une suite de volume divergeant dont on connait pas la masse. Par contre en utilisant la suite de fonctions radiales $ (u_i) $ ci-dessus avec $ c_i \to +\infty $ on a exemple avec volume divergeant et masse $\leq 8\pi $.

\smallskip

On regarde maintenant le cas d'une equation avec un operateur different du Laplacien. On ajoute un terme en gradient qu'on ne peut pas eliminer quand on cherche a savoir si l'ensemble des solutions est born\'e.

\bigskip

{\it Exemple en dimension 2 pour l'equation avec operateur different du Laplacien: sur une boule unit\'e.}

On part de :

$$ -\Delta u-\epsilon ((x_1-x_{01})\partial_1 u+(x_2-x_{02})\partial_2 u)= Ve^u, $$

$$ -\Delta u -\epsilon (x-x_0)\cdot \nabla u =Ve^u, $$

Avec condition de Dirichlet. Ici, $ \Delta = \partial_{11} + \partial_{22} $.

On utilise les notations du contre exemple de Brezis et Merle.

Le domaine $\Omega $ est la boule de rayon 1 et de centre $ x_0=(1,0) $. 

On considere $ z_i $ (par la methode variationnelle), tel que:

$$ -\Delta z_i -\epsilon_i (x-x_0) \cdot \nabla z_i =-L_{\epsilon_i}(z_i)=f_{\epsilon_i}. $$

Avec condition de Dirichlet. Par les theroemes de regularite et les injections de Sobolev $ z_i \in C^1(\bar \Omega) $.

On a:

$$ ||f_{\epsilon_i}||_1=4\pi A. $$

Par un theoreme de dualite de Stampacchia ou Brezis-Strauss, on a:

$$||\nabla z_i||_q \leq C_q,\,\, 1\leq q < 2. $$

On resout:

$$ -\Delta w_i =\epsilon_i (x-x_0) \cdot \nabla z_i, $$

Avec condition de Dirichlet.

Par les estimations elliptiques, $ w_i \in C^1(\bar \Omega) $ et $ w_i \in C^0(\bar \Omega) $ uniformement.

Par le principe du maximum:

$$ z_i-w_i\equiv u_i. $$

avec $ u_i $ la fonction du contre exemple de Brezis et Merle.

On ecrit:

$$ -\Delta z_i-\epsilon_i (x-x_0) \cdot \nabla z_i=f_{\epsilon_i}=V_i e^{z_i}. $$

Donc, on a:

$$ \int_{\Omega} e^{z_i} \leq C_1, $$

Et,

$$ 0\leq V_i \leq C_2, $$

Et,

$$ z_i(a_i)\geq u_i(a_i)-C_3 \to +\infty, \,\, a_i\to O. $$

Pour avoir un contre-exemple sur le disk unit\'e, on fait une translation de $ x\to x-x_0 $ dans le contre-exemple pr\'ec\'edent.

\bigskip

Remarques sur ce type d'equations:

\smallskip

1-Comme dans le cas du Laplacien, on peut avoir un exemple avec solutions uniformement bornees, par la methode variationnelle (comme dans le cas avec singularite au bord).

\smallskip

2-On peut ecrire ce probleme dans la boule de rayon 1 ou une ellipse et les deux problemes sont differents. Car si on passe de l'ellipse au cercle par une transformation lineaire $ (y_1,y_2)=(x_1/a,x_2/b) $, le laplacien ne se conserve pas. Si on utilise une tranformation conforme par le theoreme de Riemann, le terme $ x\cdot \nabla u $ ne se conserve pas. Il y a un probleme dans la formule de Pohozaev quand on cherche a savoir si les solutions sont compactes.

\bigskip

{\it Exemple en dimension 2 pour l'equation avec singularit\'e d'exposant positif: sur une boule unit\'e.}

Pour construire un exemple avec volume born\'es et courbures positives et uniform\'ement born\'ees avec singularit\'e au bord d'exposant positif, mais avec $ \sup $ divergeant, on part du contre exemple de Brezis-Merle (la singularit\'e est $ x_0=0 $ et $\Omega $ la boule unit\'e centr\'ee en $ (1,0) $):

$$ \Delta u_{\epsilon}=f_{\epsilon}, $$

avec condition de Dirichlet. ($\Delta =-(\partial_{11}+\partial_{22}) $).

On pose alors, $ V_{\epsilon} = |x|^{-2\beta} f_{\epsilon} e^{-u_{\epsilon}} $, $ \beta \geq 0 $. Alors on a:

$$ \Delta u_{\epsilon} = |x|^{2\beta} V_{\epsilon} e^{u_{\epsilon}} $$

avec condition de Dirichlet. 

On remarque que quand $ x\in B_{\epsilon}(a_{\epsilon}) $, $ |x| $ est de l'ordre de $ d_{\epsilon}=|a_{\epsilon}| $, car $ d_{\epsilon}=\epsilon^{s} $ avec $ s <1 $. Puis le support de $ f_{\epsilon} $ est dans $ B_{\epsilon}(a_{\epsilon}) $, donc:

$$ ||V_{\epsilon}||_{\infty} \leq C d_{\epsilon}^{-2\beta} \dfrac{1}{{\epsilon}^2} \left ( \dfrac{\epsilon}{d_{\epsilon}}\right)^{2A} $$
Ceci est la premiere condition.

Pour la deuxieme condition on a aussi un facteur $ d_{\epsilon}^{2\beta} $, dans $ w_{\epsilon} $ quand $ x\in B_{\epsilon}(a_{\epsilon}) $ et aussi, en distinguant le cas $ \{ x, \epsilon \leq |x-a_{\epsilon}|\leq d_{\epsilon} \} $ du cas $\{ x, |x-a_{\epsilon}|\geq d_{\epsilon} \} $. Pour borner $ \int_{\Omega} |x|^{2\beta} e^{u_{\epsilon}} dx $, il revient de borner:

$$ d_{\epsilon}^{2\beta}{\epsilon}^2 \left ( \dfrac{d_{\epsilon}}{\epsilon}\right)^{2A}, $$ 

Qui doit etre constant. Ce qui veut dire que l'exposant de $ \epsilon $, dans cette quantit\'e doit etre nul. Et $ d_{\epsilon}=\epsilon^{s} $ avec $ s =\dfrac{A-1}{A+\beta} $ convient.

\smallskip

{\it Exemple en dimension 2 pour l'equation avec singularit\'e logarithmique: sur une boule unit\'e.}

\smallskip
La construction de l'exemple est similaire au cas avec singularit\'es negatives et positives. Le point essentiel est de borner la quantit\'e :

$$ -\log (d_{\epsilon}/2d)\times \epsilon^2 \times (d_{\epsilon}/\epsilon)^{2A}. $$

Ceci revient a resoudre :

$$ -\log (d_{\epsilon}/2d) \times d_{\epsilon}^{2A} = \epsilon^{2(A-1)}. $$

Et de verifier qu'il est plus grand que $ \epsilon $. Ici, $ d =diam(\Omega) $. On considere alors la fonction $ f(r)=-\log(r/2d) \times r^{2A} $.

\bigskip

{\it Exemple en dimension 2 pour l'equation avec singularit\'e logarithmique interieure et poids continu: sur une boule unit\'e.}

\smallskip

On reprend les notations du premier exemple en dimension 2 avec $ x_0=0 \in \Omega = B_1(0) $, $ d=dimaetre(B_1(0))=2 $. Pour chaque $ \epsilon >0 $, il existe $ \lambda_{\epsilon} >0 $ tel que:

$$\Delta u_{\epsilon} =\lambda_{\epsilon} V_{\epsilon} \dfrac{1}{-\log \frac{|x|}{2d}} e^{u_{\epsilon}}. $$

$ \Delta = -(\partial_{11}+\partial_{22})$. On choisit $ V_{\epsilon} $ Lipschitzien tel que:

$$0<\epsilon \leq V_{\epsilon} \leq 2\epsilon, \,\, ||\nabla V_{\epsilon} ||_{\infty} \leq k\epsilon, \,\, k >0,$$

Par exemple, $ V_{\epsilon}=\epsilon(1+r^2) $ ou $ V_{\epsilon}=\epsilon(1+x_1^2) $ (pour avoir une solution non-radiale).

En utilisant la premiere valeur propre et la premiere fonction propre et une integration par parties, on a: $ 0 < \lambda_{\epsilon} \times \epsilon \leq C\lambda_1 $, $ C>0 $ independante de $ \epsilon $.

Si $ \lambda_{\epsilon} \times \epsilon \to k_0 >0 $ si $\epsilon \to 0 $, en ecrivant $ W_{\epsilon} =\lambda_{\epsilon} V_{\epsilon} $, alors $ u_{\epsilon} $ serait solution d'un probleme variationnel en dimension 2 avec condition de Dirichlet ralativement a $ 0 < a \leq W_{\epsilon} \leq b <+\infty $ et $||\nabla W_{\epsilon} ||_{\infty} \leq A $, on aurait la compacit\'e de $ u_{\epsilon} $ dans $ C^0 $ ce qui contredit la condition $ 1=\int_{\Omega}\dfrac{1}{-\log \frac{|x|}{2d}} V_{\epsilon} e^{u_{\epsilon}} dx \to 0 $ quand $ \epsilon \to 0 $, car $  0 <\epsilon \leq V_{\epsilon} \leq 2\epsilon $. (On peut utiliser la formule de Pohozaev-Rellich pour borner l'integrale de $ \dfrac{1}{-\log \frac{|x|}{2d}} W_{\epsilon} e^{u_{\epsilon}} $ et aboutir a une contradiction). 

Donc, $ \lambda_{\epsilon} \times \epsilon \to 0 $. On pose alors, $ \bar W_{\epsilon}=\frac{\lambda_{\epsilon} V_{\epsilon}}{\lambda_{\epsilon} \times \epsilon} $. $ v_{\epsilon}=u_{\epsilon}+\log(\lambda_{\epsilon} \times \epsilon) $.

On obtient:

$$ \Delta v_{\epsilon} = \dfrac{1}{-\log \frac{|x|}{2d}} \bar W_{\epsilon} e^{v_{\epsilon}}, $$

$$ 0 < 1=a_1\leq \bar W_{\epsilon} \leq b_1=2 <+\infty $$

$$ \inf_{\Omega}  v_{\epsilon} > -\infty, $$

et,

$$ \inf_{\Omega} v_{\epsilon} = \log (\lambda_{\epsilon} \times \epsilon) \to -\infty,\, \, {\rm avec} \,\, \epsilon \to 0, $$

Pour avoir une infinit\'e de solutions, on prend par exemple, $ V_{\epsilon} = \epsilon (1+sr^2) $ avec $ s \in [0,1] $ ou $ V_{\epsilon}=\epsilon(1+s x_1^2) $, $ s \in [0,1] $. Les fonctions $ V_{\epsilon} $ sont Lipschitziennes, les solutions $ u_{\epsilon} $ sont $ C^2(\bar B_1(0)) $. La fonction $ u_{\epsilon} $ blow-up a l'interieur de $ B_1(0) $ en au moins un point $ x_0 $ (sinon, elle serait born\'ee et convergerait vers $ 0 $ dans $ C^2(\bar B_1(0)) $, par les memes arguments que ci-dessus (formule de Pohozaev-Rellich et le fait que $||\nabla u_{\epsilon}||_q\leq C_q $), ce qui contredirait $ \int_{\Omega} e^{u_{\epsilon}} \to +\infty $). Finalement, on a pour $ R, \beta_0>0 $, $ \max_{B_R(x_0)} v_{\epsilon} \geq \beta_0 >0 $ avec $ B_R(x_0) \subset \subset \Omega $. On a aussi par la formule de Pohozaev-Rellich, $ \int_{\Omega} \frac{e^{v_{\epsilon}}}{-\log \frac{|x|}{2d}} dx \leq C $. Tout cela reste vrai sans la pr\'esence d'un poids (cas r\'egulier).

\smallskip

/////////////////////////////////////////

Exemples et contre-exemples:

\smallskip

On va mieux expliquer en utilisant les contre-exemples de, C.C.Chen-C.S.Lin, Journal of Diff.Geometry.1998:

\smallskip

On considere des metriques conformes (au sens des matrices):

$$ g_u=u^{4/(n-2)} g=u^{4/(n-2)} (dr^2+r^2(1+Weyl\cdot r^2+\ldots) d\theta^2), $$ 

la courbure scalaire de $ g_u $ est $ V $, les contre exemples de Chen-Lin, sont possibles si la platitude est d'ordre $ [(n-2)/2] $, avec un coefficient ne tendant pas vers 0. La courbure scalaire est deduite du tenseur de Riemann, qui est deduit de la metrique en derivant 2 fois. Donc, pour obtenir $ V $, il faut deriver 2 fois $ g_u $, donc il faut deriver 2 fois $ g $, or le developpement de $ g $ contient le tenseur de Weyl, a l'ordre 2, donc, quand on derive 2 fois $ g $, on obtient $ Weyl\not = 0 $. Donc la paltitude est possible jusqu'a l'ordre 2. Cela veut dire que du point de vue de la dimension, $ [(n-2)/2] =2 $, donc, la platitude est possible jusqu'a la dimension $ n=6 $, comme les coefficient devant le tenseur de Weyl, tendent vers $ 0 $, la dimension $ n=6 $ est incluse dans la platitude. cela veut dire que la dimension $ n=6 $ est le cas limite dans la notion de platitude. De plus l'exemple de Chen-Lin, dit que le meilleur qu'on puisse obtenir est $ \sup $ major\'e si $ \inf $ minor\'e, c'est a dire que l'in\'egalit\'e de Harnack explicite n'est pas facilement possible, c'est a dire que c'est difficile de l'obtenir.

\smallskip

Pour Chen-Lin, les contre exemples, commencent des que la platitude est d'ordre $ [(n-2)/2] $ pour $ V $, avec des coefficients ne tendant pas vers 0. Or il faut deriver $ g_u $  2 fois pour obtenir $ V $, mais $ g_u=u^{4/(n-2)} g=u^{4/(n-2)} (dr^2+r^2(1+Weyl\cdot r^2+\ldots)d\theta^2) $ (au sens des matrices), donc, il faut deriver 2 fois $ g $, or, deriver 2 fois $ g $ c'est obtenir $ Weyl\not =0$. Donc, les contre-exemples commencent des que $ [(n-2)/2]=2 $. Comme, pour le cas $ n=6 $, les coefficients devant $ Weyl $ tendant vers 0, les contre exemples, ne commencent que pour $ n\geq 7$. De plus, Chen-Lin, disent que pour la platitude $ [(n-2)/2] $, le meilleur qu'on puisse obtenir est $ \sup $ major\'e si $ \inf $ est minor\'e.(on n'a pas cela). Ce qui veut dire que le cas limite $ n=6$, ce cas est inclus dans les cas ou l'in\'egalit\'e est possible et que le meilleur qu'on puisse obtenir est $ \sup $ major\'e si $ \inf $ est minor\'e.

\smallskip

(La platitude de $ V $ est li\'ee \`a la platitude de $ \partial_{rr} (g) $, considerer la platitude de $ V $, c'est aussi considerer la platitude dans la classe: $ \partial_{rr}(g) $, il faut deriver 2 fois $ g $, c'est "dans la classe de derivation 2 fois". Si on prend $ V=1 $, il faut que cette platitude soit compatible avec la platitude de $ \partial_{rr}(g) $).

\smallskip

Donc, avec la notion de platitude et en particulier $ V=1 $ (equation de Yamabe), le cas limite est $ n=6 $ et le mieux qu'on puisse obtenir est $ \sup $ major\'e si $ \inf $ est minor\'e. On n'a pas forcement $ (\sup)^{\alpha} \times \inf \leq c $. On a (le mieux qu'on puisse obtenir): $ \sup $ major\'e si $ \inf $ est minor\'e pour le cas limite $ n=6 $, une in\'egalit\'e de Harnack implicite.

\smallskip

Donc, avec la notion de platitude, il faut considerer une classe de metriques plus petite. Pour $ n=6 $, en considerant n'importe quelle metrique $ g $ sur $ M $, le mieux qu'on puisse avoir c'est $ \sup $ major\'e si $ \inf $ est minor\'e, une in\'egalit\'e de Harnack implicite, on n'a pas forc\'ement une in\'egalit\'e explicite $ (\sup)^{\alpha} \times \inf \leq c $. Il faut diminuer la quantit\'e de metriques sur la vari\'et\'e $ M $. Par exemple pour les potentiels $ C1 $, en dimension 4, il faut prendre un espace sym\'etrique, Ricci plat. Quand on prend un espace symetrique, on inclut, les espaces euclidien pour lesquels, il y a des contre exemples de Chen-Lin \`a l'ordre: $ (n-2)/2, n\geq 5 $. Donc, il faut faire attention avec les espaces symetriques.

\smallskip

(Dans les contre-exemples de Chen-Lin (1998), il faut que: $ \int_{{\mathbb R}^n} Q(\xi+y) U_0^{2n/(n-2)} dy < 0 $, donc, les contre exemples sont valables, par exemples, qu'avec $ Q < 0 $, par exemple, le blow-up, ne peut pas etre un minimum ($ Q >0$, ca ne peut pas blow-up). D'o\`u, les resulats de Chen-Lin (1997) et L.Zhang (2007), avec $ Q=\Delta V >0 $ positive. Ca peut etre d'ordre $ (n-2)/2 $ et $ Q >0$. Les contre exemples, sont construits pour l'ordre $ \alpha \geq (n-2)/2 $, $ \alpha >1 $ et $ Q <0$).

\smallskip

Donc, avec la notion de platitude on n'a pas l'in\'egalit\'e de Harnack explicite pour $ n=6 $ et les contre-exemples commencent \`a partir de $ n\geq 7$.

\smallskip

Pour obtenir des resultats plus g\'en\'eraux concernant le $ \sup \times \inf $, voir la notion de: g\'eom\'etrisation.

\smallskip

///////////////////////////////////////////////////////////////////////////////////////////////////////////////////////////////////////

\smallskip

Sur l'article de 2012.Journal.Math.Anal.Appl: on a pour toute suite une relation:

\smallskip

$ \sup u_i \leq f(\inf u_i) $, en dimension 4 pour l'Eq. de courbure scalaire prescrite ou du type courbure scalaire prescrite en dimension 4.

\smallskip

-Le d\'enombrable et les ondes: dualit\'e onde-corpuscule. Pour les ondes on a un flot continu d'emissions, alors que pour les corpuscules on a un nombre d\'enombrable d'emissions distinctes. Les deux interpretations de la mati\`ere sont possibles en physique, donc, il suffit d'avoir la propri\'et\'e pour un nombre d\'enombrable ou une suite. Ce qui est le cas en dimension 4, pour l'Eq. de la courbure scalaire prescrite et l'Eq. du type courbure scalaire prescrite avec un potentiel proche d'une constante $ >0 $. On a encore les notions d'enroulement, de torsion (et extra-dimensions) et les valeurs $ (\sup,\inf)$, en d\'enombrable.

\smallskip

-Le continu et minorant-majorant. Dans le cas du continu, emissions continues, on utilise l'estimation a priori, qui, \`a partir du minorant $ m_0 >0 $ donne le majorant  $ M_0=c(m_0)=c(a,b, m_0, K, M) $ et cette fonction est d\'ecroissante de $ m_0 >0 $, donc, quand on prend $ m_0 >0 $ petit et voisin de $ 0 $, $ 0 < m_0 \to 0 $, $ c(m_0) $ croit, $ m_0 \searrow 0 \Rightarrow c(m_0) \nearrow $, on voit alors l'amplitude g\'en\'er\'ee par la variation du minorant $ m_0 >0$, ainsi que l'enroulement et la torsion et les valeurs, $ (majorant, minorant)=(M_0, m_0) $, $ M_0=c(m_0) $, avec l'in\'egalit\'e suivante  $ 0 < m_0\leq M_0=c(m_0) <+\infty $:

$$ 0 < m_0 \leq \inf_M u \leq \sup_K u \leq M_0=c(m_0) <+\infty.$$

\smallskip

Ici, on resout le probleme du cas $ n=4 $ de $ \sup $ major\'e si $ \inf $ est minor\'e de Chen-Lin, voir le point precedent. Dans les contre exemples de Chen-Lin, la regularit\'e est $ C^{\alpha}, \alpha =(n-2)/2, \alpha >1 $, pour ne pas avoir $ \sup $ major\'e si $ \inf $ est minor\'e, qui correspond aux cas $ n\geq 5 $. Il reste le cas de la dimension $ n=4 $. Alors, l'article de Journal.Math.Anal.Appl.2012, repond a cette question: on a bien $ \sup $ major\'e si $ \inf $ minor\'e pour des solutions avec un potentiel $ V $ voisin d'une constante, sans etre une constante forc\'ement, Lipschitzien et petit en norme Lipschitz. On traite alors, le cas de la dimension $ n=4 $ (sur toute vari\'et\'e Riemannienne, non necessairement sans bord, on donne des conditions sur $ V $: Lipschitzien et petit en norme Lipschitz, qui correspond \`a une hypothese plus faible que $ C^{\alpha}, \alpha=(n-2)/2=1, n=4 $, mais la constante de Lipschitz est petite). Dans le print avec potentiel $ C^1 $, on traite aussi, du cas $ n=4 $ avec un potentiel avec amplitude non necessaierment petite, pour un espace localement symetrique Ricci plat.

\smallskip

Comme les contre-exemples de Chen-Lin sont faits avec des suites. Pour mieux expliquer le resultat de la dimension 4, il faut utiliser des suites: On prend $ (u_i,V_i) $ avec $ |\nabla V_i|\leq A_i \to 0 $, alors, pour chaque $ m >0 $, il existe un rang tel que $ A_i \leq \epsilon_0=\epsilon_0(a,b,m), \epsilon_0 >0 $, donc on a la compacit\'e locale \`a partir d'un certain rang: $ \sup_K u_i \leq c $ si $ \inf_M u_i \geq m >0 $, pour $ i\geq i_0 $.

\smallskip

Les contre-exemples de Chen-Lin, c'est pour $ n\geq 5 $. Pour $ n=4 $, dans le print de 2024; "Harnack inequalities for equations of type...", quand le potentiel $ V $ a une amplitude quelconque ($ |\nabla V|\leq A $, $ A $ quelconque), on ne prouve pas: $ \sup $ major\'e si $ \inf $ minor\'e, ce qui veut dire qu'il est peu probable de prouver $ \sup $ major\'e si $ \inf $ minor\'e pour n'importe quelle vari\'et\'e riemannienne $ (M,g) $ et $ n=4 $. On prouve $ \sup $ major\'e si $ \inf $ minor\'e si l'amplitude de $ V $, $ A \to 0 $. Mais si $ A >0 $ est quelonque ceci est peu probable que ca soit vrai, pour n'importe qu'elle vari\'et\'e. D'o\`u le resultat avec potentiel $ C^1 $, sur les espaces localement symetriques Ricci plats. Donc, pour $ n=4$, on n'a pas de contre exemples, et les contre exemples de Chen-Lin, sont vrais pour $ n\geq 5 $, mais, le print de 2024, "Harnack inequalities for equations of type..." peut etre considerer comme un exemple illustrant l'impossibilit\'e de prouver $ \sup $ major\'e si $ \inf $ minor\'e pour $ n=4 $ avec potentiel $ V $ d'amplitude quelconque $ A >0 $ et sur n'importe quelle vari\'et\'e Riemannienne. D'o\`u le print sur les potentiels $ C^1 $, sur les espaces localement symetriques Ricci plats.

\smallskip

//////////////////

\smallskip

Explication pour la dimension $ (n=4) $. On a:

\smallskip

preprint de 2024: "Harnack inequalities for equations of type prescribed scalar curvature": dans la platitude: $ V $ tels que: $ 0 <a \leq V \leq b <+\infty, \,\, |\nabla V|\leq A < +\infty $:

\smallskip

le cas $ (n=4) $ de "Harnack inequalities for equations of type...": $ (P):$ Pour toute vari\'et\'e Riemannienne $ (M,g) $: il y a une relation entre $ \sup $ et $ \inf $: $ \Rightarrow $: $ (Q):$ cela revient \`a considerer, par blow-up, des solutions relativement \`a des potentiels $ \tilde  V_i $ d'amplitude tendant vers $ 0, |\nabla \tilde V_i| \to 0 $.  Donc:

\smallskip

$ non (Q): $ si on considere les solutions relativement \`a des potentiels d'amplitude ne tendant pas vers $ 0 $, $ V > 0, \nabla V \not =0 $ (pas de rescaling ou changement d'echelle, avec $ \lambda \to +\infty $, ou pas de blow-up): $\Rightarrow $: $ non (P) $: $ \Rightarrow $: il existe une vari\'et\'e Riemannienne $ (M_0,g_0) $ pour la quelle le $ \sup $ et l'$ \inf $ ne sont pas li\'es: $ \Rightarrow $: pour cette vari\'et\'e Riemannienne $ (M_0,g_0) $, $ \sup $ major\'e si $\inf $ minor\'e, n'est pas vrai. Comme $ \sup $ major\'e si $ \inf $ minor\'e est vrai dans le cas plat, alors, la vari\'et\'e Riemannienne $ (M_0,g_0) $ est non.loc.conf.plate.

\smallskip

//////////////////

\smallskip

Le preprint de 2024, cit\'e ci-dessus: "Harnack inequalities for equations of type...", permet d'avoir des contre -exemples, en dimension $ n\geq 4 $, en particulier, $ n=4 $, comme on vient de le voir, et $ n\geq 5 $. (potentiel $ V>0, \nabla V\not =0 $, pas de rescaling ou changement d'echelle, avec $ \lambda \to +\infty $ ou pas de blow-up).

\smallskip

(Explication: on a un ensemble de fonctions: $ E=\{V\in C^1, V>0, \nabla V \not =0, \,\, ou, \,\, 0 < c_2\leq |\nabla V|\leq c_1<+\infty \}$, et un groupe de transformations: $ G=\{g: x\to x_0+x/\lambda \} $, agissant sur $ E $, alors on doit avoir $ g\cdot E \subset E $, alors, quand $\lambda \to +\infty $, il n'est plus possible que $ G $ agisse sur $ E $: $ g\cdot E\not \subset E $, homogen\'eit\'e de l'ensemble $ E $ ou invariance de $ E $ ou invariance locale de $ E $. Donc, il n'est plus possible de faire le rescaling ou le changement d'echelle pour $ \lambda \to +\infty $ ou impossibilit\'e de faire le blow-up).

\smallskip

////////////////

\smallskip

On peut fixer chaque vari\'et\'e $ (M,g) $ et avoir des contre-exemples pour cette vari\'et\'e Riemannienne $ (M,g) $, pour $ n\geq 4 $,  avec $ E= \{ V\in C^1, V>0, \nabla V \not =0, \,\, ou, \,\, 0 < c_2\leq |\nabla V|\leq c_1<+\infty \} $: impossibilit\'e de faire le blow-up. 

Ici, comme, pour toute vari\'et\'e Riemannienne $ (M,g), n\geq 4 $, on a des contre-exemples,(on a: $ \sup $ et $ \inf $ ne sont pas li\'es, donc, on n'a pas: $ \sup $ major\'e si $ \inf $ est minor\'e), il faut determiner les vari\'et\'es sur les quelles, par exemple, en dimension $ n=4 $, $ \sup $ major\'e si $ \inf $ est minor\'e est vrai, c'est l'objet du preprint des potentiels $ C^1 $, espaces localement symetriques Ricci plat de dimension 4.(Pour ces espaces aussi, on a des contre-exemples: pour $ V >0, V\in E'=\{ V\in C^1,  V >0, 0 < c_2 \leq |\nabla V|\leq c_1 < +\infty \} $: impossibilit\'e de faire le blow-up, ici aussi).

\smallskip

Cas de la dimension $ n=4 $:

\smallskip

Dans le livre d'Aubin, pour $ n=4 $ et la courbure scalaire $ R=R_g >0 $, il y a des solutions au probleme de la courbure scalaire prescrite, c'est \`a dire en considerant un potentiel $ V >0, V \in C^1 $, avec des conditions sur $ V>0 $. On a, si considere $  E=\{R>0, 0 <c_2 \leq |\nabla R| \leq c_1 <+\infty \} $, on a: $ R=R_g \in E $, en utilisant le preprint de 2024: "Harnack inequalities for equations of type prescribed scalar curvature", on raisonne sur la courbure scalaire, on a des contre exemples: impossibilit\'e d'utiliser le blow-up: $ \Rightarrow \sup $ et $ \inf $ ne sont pas li\'es $ \Rightarrow \sup $ major\'e si $ \inf $ minor\'e n'est pas vrai. Ici, la condition $ 0 < c_2 \leq |\nabla R|\leq c_1 <+\infty $ implique qu'on est sur un espace non localement symetrique, car: $ \nabla R= 0 \Rightarrow \nabla Ricci =0 \Rightarrow \nabla Scalaire = 0 $.

Donc, on a des contre exemples dans les espaces non localement symetriques, d'o\`u, le preprint des potentiels $ C^1 $ sur les espaces localement symetriques Ricci plats.

\smallskip

////////////////////

\smallskip

On peut faire la meme chose pour l'eq. de Yamabe, sur une vari\'et\'e Riemannienne compacte sans bord, avec le preprint de 2024: "Harnack inequalities for equations of type prescribed scalar curvature", on raisonne sur la courbure scalaire, en considerant l'ensemble: $ E=\{R>0, 0 <c_2 \leq |\nabla R| \leq c_1 <+\infty \} $, avec, $ R=R_g \in E $: $ n\geq 4 $ (preprint 2024): impossiblit\'e de faire le blow-up: $ \Rightarrow $ $ \sup $ et $ \inf $ ne sont pas li\'es $ \Rightarrow $ $\sup $ major\'e si $ \inf $ minor\'e, n'est pas vrai $ \Rightarrow \sup \times \inf \to +\infty $. Comme sur une vari\'et\'e compacte sans bord, l'eq. de Yamabe a toujours des solutions, en particulier dans le cas positif, en particulier si $ R >0 $, on a alors des contre-exemples.

Pour avoir des contre exemples, on a impos\'e la condition supplementaire $ 0 <c_2 \leq |\nabla R|\leq c_1<+\infty $, il se peut que sans cette condition on n'ait pas ces contre-exemples.

\smallskip

C'est une r\'eponse \`a la question de YY.Li-L.Zhang.

\smallskip

/////////////////

\smallskip

On peut appliquer cela aussi \`a l'article:"About Brezis-Merle problem with Holderian condition", en considerant l'ensemble $ E= \{ V\in C^1(\bar B_1(0)), b\geq V\geq 0, \nabla V \not =0, \,\, ou, \,\, 0 < c_2\leq |\nabla V|\leq c_1<+\infty \} $ et $ \int_{B_1(0)} e^u dx \leq C $ avec $ bC <16\pi $ ou $ bC<24\pi$, on a des contre-exemples en dimension 2: impossibilit\'e de faire le blow-up. 

\smallskip

/////////////////

On peut appliquer cela aussi \`a l'article: "Harnack inequalities for Yamabe type equations", en considerant l'ensemble, $ R=R_g $ est la courbure scalaire: $ E=\{R >0, 0 < c_2\leq |\nabla R|\leq c_1 <+\infty\} $: impossibilit\'e de faire le blow-up.

D'apres un Theoreme d'Aubin, on a existence de solutions quand: $ \epsilon < \frac{(n-2)R}{4(n-1)} =\frac{(n-2)R_g}{4(n-1)} $. On a aussi, les solutions constantes. Pour l'eq: $ -\Delta u_{\epsilon}+\epsilon u_{\epsilon} = u_{\epsilon}^{N-1}, N=\frac{2n}{n-2}, n\geq 3,  $ sur $ (M,g), n\geq 4 $, $ \Delta =\nabla^j \nabla_j $. Rien ne prouve que les solutions constantes et les solutions d'Aubin coincident. (On a consid\'er\'e l'ensemble $ E $).

Comme on a consid\'er\'e l'ensemble $ E $: on a alors, l'existence d'une sous-suite, quand $ \epsilon_i \to 0 $, de solutions non constantes de $ -\Delta u_{\epsilon_i}+\epsilon_i u_{\epsilon_i} = u_{\epsilon_i}^{N-1}, N=\frac{2n}{n-2}, n\geq 3,  $ sur $ (M,g), n\geq 4 $, compactes sans bord, relatives \`a l'ensemble $ E $, $ \Delta =\nabla^j \nabla_j $, car, on a : $ \sup_M u_{\epsilon_i} \times \inf_M u_{\epsilon_i} \to +\infty $. 

Donc, $ R=R_g >0 $ sur $ (M,g) $ compacte sans bord de dimension $ n\geq 4 $ est une condition n\'ecessaire et suffisante (si, seulement si) pour obtenir le r\'esultat d'unicit\'e. C'est une r\'eponse au probleme 2 de Brezis-Li.

\smallskip

///////////////

\smallskip

On peut appliquer cela, au cas n\'egatif, en considerant l'ensemble: 

$ E=\{V\in C^1, V<0, \nabla V \not =0, \,\, ou, \,\, 0 < c_2\leq |\nabla V|\leq c_1<+\infty \}$: impossibilit\'e de faire le blow-up. Ici, on a la notion de potentiel hyperbolique ou de noyau hyperbolique, sans qu'on ait explicitement ce potentiel ou ce noyau.

(On a: $ \exists (u_i), \,\, u_i >0, \,\, -\Delta u_i+Ru_i=V_i {u_i}^{N-1}, \,\, V_i < 0, \,\, \exists K_0 \subset \subset M, \,\, \sup_{K_0} u_i \to +\infty $. On a la notion suivante: $ u_0>0, \,\,\Delta u_0=\delta_{x_0}, \Delta=\nabla^i(\nabla_i) $ avec $ u_0 >0 $ un objet math\'emtique et de la physique non explicit\'e).

\smallskip

//////////////////////////////////

\smallskip

Concernant le contre-exemple du probleme de Brezis-Merle et "About Brezis-Merle problem with Holderian condition":(on en a parl\'e dans les exemples du preprint "Cas d'existence de solutions d'edp", concernant l'existence de solutions au probleme variationnel en dimension 2, equation du type Liouville, avec potentiel)

\smallskip

On a: on prend $ \Omega=B_1(0) $ et $ 0 < a \leq V_{\epsilon} \leq b <+\infty $ avec $ b|\Omega|=b\pi< 1 $.

\smallskip

Soit $ \mu_{\epsilon} \geq 0 $:

$$ \mu_{\epsilon} =\inf \{||\nabla u||^2, u \in H_0^1(\Omega), \,\, \int_{\Omega} V_{\epsilon} e^u =1 \}. $$

En utilisant la condition $ b|\Omega|<1 $ et l'injection de Moser-Trudinger, on obtient (par l'absurde):

$$ \mu_{\epsilon} \geq \mu_0 >0, \,\, \forall \epsilon >0, $$

En utilisant les multiplicateur d'Euler-Lagrange, on obtient:

$$ \exists \lambda_{\epsilon}, \, \exists u_{\epsilon} \in H_0^1(\Omega), \,\, -\Delta u_{\epsilon}= \lambda_{\epsilon} V_{\epsilon} e^{u_{\epsilon}}, $$

avec condition de Dirichlet au bord.

\smallskip

On utilise la condition $ b|\Omega|<1 $ et  le principe du maximum, on obtient $ \lambda_{\epsilon} >0 $ et $ u_{\epsilon} >0 $.(par l'absurde).

\smallskip

En utilisant le Theoreme 1 de Brezis-Merle, on a : $ \lambda_{\epsilon} \not \to 0 $.(par l'absurde on aurait $ \mu_{\epsilon} \to 0 $, ce n'est pas possible). Donc: $ \lambda_{\epsilon} \geq  \lambda_0 > 0 $. En utilisant la premiere valeur propre du laplacien avec condition de Dirichlet (pour le disque unit\'e $ \lambda_1 \leq 6  $, on prend une fonction test dans le probleme variationnel $ u=1-r $.) et une integration par parties, on obtient:

$$ \lambda_{\epsilon} \times a \leq \lambda_ 1 \leq 6, $$

On pose $ W_{\epsilon} =\lambda_{\epsilon} V_{\epsilon} $:

\smallskip

On prend $ b=1/2\pi $ et $ a=2b/3 $ et $ 0 <c_2 \leq |\nabla V_{\epsilon}|\leq c_1 <+\infty $.(On a alors: $ \int_{\Omega} W_{\epsilon} e^{u_{\epsilon}} dx = \lambda_{\epsilon} \leq 18\pi <24\pi $).(on peut prendre $ b=1/2\pi $ et $ a=4b/5 $ pour avoir $ \int_{\Omega} W_{\epsilon} e^{u_{\epsilon}} dx = \lambda_{\epsilon} \leq 15\pi <16\pi $) 

\smallskip

On obtient: 

$ \int_{\Omega} W_{\epsilon} e^{u_{\epsilon}} dx = \lambda_{\epsilon} \leq 15\pi <16\pi $ ou

$ \int_{\Omega} W_{\epsilon} e^{u_{\epsilon}} dx = \lambda_{\epsilon} \leq 18 \pi <24\pi $ 

et: $ 0 <\tilde c_2 \leq  |\nabla W_{\epsilon}|\leq \tilde c_1 <+\infty $.

\smallskip

On peut alors appliquer ce qu'on a dit sur l'impossibilit\'e de faire le blow-up en utilisant l'article "About Brezis-Merle problem with holderian condition". Donc, pour une sous-suite: 

$$ \sup_{\Omega} u_{\epsilon} \to +\infty, $$

D'apr\'es Brezis-Merle, on a la compacit\'e locale, d'o\`u le blow-up au bord. Ceci est un exemple et un contre-exemple aux problemes 1 et 2 de Brezis-Merle. C'est une r\'eponse aux problemes 1 et 2 de Brezis-Merle.

\smallskip

///////////////////

\smallskip

Explication de l'impossibilit\'e d'utiliser le blow-up:

\smallskip

$ (P) $ est vraie $ \Leftrightarrow $ $ [ Non (P) \Rightarrow $ contradiction $ ] $: est vraie, $ \Rightarrow $: on a utilis\'e le blow-up: or ceci est impossible $ \Rightarrow $ $ [ Non (P) \Rightarrow $ contradiction $ ] $: n'est pas vraie, $ \Rightarrow (P) $ n'est pas vraie.

\smallskip

On considere un enemble de fonctions, comme on l'a fait precedemment, en raisonnant sur la courbure scalaire prescrite ou la courbure scalaire ou les 2 \`a la fois, ou bien sur un ensemble de fonctions, qui apparaissent dans les equations: invariance ou invariance locale des ensembles de fonctions: impossiblit\'e de faire le blow-up.

Pour les exemples et les contre-exemples, il faut determiner, explicitement des solutions, pour pouvoir dire qu'il y a des exemples et des contre exemples.

\smallskip

/////////////////////////////////////////////////////////////////////////////////////////////////////////////

\section{Operateurs d'ordre 2.}

1) In\'galit\'es de Harnack pour les solutions, les sur-solutions, les sous-solutions: ces estimations se font \'a l'aide du proc\'ed\'e d'it\'eration de De Giorgi-Nash-Moser.

\smallskip

2) Existence de solutions: formulation variationelle et th\'eor\`emes du type Lax-Milgram dans les Hilberts.

\smallskip

3) Existence des solutions pour des EDP elliptiques avec operateur Laplacien: Potentiel Newtonien, regualrit\'e en se ramenant ( par soustraction du postentiel Newtonien) aux fonctions harmoniques qu'on sait regulieres par le Theoreme de Weyl (conditions minimales $ L^1 $, mais la formulation variationelle on sait qu'on a mieux). pour la regularite au bord, on suppose la condition de Dirichlet, alors les fonctions harmoniques (on se ramene a ce cas par soustraction du potentiel Newtonnien) sont prolongebles grace \`a la formule de la moyenne qu'on etend \`a travers le demi-espace par symetrisation). La regularit\'e $ C^2 $ se fait grace a la regularit\'e du noyau Newtionnien (se fait par calcul explicite de derivees sous le signe somme), puis par soustraction, on se ramene aux fonctions harmoniques).

\smallskip

4) On dispose des theoremes d'existence et regularite pour le probleme de Dirichlet pour le Laplacien. Pour le cas d'un operateur elliptique d'ordre 2, on utilise la methode de continuit\'e pour prouver l'existence de solutions avec regularite de Schauder. Puis on compare par le principe du maximum dans $ W^{1,2}_0 $ ces solutions regulieres aux solutions obtenues par la formulation variationelle. Si le coefficient devant le terme lineaire $ c\leq 0 $, alors le principe du maximum donne l'unicit\'e des solutions, d'ou la regularit\'e de Schauder, inetreures et sur le bord. Par une alternative de Fredholm on etend au cas ou $ c $ non necessairement negatif ou nul.

\smallskip

5) R\'egualrit\'e $ W^{2,p} $ d'Agmon, elles sont font par la methode des quotients differentiels, elles sont interieures et au bord pour des fonctions tests $ C^2_c(\Omega) $ pour la regularit\'e interieures (au sens des distributions) et pour des fonctions $ C^2_0(\Omega) $ ( $ C^2 $ et nulles au bord) pour la regularit\'e au bord. En plus de la methode des quotients differentiels, Agmon utilise l'estimation a priori de fonctions regulieres dans $ L^p, p>1 $ (approximations) d'Agmon-Douglis-Nirenberg (obtenue par les noyaux).

\smallskip

Les estimations uniformes $ W^{2, p} $ peuvent etre obtenues (en plus de ce que dit Agmon) par le procede de Calderon-Zygmund, qui consiste a obtenir une estim\'ee $ L^p $ pour les derivees secondes a partir de celle du noyau Newtonien. C'est une decompostion en n-cubes de Calderon-Zygmund et l'utilisation du theoreme d'interpolation $ (L^1, L^2)\to L^p (1<p<2)$ (et par dualit\'e $ L^p, p >2 $) de Marcienewicz-Thorin. Puis on ramene l'estimation pour l'operateur $ L $ \`a une estimation pour un operateur \`a coefficient constant comme le Laplacien par exemple.

\smallskip

Les theoremes d'Agmon et de regularite sont utilises dans la construction des fonctions de Green.

\section{Fonctions de Green.}

Le monograph de Druet-Hebey-Robert est assez clair et precis. On donne quelques explications sur les fonctions de Green et sur le fait que le dernier terme dans la decomposition de cette fonction est de classe $ C^1 $. C'est essentielement du au fait que le premier terme est donn\'e explicitement et on utilise une recurrence pour prouver que les termes dela decompostion sont Sobolev et continues en dehors de la digaonale. (Ou bien comme dans le monograph de Druet-Hebey-Robert, on fait la difference des fonctions en deux points et on remarque qu'on a des fonctions Lipschitziennes, Notons comme c'est ecrit dans Ambrosio-Fusco-Pallara,pour pouvoir savoir si une fonction est Sobolev, il suffit qu'on ait des derivees directionelles).

\smallskip

Soit $ (\bar W, g) $ une vari\'et\'e Riemannienne compacte avec ou sans bord. Dans le cas sans bord, on la note $ (M, g) $.

\smallskip

Soit $ H $ la fonction definie comme dans le livre d'Aubin. Cette fonction ne depend pas de la fonction rayon d'injectiv\'e dans le cas d'une vari\'et\'e compacte sans bord. Elle depend (\`a gauche en P) du rayon d'injectiv\'e dans le cas d'une vari\'et\'e \`a bord. On sait que la fonction rayon d'injectivit\'e est continue sur $ W $ avec $ \delta(P)\leq d(P,\partial W) $, on peut alors la prolonger par $ 0 $ sur $ \partial W $. Alors $ \delta(P) $ est continue sur $ \bar W $.

\smallskip

Tout d'abord, remarquons que ces fonctions sont mesurables et integrables, car continues en dehors d'un ensemble de mesure nulle. (un point ou la diagonale).

\smallskip

On pose:

$$ \Gamma_1(P,Q)=-\Delta_Q H(P,Q),\,\, \Gamma_{i+1}(P,Q)=\int_W \Gamma_i(P,R)\Gamma_1(R,Q) dR. $$

Alors,

$$ |\Gamma_1| \leq C [d(P,Q)]^{2-n}. $$

Explicitement et par recurrence les $ \Gamma_i $ sont continues en dehors de la diagonale et pour $ P $ ou $ Q $ proche du bord, elles sont nulles (conditions au bord). On les prolonge par $ 0 $ au bord. Par le theoreme de la convergence dominee de Lebesgue, on prouve qu'elles sont $  W^{1,\infty}  $ en dehors de la diagonale en $ Q $ car elles sont Sobolev est leur derivees bornees. La fonction $ d(P,Q)=d_P(Q) $ est reguliere en $ Q $ en dehors de $ P $ et en fixant $ Q $ elle est Lipschitzienne en $ P $, on alors le fait que $ (R, Q) \to \partial_Q d_R(Q) $ est continue en dehors de la diagonale et $ L^{\infty}_{loc} $. Donc, par le theoreme de convergence dominee (enlever de petites boules) $ \Gamma_i* \partial_Q \Gamma_1 $ et $ \Gamma_i*\partial_Q H $ sont $ C^0 $ en dehors de la diagonale. Donc, les $ \Gamma_i $ et $ \Gamma_i*H $ sont $ C^1 $ en dehors de la diagonale.(Le point essentiel est la fonction distance).(Concernant la fonction distance, $ (R_i, Q) \to d(R_i,Q)=d(R_i, exp_{R_i}(r\theta)), Q=exp_{R_i}(r\theta),  \epsilon_0/2 \leq r\leq \delta(R_i)/3, \theta\in{\mathbb S}_{n-1}$ et la carte de reference est $ [B(R_0,\delta(R_0)/3),  exp_{R_0}] $ avec $ R_i \to R_0 $. On voit alors qu'en coordonnes geodesiques polaires centrees en $ R_i $ un calcul de $ \Delta_{g,Q} d(R_i, Q) \in L^{\infty} $ et $ d(R_i,.) \to d(R_0,.) $ et par les estimations elliptiques on a $ \partial_Q d(R_i, Q) \to \partial_Q d(R_0, Q) $ pour $ d(R_i,Q) \geq \epsilon_0 >0 $.)

Dans le cas o\`u $ W=M $ sans bord:

D'apres sa formule explicite $ H $ est $ C^1 $  en dehors de la diagonale. Elle est $ W^{1,1+\epsilon} $ de chaque variable pour un $ \epsilon >0 $ assez petit, par recurrence et par Fubini, Fubini-Tonnelli:

$$ \partial_P\Gamma_{i+1}(P,Q)=\int_M \partial_P\Gamma_i(P,R)\Gamma_1(R,Q)dR, $$

et,

$$ \partial_Q\Gamma_{i+1}(P,Q)=\int_M \Gamma_i(P,R)\partial_Q\Gamma_1(R,Q) dR.$$

Pour une variet\'e avec bord, on ne peut pas deriver par rapport $ P $, car l'expression de $ \Gamma_i $ et en paticulier de $ \Gamma_1 $ contient $ \delta_P $ qu'on sait seulement continue et pas forcement differentiable. Dans le cas avec bord, les fonctions $ \Gamma_i $ sont continues en $ P $ en dehors de la diagonales et $ L^{1+\epsilon}, \epsilon >0 $ seulement, mais ceci suffit pour prouver la symetrie par le theoreme d'Agmon) 

On a aussi, les estimees de Giraud. Ainsi de proche en proche, les $ \Gamma_i $ deviennent de plus en plus regulieres. on utilise le theoreme de Giraud pour prouver que $ \Gamma_k = r^{\lambda-n}* r^{2-n} $ avec $ \lambda + 2 >n $ (donc $\geq n+1 $ car $ \lambda $ est entier puisqu'on derive les $ \Gamma_i $). ( Si $\lambda >n-1 $ ($ \lambda \geq n $) c'est fini (on regarde $ \Gamma_{k-1} $), si $ \lambda =n-1 $ alors: $ \Gamma_k = r^{-1}* r^{2-n}$ et $ \partial_P \Gamma_k=r^{-2}*r^{2-n} $ et $ \partial_Q \Gamma_k=r^{-1}*r^{1-n} $ qui est en $ \log $ par Giraud et donc $\partial_P \Gamma_{k+1}$ et $ \partial_Q \Gamma_{k+1} $ sont $ C^0(M\times M) $ par Giraud et donc $ \Gamma_{k+1} $ est $ C^1(M\times M)$. (Utiliser le produit de convolution, car $\Gamma_{k+1} $ et sa deriv\'ee Sobolev sont continues pour ecrire qu'elle est $ \int \partial_P $, l'integrale d'une fonction continue donc $ C^1 $, il a convergence uniforme de la convolution, car la fonction et son Sobolev sont continues). Dans le cas avec bord, $ \Gamma_{k} $ continue implique que $ \Gamma_{k+1} $ est $ W^{1,\infty} $ en $ Q $.(On a mieux $ C^1 $ en $ Q $ car la fonction $ (R,Q) \to \partial_Q d(R, Q) $ est $ C^0 $ en dehors de la diagonale de $ W $  et $ H $ est a support compact (elle s'annule avant le bord) et donc ($ \Gamma_k* \partial_Q H $ est $ C^0 $ en dehors de la diagonale de $ \bar W $)).

Dans le cas d'un operateur $ -L=-\Delta+h $, la regualrite de $ \Gamma_{k+1} $ depend de celle de $ h $. Si par exemple $ h $ est $ C^1 $ alors comme precedemment $ \Gamma_{k+1} $ est $ C^1 $ de $ P $ et $ Q $. Bien sur le $ \Gamma_1 $ change:

$$ \Gamma_1(P,Q)=-\Delta_QH(P,Q)+h(Q)H(P,Q). $$

On peut faire plus, on remarque que par exemple que pour les varietes a bord, la parametrix est solution d'une EDP elliptique en dehors de la singularit\'e, on utilise les theoremes d'Agmon et de regularit\'e dans la construction complete de la fonction de Green dans le cas des varietes a bord, voir le monograph de F.Robert (que je n'ai pas cit\'e ici).

\section{Principe du maximum.}

\subsection{La methode "moving-plane" et "moving-sphere".}

\smallskip

La m\'ethode "moving-plane" consiste \`a rechercher, si possible, les points de sym\'etrie pour des E.D.P d\'efinies sur des domaines ayant des axes de sym\`etries, puis, de caract\'eriser ces solutions. On part de "l'infini", un point tr\'es loin, puis, on consid\`ere la fonction et son sym\'etris\'ee par rapport au plan contenant ce point, puis on ram\`ene, le plan jusqu'\'a l'annulation de la diff\'erence entre cette fonction et son sym\'etris\'ee, si c'est le cas, le plan limite est le plan de sym\'etrie. En plus de la determination d'une barriere, cette technique utilise le principe du maximum et le lemme de Hopf. Ici, on suppose les solutions dans $ C^{2, \alpha} $, $ \alpha >0 $.

Cette technique de symetrie a \'et\'e introduite par A.D. Alexandrov.

Notons que certains des auteurs suivants, utilisent les proprietes du au principe du maximum.

\smallskip

Dans les deux points particuliers suivants l'id\'ee essentielle est qu'on doit utiliser les propri\'et\'es du principe du maximum:

\smallskip

1-Le principe du maximum lui meme et le lemme de Hopf.

\smallskip

a) Principe du maximum sur un domaine $ \Omega \subset \subset {\mathbb R}^n, n\geq 2 $ et lemme de Hopf.

\smallskip

b) Principe du maximum pour une vari\'et\'e compacte \`a bord, $ [\lambda, t]\times {\mathbb S}_{n-1}, n\geq 2 $ et lemme de Hopf. L'operateur est global en la variable de la sph\`ere, le principe du maximum s'applique localement dans des domaines de la forme $ ]\lambda,t[\times U $ avec $ U $ un ouvert de carte de la sphere(autour d'un point $(t_0,\theta_0)$ quelconque, on s'est ramen\'e a un ouvert particulier de $ {\mathbb R}^n $, un $n$-cube ou $n$-parallelepipede par exemple pour appliquer le lemme de Hopf).

\smallskip

2-Le minimum est atteint sur le bord, (YY.Li-L-Zhang). Et principe du maximum de Hopf.

\bigskip

{\it Exemple 1: Gidas-Ni-Nirenberg.}

\smallskip

Sur la boule unit\'e $ B $ de $ {\mathbb R}^n $, $ n \geq 3 $, on consid\`ere le probl\`eme suivant:

\begin{displaymath}  \left \{ \begin {split} 
      -\Delta u  & = u^{(n+2)/(n-2)-\epsilon}              \,\, &&\text{dans} \!\!&&B_1(0) \subset {\mathbb R}^n, \\
                  u & >0                 \,\, &&\text{dans} \!\!&&B_1(0) \\
                  u & =  0   \,\,             && \text{sur} \!\!&&\partial B_1(0).               
\end {split}\right.
\end{displaymath}

O\`u $ \epsilon >0 $ est tr\'es petit.

\smallskip

Gidas-Ni-Nirenberg: La solution $ u $ du probl\`eme pr\'ec\'edent est radiale et strictement d\'ecroissante. Il y a des resultats similaires avec potentiel $ V $ radial strictement decroissant.

\bigskip

{\it Exemple 2: De Figueiredo-Lions-Nussbaum et Chen-Li, Ma-Wei.}

\smallskip

Sur un ouvert regulier en dimension 2, on considere:

\begin{displaymath}  \left \{ \begin {split} 
      -\Delta u  &= u^{q},\,\, 1< q < +\infty  \,\, &&\text{dans} \!\!&&\Omega \subset {\mathbb R}^2, \\
                  u & = 0             \,\,             && \text{dans} \!\!&&\partial \Omega,\\                               
\end {split}\right.
\end{displaymath}

Et aussi, 

\begin{displaymath}  \left \{ \begin {split} 
      -\Delta u  &= V e^{u},  \,\, &&\text{dans} \!\!&&\Omega \subset {\mathbb R}^2, \\
                  u & = 0             \,\,             && \text{dans} \!\!&&\partial \Omega,\\                               
\end {split}\right.
\end{displaymath}

avec,

$$ 0\leq V \leq b ,\,\, \text{sur} \,\, \bar \Omega, $$

et,

$$ |\nabla V| \leq A ,\,\, \text{sur} \,\, \bar \Omega. $$

Pour l'equation avec l'exponentielle, on suppose de plus que:

$$ \int_{\Omega} e^u dx \leq C, $$

Alors les solutions sont uniform\'ement born\'ees sur $\bar \Omega $.

Un point essentiel est de borner les solutions au voisinage du bord. Cela se fait par la m\'ethode "moving-plane".
                  
Estimations de De Figueiredo-Lions-Nussbaum: ce sont des estimations uniformes au voisingage du bord, on les prouve grace a la m\'ethode moving-plane. C'est essentiellement une transformationde Kelvin, pour rendre le voisinage du bord convexe, et une application du principe du maximum (suivant la direction normale).

Par la methode moving-plane on prouve que les points crititques sont loin du bord, en se placant  suivant la normale, apres avoir utiliser la transformee de Kelvin. Voir les articles de Chen-Li et Ma-Wei. En particulier, Chen-Li prouvent que si $ \log V $ est uniform\'ement $ C^1 $ et les solutions born\'ees dans $ L^1_{loc} $, alors elles sont uniform.born\'ees au voisinage du bord.

\smallskip

Pour Chen-Li, le potentiel $ V $ doit etre $ C^1 $ et uniformement $ C^1 $.(Ils utilisent la formule de la moyenne et la preuve est locale au voisinage de chaque point du bord pour prouver la compacit\'e au voisinage du bord).

\smallskip

{\bf Remarque:} sur le fait que la r\'egularit\'e du domaine dans (De Figueiredo-Lions-Nussbaum et Chen-Li) est au moins $ C^{2,\beta}, \beta >0 $. Soit $ \alpha_y=<\nu_1|\nu(y)> $, l'angle avec la normale qui doit etre $ >0 $, avec $ y\in \partial \Omega $, qu'on devrait obtenir pour appliquer la m\'ethode moving-plane. Ici, $ \nu(y) $ est la normale au bord en $ y $: Pour definir l'angle, il suffit d'avoir une tangente  a la courbe de $ y \to \nu(y) $, et donc une d\'eriv\'ee de $ \nu(y) $, (l'application de carte $ \phi $ est au moins deux fois d\'erivable). Pour avoir une branche continue de 'cones', il suffit que $ y\to (\nu(y))'$ soit continue, (l'application de carte $ \phi $ est $ C^2 $, $ \phi''$ continue). Pour que l'angle $\alpha_y >0 $ soit  strictement positif, il suffit que la tangente a la courbe de $ \nu(y) $ soit non nulle, il suffit alors d'avoir la regularit\'e $ C^{0,\beta}, \beta >0 $ de la tangente (forme arrondie: continue et arrondie: $ C^{0,\beta} $).(Donc, il suffit que l'application de carte soit $ C^{2,\beta}, \beta >0 $).

\smallskip

Donc: pour De Figueiredo-Lions-Nussbaum et Chen-Li, la regularit\'e du domaine doit etre $ C^{2,\beta}, \beta >0 $.

\smallskip

Donc; pour Chen-Li et Ma-Wei, la regularit\'e du potentiel $ V $ est $ C^1 $. Il doit etre $ C^1 $ et uniform\'ement $ C^1 $.

\smallskip

/////////////////////////////////////////////////////////////////////////////////////////////////////////////////////////////////////////

En utilisant la transformee de Kelvin et le principe du maximum, on classifie les solutions de l'equation :

\smallskip

{\it Exemple 3: Caffarelli-Gidas-Spruck et Chen-Li.}

\smallskip

Caffarelli-Gidas-Spruck. Les solutions de:

$$ \begin{cases} -\Delta u=u^{(n+2)/(n-2)}, \,\, u >0 \,\, {\rm dans } \,\, {\mathbb R}^n \\  

u(0)=1,\,\nabla u(0)=0.\\

 \end{cases}  $$

Sont:

$$   u(x) = (1+\gamma |x|^2)^{(2-n)/2}, \,\,\gamma >0 $$

La preuve de Chen-Li  est plus courte et elle utilise les memes arguments que ceux de Caffarelli-Gidas-Spruck, c'est essentielement une transformation de Kelvin pour avoir un comportement asymptotique  et aussi l'utilisation du principe du minimum pour des fonctions regulieres singulieres en un point puis l'argument moving-plane bas\'e sur le principe du maximum.

Il y a une autre preuve simplifi\'ee du r\'esultat de Caffarelli-Gidas-Spruck obtenue par YY.Li et L.Zhang.

On presente le resultat de Classification par la methode moving-plane des solutions d'une EDP elliptique nonlineaire  en dimension 2:

Chen-Li: pour $ K >0 $, les solutions de:

$$ \begin{cases}

-\Delta u=Ke^u, \,\, {\rm dans } \,\, {\mathbb R}^2 \\
        u(0)=0,\,\nabla u(0)= 0, \\
            \int_{{\mathbb R}^2 } e^u <+ \infty. 
             
 \end{cases}  $$
 
Sont:

$$   u(x) = -2\log(1+\gamma |x|^2), \,\,\gamma= \sqrt {(K/8)} $$

Cette preuve utilise des estimees asymptotiques specifiques a la dimension 2 et le principe du maximum et lemme de Hopf.

\smallskip

{\it Exemple 4: Han.}

\smallskip

Si on remplace $ B $ par un ouvert r\'egulier quelconque, not\'e $ \Omega $, en maintenant la m\^eme \'equation, les solutions ne sont pas forc\'ement radiales.

\begin{displaymath}  \left \{ \begin {split} 
      -\Delta u  & = u^{(n+2)/(n-2)-\epsilon}              \,\, &&\text{dans} \!\!&&\Omega \subset {\mathbb R}^n, \\
                  u & >0                 \,\, &&\text{dans} \!\!&&\Omega \\
                  u & =  0   \,\,             && \text{sur} \!\!&&\partial \Omega.               
\end {split}\right.
\end{displaymath}

O\`u $ \epsilon >0 $ est tr\'es petit.

\smallskip

Han Il existe un voisinage du bord $ \omega $ ne d\'ependant que de la g\'eom\`etrie du domaine $ \Omega $ et de la dimension $ n $, ainsi qu'une constante $ c=c(\Omega,n) >0 $ telle que:

$$ ||u||_{L^{\infty}(\omega)} \leq c. $$

O\`u $ u $ est la solution du probl\`eme pr\'ec\'edent.

\smallskip

{\it Exemple 5: Chen-Lin.}

\smallskip

En utilisant la m\'ethode moving-plane, Chen et Lin ont prouv\'e une estimation a priori sur un ouvert $ \Omega $ de $ {\mathbb R}^n $, $ n \geq 3 $:

Si $ u $ est solution de:

$$ -\Delta u = V u^{(n+2)/(n-2)}, $$

avec,

$$ 0 < a \leq V(x) \leq b < + \infty,\,\, ||\nabla V||_{\infty} \leq A, $$

$$ ||\nabla^{\alpha} V|| \leq C_{\alpha} ||\nabla V||^{\beta(\alpha)}, \,\,\, {\rm avec} \,\,\, \alpha \leq n-2. $$

alors,

$$ \sup_K u \times \inf_{{\Omega}} u \leq c= c(a, b, A, C_{{\alpha}}, n,  K, {\Omega}). $$

La preuve de Chen-Lin consiste a utiliser le principe du maximum et le lemme de Hopf pour prouver que la symetrisee d'une certaine fonction est positive et que la derivee de cette fonction est strictement positive alors que par le raisonnement par l'absurde ("blow-up")  cette fonction a maximum local (point critique), ce qui est absurde.

\smallskip
Sur l'article de C.C.Chen. C.S.Lin. Communications on pure and applied math, 1997. Ils supposent les potentiels $ V $, $ C^1$:

Ceci est du au fait suivant: Ils cherchent a avoir comme donn\'ee, dans la convergence des fonctions blow-up, $ R_i^{n-2} ||v_i-U||_{C^2(B_{2R_i}(0))} \to 0 $. En faisant le rescaling, ou changement d'echelle, $ x=R_i y $, $ y\in B_3(0) $, il suffit et il faut que, la difference, $ w_i(y)=v_i(R_iy)-U(R_iy) $  tend vers $ 0 $ et appliquer les estimations de Schauder, aux fonctions, $ w_i $, $ \Delta (w_i-U) \to 0$, or, $ -\Delta (w_i-U)=R_i^2(V_i(y_i+R_iy)-V_i(y_i))(v_i)^{N-1} +... \to 0 $. C'est a dire qu'il faut que $ R_i^2(V_i(y_i+R_iy)-V_i(y_i)) \to 0 $, or cette difference n'est que la deriv\'ee de $ V_i $ au point blow-up: $ \frac{V_i(y_i+h)-V_i(y_i)}{h} $ avec $ h=\frac{1}{R_i^2} $. Donc, il faut que $ V_i $ soit derivable en $ y_i $ et $ V_i'(y_i)=0 $ et par le theoreme des accroissement fini, $ V_i'(c_h) \to 0 $ avec $ 0 < c_h< h=\frac{1}{R_i^2} $ et la deriv\'ee en $ 0 $ doit etre continue. Donc, $ V_i $ doit etre $ C^1$ en $ y_i $, comme le point $ y_i$ est quelconque, il faut que $ V_i $ et sa limite $ V $ soient $ C^1$. Il faut que pour toute suite $ (V_i) $ convergeant dans $ C^{\alpha} $ vers $ V $, il faut que $ V $ soit $ C^1 $, pour que ceci soit possible, il faut que $ \alpha =1 $, il faut que $ (V_i) $ soient $ C^1 $ et convergent dans $ C^1$ vers $ V $, $ C^1$. Donc, il faut que $ (V_i) $ soient $ C^1 $ et convergent vers $ V \in C^1 $ dans $ C^1$.

\smallskip

a) On a la condition de depart pour les $ V_i $.(Par exemple: $ 1+a_i|x|^{\alpha}\to 1, a_i \to 0 $).

\smallskip

b) On a la condition: \underbar {pour toute suite} $ V_i $ dans $ C^{\alpha}$ qui converge dans $ C^{\alpha}$ vers $ V $, $ V $ doit etre $ C^1$.

Pour que ce soit possible il faut que $\alpha=1$, $ V_i $ soit $ C^1 $ et converge dans $ C^1 $ vers $ V $. On a alors $ V $ est $ C^1 $.

//////////////////////////////////////////////////////////////////////////////////////////////////////////////////////////////////
\smallskip

{\it Exemple 6: Brezis-Li-Shafrir.}

En dimension 2, la methode "moving-sphere" donne pour les solutions de:

$$ -\Delta u= Ve^u, $$

On suppose $  0 < a \leq V(x) \leq b < + \infty $ et que $ V $ est uniform\'ement Lipschitzienne,

$$ \sup_K u + \inf_{\Omega} u \leq c=c(a,b, ||\nabla V||_{\infty}, K, \Omega). $$

////////////////////////////////////////////////////////////////////////////////////////////////////////////////////////////////////

{\it Exemple 7: Korevaar-Mazzeo-Pacard-Schoen: leur article n'est pas correct:}

\smallskip

Sur un ouvert $ \Omega $ de ${\mathbb R}^n $, $ n\geq 3 $, on consid\`ere l'\'equation de Yamabe:

$$ -\Delta u = u^{(n+2)/(n-2)}, \,\,\, u >0 . $$

Il y a l'article de C.C.Chen et C.S.Lin. Chen-Lin prouvent le r\'esultat (l'inegalit\'e de Harnack) pour une courbure prescrite generale et qui peut etre constante, c'est ce qu'on a ecrit avant.

{\bf Sur l'article de Korevaar-Mazzeo-Pacard-Schoen:}

\smallskip

Concerant la proposition avec l'in\'egalit\'e $ u(x)\leq d(x,\Lambda)^{(2-n)/2} (\inf_{\partial B(0,3/4)} u)^{-1} $.

Ce n'est pas correct. Ils ont ecrit n'importe quoi.

\smallskip

1)(/////) Le "blow-up" est mal ecrit: il n'est pas explicit\'e correctement. Ici, ils ont ecrit n'importe quoi. Ils n'expliquent pas le blow-up.

\smallskip

2)(/////) Ils n'expliquent pas le processus du commencement "moving-plane".

la phrase "strating with $ t_1 $ very positive, and continuing as far as possible"  veut dire qu'il existe un $ \lambda >>1$ tres grand, puis on decroit les valeur de $ \lambda $, jusqu'a ce que le processus se termine. Ca va dans le sens decroissant et non croissant.

Mais ce n'est pas ce qu'il faut dire: il faut dire d'abord, \underbar { qu'il existe un rang } $ \nu $ qui enclenche le processus "moving-plane". 

Le cas de la dimension 2 et le cas de la dimension $ n\geq 3 $ ont des preuves differentes. Ils ne disent pas comment faire et ils ne parlent meme pas de rang.

\smallskip

3)(/////) Ils n'expliquent pas le $ t_1= \sup \{\lambda...\} $ de la technique "moving-plane".

Ils ecrivent: $ \inf v(t_i,\theta) < \sup v(t, \theta), t >-2-t_i $, $ t_i=-\log \frac{\lambda^{-1}}{16} $. Premierement ce n'est pas correct, et ici, ils prennent le sup egal \`a $ -1 $, qui dit que le sup est egal a $ -1 $. Puis parlent de theoreme de la moyenne, qui n'a rien a avoir ici avec ce qui est voulu.

\smallskip

4)(/////) Ils ne disent pas s'ils considerent une metrique sur $ {\mathbb R}\times {\mathbb S}_{n-1}$, en fait ils en parlent au debut de l'article avec "cylindrical metric", mais n'exliquent pas comment on utlise le principe du maximum et le lemme de Hopf.

\smallskip

5)(/////) Ils ne disent pas comment appliquer le principe du maximum:

\smallskip

5a)(///) Ils ne disent pas si, il faut considerer une vari\'et\'e Riemannienne $ {\mathbb R}\times {\mathbb S}_{n-1}$ avec le lemme de Hopf o\`u il faut considerer la deriv\'ee normale, donc un produit scalaire, donc une metrique Riemannienne, comme c'est ecrit dans le livre d'Aubin (en ce qui concerne le principe du maximum et le lemme de Hopf sur les vari\'et\'es \`a bord). 

5b)(///) Ils ne disent pas s'il faut un proc\'ed\'e local.

\smallskip

6)(/////) Ils n'expliquent pas les points 39), 40) et 41) et 43) de l'article de Brezis-Li-Shafrir.

Par exemple, pourquoi ils n'expliquent pas le point 39) de Brezis-Li-Shafrir ? de meme pourquoi, Y.Y.Li repete le meme argument dans l'article de 1999, comm.math.phys. ? :

\smallskip

C'est parce qu'ils ont consid\'er\'e la sphere comme une vari\'et\'e constitu\'ee de 2 cartes:

Soit $ u_n $ la fonction definie sur $ [\lambda, 0] \times S^{n-1} $, sur $ S^{n-1}$ on a 2 cartes $ (U_1,\phi_1) $ et $ (U_2,\phi_2) $. La differentiabilit\'e en $ (t_0,\theta_0) $ est definie comme ceci:

\smallskip
 a0) Si $ \theta_0 \in U_1 $ alors, on peut consid\'erer:
 
 $$ \nabla u_n(t_0,\theta_0)= [\partial_t [u_n(t, \phi_1(\theta_{01}, \theta_{02}, \ldots, \theta_{0(n-1)}))], \ldots, \partial_{\theta_j} [u_n(t, \phi_1(\theta_{01},\theta_{02}, \ldots, \theta_{0(n-1)}))]] $$

\smallskip
a1) Si $\theta_0 \in U_2 $ alors, on peut consid\'erer:

$$ \nabla u_n(t_0,\theta_0)= [\partial_t [u_n(t, \phi_2(\theta_{01}, \theta_{02}, \ldots, \theta_{0(n-1)}))], \ldots, \partial_{\theta_j} [u_n(t, \phi_2(\theta_{01},\theta_{02}, \ldots, \theta_{0(n-1)}))]] $$ 

Il y a deux fonctions: $ u_{n1}=u_n(t, \phi_1(\cdot)) $ et $ u_{n2}=u_n(t, \phi_2(\cdot)) $. Selon qu'on se place dans les ouverts de cartes $ U_1 $ ou $ U_2 $. Quand on utilise la definition de la differentiablit\'e en un point $ (t_0,\theta_0) $.

\smallskip

Par exemple en dimension 2, sur le cercle, "par un dessin", on a 2 applications de cartes, "qui n'en est qu'une", c'est une seule fonction: $ \theta \to e^{i\theta}, \theta \in {\mathbb R} $. On a: $ \theta \in ]0, 2\pi[ \to (\cos \theta, \sin \theta) $ et $ \theta \in ]-\pi,\pi[ \to (\cos \theta, \sin \theta) $.
\smallskip

////////////////////

\smallskip

{\bf Remarque importante:} Dans l'article de Brezis-Li-Shafrir, ils parlent de coordonn\'ees polaires, $ (r,\theta) \to (r\cos \theta, r \sin \theta) $. Ce n'est pas cela les coordonn\'ees polaires, car, un systeme de coordonn\'ees est definit, par le fait que tout point, a, un voisinage ouvert hom\'eomrphe \`a un ouvert (connexe, definition d'application: tout point a une unique image, ici un point peut avoir 2 images differentes, ce qui n'est pas possible) de $ {\mathbb R}^2 $, or ici il y a un probleme au point $ (1,0)=e^{i0}=e^{i2\pi} $, il faut sortir du domaine de definition $ [0, 2\pi] $ au voisinage de $ 0 $. Il faut considerer des ouverts (connexes), or en $ (1,0) $ il n'est pas possible de definir un voisinage de $ 0 $. Si on sort de $ [0,2\pi] $, alors il faut considerer 2 cartes comme on l'a dit precedemment, $ ]-\pi,\pi[ $ et $ ]0, 2\pi[ $. Donc, on definit 2 fonctions, alors qu'il faut definir une seule fonction.

Ce qu'ont fait Brezis-Li-Shafrir, c'est un raisonnement avec le revetement : $ p: \theta \to e^{i\theta}, \theta \in [0,2\pi[ $, l'intervalle est semi-ouvert, semi-ferm\'e, il n'est pas ouvert, alors que pour definir un systeme de coordonn\'ees, il faut avoir des ouverts(connexes). Ce n'est pas un raisonnement avec les coordonn\'ees polaires, mais un raisonnement avec le revetement, l'application $ p $ doit etre bijective, pour pouvoir avoir une seule fonction, or pour qu'elle soit bijective, il faut un intervalle semi-ouvert, semi-ferm\'e. Donc, l'intervalle n'est pas ouvert, on ne peut pas parler de coordonn\'ees.

\smallskip

(On a: $\theta \to (\cos \theta, \sin \theta) $ est une carte particuliere du cercle, alors que pour definir un systeme de coordonn\'ees, il ne faut pas que ca depende d'une carte particuliere.).

\smallskip

Brezis-Li-Shafrir, ce n'est pas un systeme de coordonn\'ees, il faut des ouverts et si on considere des ouverts, il faut sortir de l'intervalle, c'est \`a dire considerer 2 fonctions, alors qu'il faut considerer une seule fonction. Brezis-Li-Shafrir, c'est un raisonnement avec le revetement.

\smallskip

(Par exemple: 

1) si on considere 2 cartes, ou 2 fonctions, on n'a plus une fonction du cercle dans $ {\mathbb R} $, car, comme, l'application, la definition d'une fonction, c'est un point a une seule image, or ici, on aurait des points qui auraient 2 images distinctes, ici en dimension 2, on a: $ p(\theta)=e^{i\theta} $ avec 2 inverses, $ \phi_1, \phi_2 $: $ \phi_1(z=e^{i\theta})= \theta $, si $ \theta \in ]0,2\pi[ $ (origine 0) et $ \phi_2(z=e^{i\theta'})=\theta' = \theta-\pi $, si $ \theta' \in ]-\pi,\pi[ $ (origine $ -\pi $). Avec, $ p(\phi_1(e^{i\theta}))=e^{i\theta}, p(\phi_2(e^{i\theta'}))=e^{i\theta'}=-e^{i\theta} $. 

\smallskip

2) si on prend un intervalle ferm\'e en 0 ou $ 2\pi$, ou on \'etend l'intervalle au dela de 0 ou $ 2\pi$, on a dit qu'on n'avait plus d'ouverts connexes, on n'a plus d'applications de cartes, dans ce cas la il n'est pas possible de definir la metrique dans des cartes, par transport, ou pull-back, soit ce n'est pas un ouvert, soit 2 expressions distinctes.

Donc, Brezis-Li-Shafrir, ce n'est pas cela les coordonn\'ees polaires, c'est un raisonnement avec revetement.)

\smallskip

En dimension $ n\geq 3 $, on ne peut pas faire avec la sphere, ce raisonnement avec un revetement, car il y a 2 applications de cartes distinctes, il n'y a pas une application $ p $ semblable \`a : $ p :\theta \to e^{i\theta},\theta \in [0,2\pi[ $ avec $ p $ bijective, pour la sphere $ S^{n-1}$, $ n-1\geq 3-1=2 $, pas de parmetrisation globale en dimension $ n-1\geq 3-1=2 $.

\smallskip

Ils confondent carte particuliere en "coordonn\'ees" spheriques ou hyperspheriques et systeme de coorodonn\'ees polaires. En dimension $ n\geq 3 $, on ne peut pas parametrer globalement la sphere par les "coordonn\'ees" spheriques ou hyperspheriques. On ne peut pas faire le raisonnement du revetement de la dimension 2. Les "coordonn\'ees" spheriques ou hyperspheriques, sont des cartes particulieres du systme de coordon\'ees polaires.

\smallskip

//////////////////

\smallskip

Il y a deux fonctions: $ u_{n1}=u_n(t, \phi_1(\cdot)) $ et $ u_{n2}=u_n(t, \phi_2(\cdot)) $. Selon qu'on se place dans les ouverts de cartes $ U_1 $ ou $ U_2 $, quand on utilise la definition de la differentiablit\'e en un point $ (t_0,\theta_0) $.

\smallskip

C'est pour cela qu'ils(Korevaar-Mazzeo-Pacard-Schoen) considerent la fonction blow-up, et utilisent la converence dans $ C^2$. Alors qu'il n'y a pas besoin de convergence dans $ C^2$.

La convergence dans $ C^2$ utilise la notion de cartes locale pour la vari\'et\'e compacte sans bord qui est la sphere $ S^{n-1}$. Comme la fonction limite est radiale, la deriv\'ee par rapport \`a $ t $ ne dependera plus de la variable angulaire.

\smallskip

Comme pour l'integration (Schoen a utilis\'e l'integration dans son article de 1984, probleme de Yamabe), on a bien une fonction definie sur toute la vari\'et\'e produit $ [\lambda,0] \times S^{n-1} $, en particluier la partie angulaire, mais l'integration est definie de maniere globale mais \`a l'aide cartes locales, on peut alors considrer les fonctions: $ u_{n1}=u_n(t, \phi_1(\cdot)) $ et $ u_{n2}=u_n(t, \phi_2(\cdot)) $, mais, on a une notion globale, l'int\'egrale, qui inclut le local.

\smallskip

Aussi sur le fait d'avoir 2 cartes et donc 2 fonctions $ u_{n1}, u_{n2} $ et aussi l'utilisation de blow-up:

\smallskip

C.C.Chen et C.S.Lin, ont eu le meme probleme, apres utilisation de la fonction blow-up, il y a bien convergence uniforme sur une partie et il reste une autre partie sur laquelle, il n'y a pas de convergence, pour laquelle on ne connait pas le signe de la deriv\'ee. C'est pour cela qu'ils utilisent l'hypothese de l'absurde, sur la partie o\`u il n'y a pas de convergence ni de signe constant de la deriv\'ee. Voir leur article de: Comm.Pure.Appl.Math. 1997.

(Pour C.C.Chen et C.S.Lin. Comm.Pure.Appl.Math.1997. Cette incomprehension, leur a permis de definir une hypothese de l'absurde et de prouver par une autre methode, moving-planes en coordonn\'ees cart\'esiennes, l'in\'egalit\'e optimale $ \sup \times \inf $ dans le cas de la boule unit\'e de $ {\mathbb R}^n $).

\smallskip

Ici: (solution \`a ces problemes et en particulier la notion de coordonn\'ees polaires): il suffisait de dire qu'on fixe la variable angulaire $ \theta $ (variable globale sur la sphere unit\'e, $ S^{n-1} \ni \theta=(\theta_1, \theta_2, \ldots,\theta_n) \in {\mathbb R}^n, |\theta|=1 $, la sphere n'est pas seulement determin\'ee par 2 cartes (c'est ce qu'on a l'habitude de voir par exemple pour le cercle ou la sphere), c'est aussi une partie de $ {\mathbb R}^n$), donc, on fixe la variable angulaire $ \theta $, car on derive par rapport \`a $ t $ et non $\theta $, car l'accroissement des fonctions, ou la methode moving-plane, ou moving-sphere, se fait par rapport \`a $ t $, la variable radiale.

\smallskip

Ici: (en ce qui concerne le blow-up):

\smallskip

On utilise la technique de Schoen avec modification. On met un facteur $ c_i $ dans l'expression des $ L_i= \dfrac{l_i}{c_i^{1/2(n-2)}} [u_i(a_i)]^{2/(n-2)} $, et le majorant, $ \beta_i \to 1 $, de la fonction blow-up, pour avoir $ 0 $ comme point critique de la fonction limite. Il suffit que $ 0 $ soit point critique de la fonction limite et non de chaque terme de la suite. On n'a pas \`a changer de suite de points, comme le font Brezis-Li-Shafrir: $ x= a_i+M_i^{-2/(n-2)} y \to  \tilde x= [a_i+M_i^{-2/(n-2)} y_i]+\tilde M_i^{-2/(n-2)} \tilde y $, dans le cas 1, du blow-up, de l'article de Brezis-Li-Shafrir. Cette id\'ee de Brezis-Li-Shafrir, a \'et\'e reprise par C.C.Chen-C.S.Lin (Comm.Pure.Appl.Math.1997) et L.Zhang (JFA, 2002, Zhang, l'ecrit bien dans son article et explique la technique de Brezis-Li-Shafrir, reprise par Chen-Lin).

\smallskip

Pour ce qui nous concerne: On modifie un peu la technique de Schoen, pour que $ 0 $ soit seulement point critique et maximum de la fonction limite. Il n'est pas necessaire qu'il soit maximum local de tous les termes de la suite, comme dans Brezis-Li-Shafrir, Chen-Lin, L.Zhang. la technique de Brezis-Li-Shafrir, reprise par Chen-Lin, Zhang, avec la technique de Schoen, marche dans le cas d'un ouvert de $ {\mathbb R}^n $, elle n'est plus possible sur les vari\'et\'es ouvertes.

\smallskip

Ceci, cette modification de la technique de Schoen, nous permet de faire le blow-up, dans le cas general, sur les vari\'et\'es. On r\'esout ainsi, le blow-up, sur les vari\'et\'es: les vari\'et\'es ouvertes, non necessairment compacte sans bord.

\bigskip

////////////////////////////////////////////////////////////////////////////////////////////////////////////////////////////////

\smallskip

{\it Exemple 8: C.C.Chen-C.S.Lin.}

\smallskip

En dimension 2, la technique "blow-up" et une in\'egalit\'e g\'eom\'etrique sur les courbures int\'egrales a permis a CC.Chen et CS.Lin de prouver pour l'equation suivante, une ine\'egalit\'e optimale $ \sup +\inf $:

$$ -\Delta u= Ve^u, $$

avec, $  0 < a \leq V(x) \leq b < + \infty $ et que $ V $ est uniform\'ement H\"olderienne ($ C^{\alpha}, \alpha \in ]0,1] $), alors:

$$ \sup_K u + \inf_{\Omega} u \leq c=c(a,b, \alpha,  ||V||_{C^{\alpha}}, K, \Omega). $$

Ce th\'eor\`eme peut etre prouv\'e par la methode "moving-sphere" de YY.Li-L.Zhang (point suivant).
/////////////////////////////////////////////////////////////////////////////////////////////////////////////////////////////////
\smallskip

{\it Exemple 9: Li-Zhang.}

\smallskip

La methode "moving-sphere" donne:

\smallskip

Sur une vari\'et\'e Riemannienne quelconque $ (M,g) $ (non n\'ec\'essairement compacte), de dimension $ n\geq 3 $ et de courbure scalaire $ S_g $, on consid\`ere l'\'equation de Yamabe:

$$ -\Delta_g u + \dfrac{n-2}{4(n-1)} S_g u = u^{(n+2)/(n-2)}, \,\,\, u >0 . $$

Alors: en dimensions $ 3, 4 $ avec en dimension 4 l'operateur conforme $ L_g=-\Delta+\frac{1}{6} S_g $ doit etre coercif dans $ H_0^1 $ (c'est une sorte de condition sur la courbure scalaire), pour tout compact $ K $ de $ M $, on a:

$$ \sup_K u \times \inf_M u \leq c = c(K, M, g, n = 3 \,\, {\rm ou}\,\, n= 4). $$

Apres avoir utiliser la transformation de Kelvin, Li-Zhang prouvent que la positivite de la fonction obtenue par symetrisation implique la positivite de la fonction limite  (obtenue par symetrisation) et dans ce cas elle serait constante, ce qui n'est pas possible par le resultat de Caffarelli-Gidas-Spruck. (La fonction limite est due au "blow-up").

En dimension 4, l'operateur doit etre coercif pour avoir une fonction de Green avec condition de Dirichlet et avec certaines propri\'et\'es: $ C_1 d(x,y)^{2-n} \geq G(x,y) \geq C_2 d(x,y)^{2-n}, C_1, C_2 >0$, voir le print de Frederic Robert sur les fonctions Green pour des operateurs coercifs sur une vari\'et\'e compacte \`a bord.

\smallskip

Dans la methode Y.Y.Li-L.Zhang, il faut supposer l'operateur $ L_g $ coercif, pour $ n=3,4$:

Pour $ n=3 $ c'est du \`a :

[Preuve]: \{ Hypothese de l'absurde: (Etape 1: blow-up analysis)\ldots (Etape n) \ldots (Fin) \}

a) Il y a un probleme par rapport \`a la preuve: Toutes les etapes sont cens\'ees etre vraies, jusqu'a obtenir une contradiction.

\smallskip

b) Ils utilisent comme \underbar {donn\'ee}: comparer le minimum sur la boule ferm\'ee et le minimum sur le bord.(Donc dans ce cas il faut voir si ce n'est pas en contradiction avec le principe du maximum).

\smallskip

c) Si l'inf peut etre atteint a l'interieur, alors (il se peut que l'operateur soit coercif et les solutions seraient constantes et born\'ees), on ne peut pas commencer la preuve: l'etape 1, le blow-up analysis n'est pas possible. "il n'ya rien \`a prouver". Ce qui n'est pas vrai.

\smallskip

c-1) Si l'inf est atteint \`a l'interieur, les solutions sont constantes et ceci bloque l'etape 1: blow-up analsyis.

c-2) Supposons que la preuve est vraie: alors toutes les etapes sont vraies, jusqu'a la fin et on obtient la contradiction. Pour que ce soit vrai: il faut chaque etape soit vraie, en particulier l'etape 1: le blow-up analysis. Or en supposant l'inf atteint \`a l'interieur, (l'operateur peut etre coercif), les solutions sont alors constantes et born\'ees, ce qui fait que l'etape 1: le blow-up analysis n'est pas vrai.

\smallskip

d) Pour $ n=4 $ c'est du au meme fait que $ n=3 $ + \`a la fonction de Green. Pour la fonction de Green, ils utilisent le theoreme 1.36 du livre d'Aubin (il y a la distance geodesique de l'espace ambiant de depart $ M $ et la distance de l'espace, ou ouvert de carte $ d_{\Omega} $, le theoreme, 1.36 du livre d'Aubin dit qu'on peut toujours se ramener \`a cette distance de l'espace ambiant de depart).(De plus, ici, on a besoin de comparer des quantit\'es en fixant 2 points, la carte exponentielle sert quand on compare des quantit\'es en fixant un point, l'autre c'est l'origine. Alors que pour, Y.Y.Li-L.Zhang, ils veulent comparer des quantit\'es en fixant 2 points, pour cela, le theoreme 1.36 du livre d'Aubin est necessaire). On ne peut pas utiliser les valeurs propres des petites boules geodesiques, car, 1) on ne se retrouve plus dans l'espace ambiant de depart qui est $ M $, mais dans des petites boules, alors qu'ils utlisent $ \exp_{x_k}(y) $, donc c'est par rapport \`a l'espace ambiant de depart. 2) En comparant les rayons des petites boules geodesiques, pour les quelles, les valeurs propres sont grandes, il se peut qu'ils soient (ces rayons), plus ou moins grands que les rayons des boules $ B(x_k,2\epsilon_k) $, pour les quels on veut appliquer le principe du maximum. (Il faut comparer les rayons des differentes boules).

\smallskip

e) Donc, pour $ n=3,4$ l'operateur $L_g$ doit verifier le principe du maximum, donc doit etre coercif. Si on considere un autre operateur elliptique pas forc\'ement l'operateur conforme, le nouvel operateur doit verifier le principe du maximum, donc, doit etre aussi coercif, pour la dimension 3, 4, c'est a dire dans la preuve de YY.Li-L.Zhang. l'operateur doit etre coercif.
\smallskip

{\bf Remarques: sur la dimension 3:}

\smallskip

Cas plat: dans l'article "Majorations du type..." de 2004: le potentiel $ V $ est Lipschitzien et le bord du domaine born\'e est quelconque (le bord n'est pas n\'ecessairement regulier).

\smallskip

Alors que dans la methode de C.C.Chen-C.S.Lin: le potentiel doit etre $ C^1 $ et uniformement $ C^1$. On l'a dit deja.

\smallskip

Alors que dans la methode de Li-Zhang: le bord du domaine doit etre regulier. C'est du \`a: 

\smallskip

1) In\'egalit\'e de Poincar\'e: le domaine doit etre au moins Lipschitzien. Quand on considere un probleme elliptique ou le principe du maximum.

\smallskip

2) Le principe du maximum de Hopf: les solutions doivent etre $ C^2(\Omega)\cap C^1(\bar \Omega) $, pour avoir la regularit\'e des solutions jusqu'au bord, il faut qu'il y ait une notion de bord regulier.

\smallskip

3) L'operateur doit etre coercif: on l'a deja dit, mais ceci inclus le point 1) et 2).

Dans la methode Li-Zhang, dans le cas plat aussi: il y a une contrainte: il faut que le bord soit regulier. Operateur coercif. Il faut que le principe du maximum soit possible.(D'o\`u la citation de l'article de Berestycki-Varadhan-Nirenberg).

Dans la methode de Li-Zhang, dans le cas plat aussi: ils ont cherch\'e a comparer le minimum sur le bord et \`a l'interieur (des boules):

a) il se peut que le minimum soit atteint \`a l'interieur des boules, et aussi: 

b) il se peut que le minimum soit un minimum global (d'o\`u par exemple, la necessit\'e du principe du maximum sur tout l'ouvert $ \Omega $, ou espace ambiant).

c) le principe du maximum doit etre vrai localement partout, d'o\`u il faut une condition globale, il doit etre vrai sur l'espace ambiant.(et apres on l'applique localement. On considere des conditions sur les fonctions du type: $ C^2(\Omega') \cap C^1(\bar \Omega') , \forall \Omega'\subset \subset \Omega $: alors soit on a une infinit\'e de conditions, soit on considere une condition simple qui est $ C^2(\Omega) \cap C^1(\bar \Omega)$, qui englobe tout, et puisque on a pris $ C^1(\bar \Omega)$, il faut que le bord soit regulier, car une fonction reguliere, bord compris, implique, que le bord est regulier).

Tout cela necessite un bord regulier pour l'espce ambiant $\Omega $.

\smallskip

Cas non plat: Concernant l'article de Math.Aeterna. 2011. Tout d'abord, c'est le cas non plat.

L'operateur n'est pas necessairement coercif et le bord n'est pas necessairment regulier.

Difference avec l'article de 2004: "Majorations du type..."

Dans le cas non plat, il y a moins de possibilit\'es. Dans le cas plat, il y a plus de compacts, les ferm\'es born\'es  sont compacts. Alors que sur une vari\'et\'e Riemannienne, il y a moins de possibilit\'es de caracteriser les compacts: D\'enombrabilit\'e. Possibilit\'es. Aussi, potentiel Lipschitz pour une metrique Riemannienne n'implique pas potentiel Lipschitz pour la distance de $ {\mathbb R}^3 $. Dans l'article de Math.Aeterna 2011, en dimension 3, le potentiel est Lipschitz pour la distance geodesique, (les solutions) dependent de la metrique Riemannienne, aussi pour le cas o\`u la metrique Riemannienne est euclidienne. Alors que pour le cas d'un ouvert de ${\mathbb R}^3$, dans l'article de 2004, le potentiel est Lipschitz pour la distance usuelle de $ {\mathbb R}^3 $, (les solutions) dependent de la structure d'espace metrique pour la distance usuelle. Ce ne sont pas les memes ensembles de solutions.

\smallskip

///////////////////////////////////////////////////////////////////////////////////////////////////////////////////////////////////

\smallskip

La preuve de Li-Zhang s'adapte en dimension 3 sur un ouvert $ \Omega $ de bord regulier, ici, on a pris $ \Omega =B_1(0)$, la boule unit\'e de $ {\mathbb R}^3 $, quand on suppose les courbures prescrites uniform\'ement entre deux constantes positives $ a, b $ et $ s$-h\"olderiennes, $ s\in ]1/2,1] $, de constante de H\"older $ A $,  on a une in\'egalit\'e du type:

$$ (\sup_K u)^{2s-1} \times \inf_{\Omega} u \leq c= c(a,b, A, s, K, \Omega) $$. 

///////////////////////////////////////////////////////////////////////////////////////////////////////////////////////////////////

\subsection{Principe du maximum dans $ W_0^{1,1}(\Omega) $}

Le principe du maximum dans  $ W_0^{1,1}(\Omega) $ est bas\'e sur l'inegalite de Kato dont on peut trouver la preuve dans le livre d'Otared Kavian.

\smallskip

Kato: si $ u \in L^{1}(\Omega) $ et $ \Delta u \in L^1(\Omega) $, alors on a au sens des distributions:

$$ \Delta u^+ \geq \chi_{ \{ u\geq 0 \}} \Delta u. $$

On remarquera qu'en utilisant la fonction $ f_{\epsilon}(t)={\sqrt {(t^2+\epsilon^2)}} $, on prouve que $ Trace (|u|)=|Trace(u)| $ et $ Trace (u^+)=(Trace (u))^+ $ et donc si $ u \in W_0^{1,1} $ alors $ u^+\in W_0^{1,1} $ et appliquer a $ u^+ $ ce qui est consider\'e comme le principe de comparaison dans le monograph de Brezis-Marcus-Ponce. (il suffit de consider la fonction $ -\Delta \xi=1 $ avec condition de Dirichlet).

\smallskip

{\it Exemple 1: Brezis-Merle.}

On consid\`ere deux suites de fonctions $ (u_i, V_i) $ solutions (au sens des distributions) de:

$$ -\Delta u_i =V_i e^{u_i}, $$

avec,

$$ V_i \geq 0, \,\, ||V_i||_{L^{p}(\Omega)} \leq C_1, \,\,  ||e^{u_i}||_{L^{p'}(\Omega)} \leq C_2, $$

avec, $ 1 < p \leq + \infty $. Alors, on a, ou bien

$$ \forall \,K \subset \subset \Omega, \,\, \sup_K |u_i|  \leq c=c(K, C_1, C_2, \Omega) $$

ou bien,

$$ u_i \to -\infty, \,\, {\rm sur \,\, tout \,\,  compact \,\,  de} \,\, {\Omega}, $$

ou bien,

Il existe un ensemble fini de points $ S $ tel que $ u_i \to -\infty $ uniform\'ement sur tout compact de ${\Omega}-S $ et $ u_i \to +\infty $ sur $ S $ et, au sens faible, on a,

$$ V_ie^{u_i} \to \sum_j\alpha_j \delta_{a_j}, \,\, {\rm avec } \,\, \alpha_j \geq  4\pi/p', \,\, {\rm et} \,\, a_j \in  S. $$ 

Comme corollaire de ce r\'esultat, on a:

$$  \sup_K u  \leq c=c(\inf_{\Omega} V, ||V||_{L^{\infty}(\Omega)}, \inf_{\Omega} u, K, \Omega). $$

Cela veut dire que le maximum local est fonction du minimum.

Concernant la compacit\'e globale:

On a aussi la compacit\'e globale de Chen-Li, cit\'e dans la section 4.1. Le resultat de Chen-Li utilise le fait qu'on a compacit\'e au voisinage du bord lorsque $ ||\nabla \log V||\leq A $, puis il etend ce resultat lorsque $ ||\nabla V||\leq A_1 $, pour cela il utilise l'extension des resultats de Brezis-Merle (Theorem 1 de Brezis Merle), qui reste vrai dans des domaines Lipschitzien d\`es qu'on a la regularit\'e des solutions dans $ W_0^{1,2} $, car le principe du maximum est valable dans ce cas (on utilise l'integration par parties, qui est vrai des que la regularit\'e du bord est Lipschitzien. Pour l'inegalit\'e de Sobolev aussi et la resolution d'un probleme variationnel dans $ L^2 $. Lipschitz suffit). Ceci pour la fonction $ u_2 $ du debut de preuve de Chen-Li. (et aussi l'extension des fonctions harmoniques et la formule de Poisson, pour $ u_1 $ qui necessite une application conforme). Voir la preuve du corollaire de l'article de Chen-Li.

\smallskip

(Voir l'article de Sweers-Nazarov. Journal.Diff.Equations. 2007. Pour les conditions de regularit\'e $ W^{2,p} $ des problemes sur des ouverts Lipschitziens.)

\smallskip

Remarque sur la preuve de Chen-Li : pour etendre la partie $ u_1 $ harmonique, il faut supposer le domaine analytique, pour pouvoir utiliser une transformation conforme qui reste invariante par le Laplacien. Principe de symetrisation de Schwarz. Donc, le resultat de compacit\'e reste vrai avec la regularit\'e smooth lorsque on suppose  $ ||\nabla \log V||\leq A $. Mais la regularit\'e du domaine doit etre suposs\'ee analytique lorsqu'on passe a $ ||\nabla V||\leq A_1 $.

Dans leur preuve Chen-Li, utilisent le fait que l'operateur est invariant par appplication de carte, ceci est possible si cette application est conforme, elle preserve le Laplacien. Puis, symetrise la fonction en symetrisant un probleme de Dirichlet, puis soustrayent les valeurs aux bords et ils obtiennent l'image de $ v_1 $. Alors $ u_1 $ est l'image de $ v_1 $ par l'application de carte. Maintenant pour construire $ v_1 $ ils utilisent une symetrisation d'un probleme de Dirichlet, qui requiert les solutions dans $ W^{2,p}\cap C^2(B_{\epsilon}) \cap C^1(\bar B_{\epsilon}), p >2 $ (la formule de representation de Green reste valable, dans ce cas, voir la preuve dans Gilbarg-Trudinger). Puis, ils utilisent la formule integrale de Poisson (qui necessite d'avoir l'operateur Laplacien).

\smallskip

a) Pour utiliser la formule de Poisson, on conserve le Laplacien : transformation conforme $\phi $.
\smallskip

b) Ils symetrisent $ uo\phi $ ils obtiennet une fonction $ u_v \in C^1(\bar B_{\epsilon}(0)) \cap W^{2,p} $.
\smallskip

c) Ils resolvent : $ -\Delta v_1=-\Delta u_v $ avec condition de Dirichlet sur $ B_{\epsilon}(0) $.
\smallskip

d) Ils utilisent la formule integrale de Poisson pour $ v_1-u_v \in W^{2,p}\cap C^2(B_{\epsilon}) \cap C^1(\bar B_{\epsilon}), p >2 $.
Sur le bord, il n'y a que les valaurs de $ u $.  

(Ce travail revient \`a symetriser une fonction harmonique qui n\'ecessite le theoreme de symetrisation de Schwarz
, qui necessite une application conforme, donc un domaine de depart $ \Omega $ analytique).

\bigskip

{\it Exemple 2: Equation avec poids continu:}

\smallskip

Sur un domaine born\'e $ \Omega \subset \subset {\mathbb R}^2 $ avec $ d=diametre(\Omega) $ et $ 0 \in \Omega $, on consid\`ere des fonctions $ (u, V) $ solutions (au sens des distributions) de l'equation du type Liouville:

$$ -\Delta u =\frac{1}{-\log \frac{|x|}{2d}} V e^{u} \,\, {\rm dans }\,\, \Omega $$

avec,

$$ u\in L^{\infty}_{loc}(\Omega), \,\, 0< a \leq V \leq b< +\infty, \,\, 0 \in \Omega. $$

Alors on a l'in\'egalit\'e de Harnack implicite:

\smallskip

Pour tout compact $ K $ de $ \Omega $:

$$ \sup_K u  \leq c(a, b, K, \Omega, \inf_{\Omega} u) $$

et,

$$ \inf_{\Omega} u\geq m >-\infty \Rightarrow  \sup_K u \leq c=c(a, b, K, \Omega, m). $$

\bigskip

\section{Concernant les in\'egalit\'es du type Harnack}

\bigskip

Fixant un minorant $ m >0 $, on a une relation entre $ \sup_K v $ et $ m $ pour toute $ v >0 $ relativement a un $ W $ avec des conditions a priori sur $ W $, $ v $ et $ W $ sont li\'es par une equation ou plusieurs equations. Or, en prenant le minorant $ m=\inf_M u >0 $, on a l'inegalit\'e pour $ v>0 $ or $ u>0 $ est d\'eja une solution du probleme, donc $ \sup_K u $ est fonction de $ \inf_M u $ et des autres param\`etres.

\bigskip

Pour l'equation de la courbure scalaire prescrite en dimensions 3 et 4, et equation du type courbure scalaire prescrite en dimension 3:

\bigskip

a) En dimension 3, on a la dependance explicite du $ \sup $ en fonction du $ \inf $. Ici l'operateur n'est pas necessairement  coercif, il est suppos\'e non-coercif.

$$ (\sup_K u)^{1/3} \times \inf_M u \leq c, $$

\bigskip

b) En dimension 4, on a une estimation a priori. En plus on a des exemples pour ce cas. Dans le cas ou la courbure scalaire prescrite est constante et l'operateur conforme est non coercif, on a une in\'egalit\'e de Harnack implicite. Cette in\'egalit\'e de Harnack implicite est differente de celle qui est explicite, ce ne sont pas les memes fonctions: $ c: m\to c(m) $.

\smallskip

En general dans le cas d'une vari\'et\'e Riemannienne $ M $ de dimension 4, on a une estimation a priori. On peut se ramener a une estimation du $ \sup $ en fonction du $ \inf $ en ecrivant:

$$ \sup_K u \leq c(a, b, \inf_M u, K, M), $$

pour toute $ u >0 $ solution de l'equation de la courbure prescrite en dimension 4 sur $ M $, relativement a $ V $ (Lipschitzienne) verifiant:

$$ ||\nabla V||_{\infty} \leq \dfrac{3a {\sqrt a}}{32e^2{\sqrt 2}}\inf_M u, $$

et,

$$ 0 < a \leq V \leq b < +\infty, $$

\smallskip

En effet, $ u >0 $ et $ V $ sont li\'es par l'equation, on peut supposer que le gradient de $ V $ et $ \inf_M u $ sont aussi li\'es, ici on a deux equations liant $ u>0 $ et $ V $, une contrainte de plus. On a des exemples quand $ ||\nabla V_i||_{\infty} \leq A_i \to 0 $ et $ \inf_M u_i \geq m >0 $, on a une estimation locale uniforme.

\smallskip

Ici, la constante $ k(a)=\dfrac{32e^2{\sqrt 2}}{3a {\sqrt a}}$ entre $ \min_M u $ et la norme du gradient est plus grande. Il faut faire la preuve pour $ \min_M u \geq m >0 $, car on doit eliminer des termes en $ o(1) $.

\footnote{ Pour la preuve, dans le cas non plat, comme c'est une estimation locale, on prend $ x_0 \in M $, on fait un changement de metrique conforme $ \tilde g = e^f g=\phi^3 g $ de sorte  que $ \tilde Ricci_{x_0} = 0 $, on a (voir le livre d'Aubin) $ f(x_0)=0 $ et $ \phi(x_0)=1 $. On prouve l'estimation en supposant $ \min_M u \geq m >0 $ pour la nouvelle fonction $ v=u /\phi $ pour $ m>0$ (car, voir l'article de la preuve du cas general en dimension 4, il y a des termes $ \tilde Ricci_{y_i} $ et $ \tilde R_{\tilde g}(y_i) $, des termes en $ o(1) $), comme $ \phi(x_0)=1 $, $ \min_{B_r^{\tilde g}(x_0)} v \geq (\min_M u)/2 \geq m/2 >0 $. Comme on fait un changement de metrique conforme, au voisinage de $ x_0 $, tout est de "l'ordre" de $ u $ et $ m $ car $\phi(x_0)=1 $, $ f(x_0)=1 $, l'exponentielle (en norme de matrices) est $ (1+\epsilon)-$ Lipschitzienne, on a chang\'e de metrique, les majorants restent de l'ordre de $ 1+\epsilon $ au voisinage de $ x_0 $ car $ \phi(x_0)=1 $, ($ f(x_0)=0$) et la constante $ k(a)$ entre $ \min_M u $ et la norme du gradient de $ V $ est plus grande. 

\smallskip

On prouve alors, en supposant l'inegalit\'e entre le $ \min_M u $ et le gardient de $ V $: $ \min_M u \geq m >0  \Rightarrow \sup_{B_r^{\tilde g}(x_0)} v \leq c \Rightarrow \sup_K u \leq c(a, b, m, K, M) $. Apres on remplace $ m $ par $ \min_M u $.}

\smallskip

{\bf Remarque:} i) Dans le cas d'un ouvert $\Omega $ de $ {\mathbb R}^4 $, la deuxieme contrainte est "v\'erifi\'ee" (apr\'es blow-up), dans l'article de C.C.Chen et C.S.Lin (Comm.on pure.appl.math.1997), si on suppose $ V \in C^2 $ uniform\'ement. Ils obtiennent une inegalit\'e entre $ \nabla V(y_i) $ et $ [u_i(y_i)]^{-1} $, dans ce cas ils obtiennent l'in\'egalit\'e optimale $ \sup \times  \inf $. (les points $ y_i $ sont les points blow-up).

ii) Dans le cas plat on s'int\'eresse particulierment au cas $ V $ Lipschitzien et non $ C^1 $ pour obtenir l'inegalit\'e optimale $ \sup \times \inf $ avec la contrainte supplementaire: $ ||\nabla V||_{\infty} \leq \dfrac{3a \sqrt a}{8e^2 \sqrt 2} \inf_{\Omega} u $.  On obtient alors un r\'esultat du type Brezis-Li-Shafrir, en dimension 4 avec deux contraintes. Notons que dans l'article de C.C.Chen-C.S.Lin (comm.on pure.appl.math. 1997) prouvent que $ \sup \times \inf \to +\infty $ avec $ V=1-Kr $ ($\rho =1 $, $ V $ Lispchitzien et non $ C^1 $), on voit alors que la deuxieme contrainte est suffisante pour obtenir l'inegalit\'e optimale $ \sup \times \inf $, on en d\'eduit la majoration uniforme locale des normes $ L^4 $ et $ H^1 $ ($ N=\frac{2n}{n-2}=4 $). Dans le cas non plat, on obtient un r\'esultat du type Brezis-Merle, en dimension 4 avec deux contraintes (Une in\'egalit\'e de Harnack implicite).

\bigskip

iii) Dans le cas non plat, un exemple, ce resultat s'applique lorsque $ V\equiv constante >0 $. On retrouve, en particulier, le r\'esultat de YY.Li et L. Zhang, eux, obtiennent l'in\'egalit\'e optimale $ \sup \times \inf $ dans le cas positif et ici, en particulier on a une in\'egalit\'e de Harnack implicite dans le cas ou l'operateur conforme n'est pas forcement coercif. 

\bigskip

iv) Dans le cas non plat, un autre cas, on l'a dit ci-dessus, c'est celui o\`u $ ||\nabla V_i ||_{\infty} \leq A_i \to 0 $ et $\min_M u_i \geq m >0 $, on obtient une estimation locale uniforme.

\bigskip

v) Dans le cas non plat, on a un r\'esultat analogue a celui de Brezis-Merle en dimension 2, en dimension 4: pour deux suites $ (u_i), (V_i) $ solutions de l'equation de la courbure scalaire prescrite en dimension 4, $ 0 < a \leq V_i \leq b $ et $||\nabla V_i||_{\infty} \leq A_i \to 0 $, alors $ u_i\geq m >0 \Rightarrow \sup_K u_i \leq c(a, b, (A_i)_i, m, K, M) $. En particulier, si on prend $ m=\inf_M u_i >0 $, on a l'in\'egalit\'e de Harnack implicite, $ \sup_K u_i \leq c(a,b, (A_i), \inf_M u_i, K, M) $ (La fonction $ m\to c $ est decroissante de $ m >0 $ et change en une autre fonction $ \tilde c $ si on remplace les suites $ (u_i,V_i) $ par une d'autres suites $ (v_i,W_i) $ avec des hypoth\`eses similaires, et on obtient: $ \sup_K v_i \leq \tilde c(\tilde a, \tilde b, (\tilde A_i), \inf_M v_i, K, M) $).

\smallskip

vi) Voir la section 1 pour des exemples de $ (u_i,V_i) $ avec les deux contraintes dans le cas plat et non plat, avec $ \inf $ ne tendant pas vers 0 forc\'ement ou n'est pas forc\'ement minor\'e par $ m >0 $ ou tendant vers 0 ou minor\'e par $ m >0 $.

\footnote{Concernant l'exponentielle dans le preuve du r\'esultat ci-dessus. On commence par changer la metrique $ d_g \to d_{\tilde g} $, on a un facteur $ (1+\epsilon) $, puis on considere $ \tilde \exp_y $, l'exponentielle pour la metrique $ \tilde g $, construite a partir d'une carte $ (\Omega, \psi) $ normale en $ x_0 $ pour la metrique $ \tilde g $, on passe de $ d_{\tilde g} $ \`a $ ||\psi o \tilde \exp ||_{{\mathbb R}^4}$, on a un majorant $ (1+\epsilon) $, puis on utlise le fait que $ \dfrac{d}{dz} [z\to \tilde \exp_{x_0} (z)]_{|z=0} =id_{{\mathbb R}^4} $, c'est \`a dire que c'est vrai pour la differentielle de $ \psi o \tilde \exp $ et par continuit\'e des differentielle, (en $ y $ et $ z $ assez petits), on a un majorant $ (1+\epsilon) $, puis on utilise la norme matricielle $ \sup $ des matrices et on ecrit que: $ \psi o \tilde \exp_y(z_1)-\psi o \tilde \exp_y(z_2) =\int_0^1 \partial_s [\psi o \tilde \exp_y(s(z_1-z_2)+z_2)] ds=\int_0^1 [(\partial_z(\psi o \tilde \exp_y)(s(z_1-z_2)+z_2))\cdot (z_1-z_2)] ds $ et on obtient; $ ||\psi o\tilde \exp_y(z_1)-\psi o \tilde \exp_y(z_2)||_{{\mathbb R}^4} \leq (1+\epsilon) ||z_1-z_2||_{{\mathbb R^4}} $, uniform\'ement en $ y $ voisin de $ x_0 $ et $ |z_1|\leq r, |z_2|\leq r, r >0 $ assez petit.}  

\smallskip

La formulation pr\'ec\'edente en dimension 4 s'applique lorsqu'on suppose en $ x_0 $ point critique de $ V>0 $ $ V\in C^1(M) $. on fait la preuve avec $ \min_M u \geq m >0 $ pour $ v=u/\phi$ avec $ \phi>0 $ duc changement de metrique conforme telle que $ \tilde Ricci_{\tilde g}(x_0)=0 $, on obtient:

$$ \exists \,c,R >0,\,\, \sup_{B_R(x_0)} v \leq c, $$

Donc,

$$ \exists \, c,R >0,\,\, u(x_0) \leq \sup_{B_R(x_0)} u \leq c,\, \forall\,  u>0 $$

Donc, on obtient une in\'egalit\'e du type:

$$ u(x_0)\leq c(x_0,V, \inf_M u, M, g), $$

pour toute solution $ u >0 $ de l'equation de la courbure scalaire prescrite en dimension 4 relativement \`a un $ V >0 $, $ V\in C^1(M) $ fix\'e ayant le point $ x_0 $ comme point critique, $ \nabla V(x_0)=0 $. Ou bien $ V \in C^{1,\alpha}(M) $, $ 0 < \alpha \leq 1 $, $ 0 < a \leq V \leq b < +\infty $ et $ ||V||_{C^{1,\alpha}}\leq A $ et $ \nabla V(x_0)=0 $, on obtient:

$$ \forall \,\, u >0,\,\, u(x_0)\leq c(a, b, A, \alpha,x_0,\inf_M u, M, g). $$

Il y a des exemples de solutions de l'equation de la courbure prescrite en dimension 4 avec $ x_0 $ point o\`u $ V $ est maximum et $ M $ compacte sans bord dans le livre d'Aubin.

Quand on prend les notations de YY.Li et L.Zhang: $(M,g)=(B_1(0)\subset {\mathbb R}^4, g) $ et $ 0 < V\in C^1(B_1(0)) $  avec $ \nabla V(0)=0 $, on obtient l'estimation a priori dans un voisinage de $ 0 $ en supposant le minimum unifrom\'ement minor\'e et aussi:

\be u(0) \leq c(V, B_1(0), g, \inf_{B_1(0)} u), \ee

pour toute solution $ u >0 $ de l'equation de la courbure scalaire prescrite en dimension 4 relativement a $ V $ avec les conditions precedentes ou $ V\in C^{1,\alpha}(B_1(0))$ avec $ 0 < \alpha \leq 1 $, $ 0 < a \leq V \leq b < +\infty $ et $ ||V||_{C^{1,\alpha}}\leq A $ et $ \nabla V(0)=0 $.

\be \forall \,\, u >0,\,\, u(0)\leq c(a, b, A, \alpha,\inf_{B_1(0)} u, B_1(0), g).\ee

Le resultat precedent est semblable au resultat de G.Tarantello en dimension 2 avec singularit\'e interieure. Ici, on peut prendre par exemple, $ V(x)=1+\epsilon |x|^{1+\beta} W(x) $ avec $ \beta >0 $, $ 0 <\epsilon \to 0 $ et $ W \in C^1 $, $ 0\in M $, $ (M,g)=(\Omega=B_1(0),g) $. Ou bien $ V(Q)=1+\epsilon [d(Q,P)]^{1+\beta} W(Q) $ avec $  \beta >0 $, $ 0<\epsilon \to 0 $, $ W \in C^1(M) $, $ P\in M $. Ici on obtient une in\'egalit\'e du type:

$$ u(0) \leq c(\epsilon,\beta, W, g, \inf_{B_1(0)} u). $$

C'est un r\'esultat semblable a celui de G.Tarantello en dimension 2, elle considere $ -\Delta u = |x|^{\beta} W(x) e^u, \beta >0 $ et obtient le resultat suivant: $ u(0)+\inf_{\Omega} u \leq c $.

Ici, en dimension 4, on considere: $ V(x)=1+\epsilon |x|^{1+\beta} W(x) $ avec $ \beta > 0 $, $ 0 <\epsilon \to 0 $ et $ W \in C^1 $, $ 0\in M $, $ (M,g)=(\Omega=B_1(0),g) $.

\smallskip

Pour l'equation de Yamabe en dimensions 5 et 6:

\smallskip

c) En dimension 5, il y a une relation explicite entre $ \sup_K u $ et $ \inf_M u $.

$$ (\sup_K u)^{1/7} \times \inf_M u \leq c, $$

d) En dimension 6, on a une estimation a priori et une relation implicite entre $ \sup_K u $ et $ \inf_M u $ (in\'egalit\'e de Harnack implicite). Pour toute solution $ u >0 $ de l'equation de Yamabe en dimension 6:

$$ \sup_K u \leq c(K, M, g, \inf_M u). $$

Dans le cas d'une vari\'et\'e compacte sans bord $(M,g) $ de dimension 5 ou 6 et d'invariant de Yamabe $\mu_g >0 $, le probleme ne se pose pas, car Druet, YY.Li et L.Zhang et F.C.Marques prouvent que les solutions sont uniform\'ement born\'ees, major\'ees et minor\'ees uniform\'ement par des constantes positives.

\smallskip

Pour l'equation de Yamabe et du type Yamabe en dimension 4, on a l'in\'egalit\'e:

$$ (\sup_K u)^{1/3} \times \inf_M u \leq c, $$

ici, il n'y a pas de condition sur la courbure scalaire, on considere le cas ou les operateurs, conforme $ L_g=-\Delta+\frac{1}{6}S_g $ ou $ -\Delta +h $ ne sont pas necessairement coercifs. Un exemple du cas non coercif, dans le probleme de Yamabe: le cas negatif et le cas nul sur une vari\'et\'e compacte sans bord, si $ S_g\equiv -1 $ ou $ S_g \equiv 0 $, l'operateur conforme est non coercif (la fonction $ u\equiv 1 $ rend la fonctionnelle de Yamabe negative ou nulle). Dans le cas n\'egatif et sur une vari\'et\'e compacte sans bord, on resout l'equation de la courbure scalaire prescrite avec une courbure prescrite cutoff $ -\bar M \leq f \leq 0 $, $ \bar M >0 $ avec $ f\equiv 0 $ dans une boule $ B_r(P) $ petite (cas plat et non plat), (ceci est fait dans le livre d'Aubin), puis par un th\'eoreme du livre d'Aubin, il existe un voisinage de  $ f $, $ V_f $ sur lequel on peut resoudre le probleme de la coubure scalaire prescrite, on se ramene \`a la boule $ B_{r/2}(P) $ au probleme de Yamabe avec $ f+c(\bar M)=c(\bar M) >0 $ et $ u >0 $ la solution, avec le fait qu'on se sait pas si l'operateur conforme $ L_g=-\Delta+\frac{1}{6} S_g $ est coercif dans $ H_0^1(B_{r/2}(P)) $.(Si $c(\bar M) \to 0 $  quand $ \bar M \to +\infty $, on remplace $ u $ par $ v=u \cdot \sqrt {c(\bar M)} $, en dimension 4 par exemple, comme on a fait dans certains exemples). Ici, $ v >0 $ solution de l'equation de Yamabe en dimension 4 avec le fait qu'on sait pas si l'operateur conforme est coercif.

\smallskip

Concernant l'article: "Some uniform estimates for scalar curvature type equations". En dimension 4, on a en plus de l'estimation a priori, une dependance implicite du $ \sup_K u $ en fonction du $ \inf_{\Omega} u $ (in\'egalit\'e de Harnack implicite):

$$ \sup_K u \leq c(a, b, A, \alpha, K, \Omega,\inf_{\Omega} u) ,$$

pour le theoreme 3

\bigskip

et,

$$ \sup_K u \leq c(a, V, K, \Omega,\inf_{\Omega}u). $$

pour le theoreme 4.

\bigskip

{\bf Remarque:} On a bien une fonction $ c $ de l'inf car, pour chaque $ m $ on a une constante  $ c_0 $ majorant de $ \sup_K u $ et l'ensemble des majorants est lui meme minor\'e par $ \sup_K u $ et $ m $, $ c $ est alors l'inf des majorants de $ \sup_K u $, qui est unique. On voit alors que cette correspondance $ m  \to c $ est bien une fonction et aussi la correspondance $ (a, b, K, \Omega, m) \to c $ est une fonction. Apres on remplace $ m $ par $ \inf_{\Omega} u $ comme on l'a dit ci-dessus. On a bien la correspondance $ (a, b, K, \Omega, \inf_{\Omega}u) \to c $ est une fonction. 

\smallskip

De plus la fonction $ m\to c $ est d\'ecroissante en $ m >0 $. On a alors, si $\inf_{\Omega} u \geq m >0 $ (ou $\inf_M u \geq m >0 $, si on considere une vari\'et\'e Riemannienne $ M $), alors $ c (a, b, K, \Omega, \inf_{\Omega} u) \leq c(a, b, K, \Omega, m) <+\infty $, ou $ ( c(a, b, K, M, \inf_M u) \leq c(a,b, K, M, m) $, si on est sur une vari\'et\'e Riemannienne $ M $), on obtient alors l'estimation a priori lorsque l'inf est uniformement minor\'e. Comme toute fonction derivable est continue sauf sur un ensemble denombrable de points, on a des ensembles ou la fonction $ m\to c $ est continue.

\smallskip

//////////////////////////

\smallskip

Une consequence de l'in\'egalit\'e: $ (\sup_K u)^{1/3} \times \inf_M u \leq c $. On regarde par exemple $ n=4 $: On a:

$$ \inf_{B_{2R}(0)} u \geq c_1 \int_{B_R(0)} u^3 dx = c_1 \int_0^R (\int_{\partial B_r(0)} u^3 d\sigma_r) dr, $$

$$ \sup_{B_R(0)} u \times \inf_{B_{2R}(0)} u \geq c_2 \int_{B_R(0)} u^4 dx, $$

donc,

$$ c\geq (\sup_{B_R(0)} u)^{1/3} \times \inf_{B_{2R}(0)} u \geq c_3 (\int_{B_R(0)} u^4 dx)^{1/3} \times (\inf_{B_{2R}(0)} u )^{2/3} \geq $$

$$ \geq c_4 (\int_{B_R(0)} u^4 dx) ^{1/3} \times (\int_0^R (\int_{\partial B_r(0)} u^3 d\sigma_r) dr )^{2/3}, $$

On obtient:

$$ (\int_{B_R(0)} u^4 dx) ^{1/3} \times \left (\int_0^R (\int_{\partial B_r(0)} u^3 d\sigma_r) dr \right )^{2/3} \leq c_5. $$

On a:

$$ [Vol(u^2\cdot g, B_R(0))]^{1/3} \times \left [\int_0^R (Aire(u^2\cdot g, \partial {B_r(0)})) dr \right ]^{2/3} \leq c_6.$$

On a:

$$ [Vol_{(u^2\cdot g)}(B_R(0))]^{1/3} \times \left [\int_0^R (Aire_{(u^2\cdot g)}( \partial {B_r(0)})) dr \right ]^{2/3} \leq c_6.$$

C'est une in\'egalit\'e isoperimetrique entre volume conforme et aire conforme pour la metrique conforme $ u^2\cdot g $. Le terme de l'aire est une moyenne des aires conforme et l'autre terme et le volume conforme.

\smallskip

Quand on a l'in\'galit\'e optimale $ \sup \times \inf $ ca implique la majoration uniforme du volume local. 

\smallskip

Ici, $ n=4 $, on a une relation entre le volume conforme local et l'aire conforme locale, c'est une in\'egalit\'e isop\'erimetrique pour la metrique conforme $ u^2\cdot g $.

\smallskip

Il est possible  d'avoir des resultats de ce type pour $ n=3,5 $.

\smallskip

Lorsqu'on a l'in\'egalit\'e optimale et l'operateur conforme est coercif et $ V=1$, Y.Y.Li-L.Zhang, bornent le volume local dans ce cas.

\smallskip

Ici, on suppose qu'on considere l'eq. suivante ($ R_g =-6 $):

$$ \Delta u-u=u^3, u >0, $$

avec $ \Delta =-\nabla^i \nabla_i $.

\smallskip

l'operateur conforme n'est pas necessairement coercif. Donc, on n'est pas dans le cas Li-Zhang, on n'a pas l'in\'egalit\'e optimale. Si de plus on suppose que le laplacien est coercif, petites boules geodesiques, alors on peut appliquer ce qu'on a dit precedemment, l'in\'egalit\'e isoperimetrique ci-dessus. c'est un cas, ou n'a pas l'in\'egalit\'e optimale, o\`u on peut appliquer ce qu'on a dit avant.

\bigskip

///////////////////////////////////////////////////////////////////////////////////////////////////////////////////////////////

\bigskip

\bigskip

\bigskip

{\bf Remarques importantes:}

\smallskip

On a vu dans l'article. Bull.Sci.Math.2006, (voir aussi, le print: Quelques remarques sur les vari\'et\'es, fonctions de Green et formule de Stokes): qu'on peut avoir, selon le point de vue une in\'egalit\'e de Harnack \`a 3 parametres ou 2 parametres. 

Si on suppose la condition de Brezis-Merle ou Brezis-Li-Shafrir: $ \inf_M u >0 $, on a, \`a partir de l'in\'egalit\'e \`a 3 parametres:

$$ \sup_M u\times \inf_K u \geq c(\sup_K u) >0, $$
 avec $ m \to c(m) >0 $ une fonction croissante de $ m >0 $, $ \sup_K u\geq \inf_K u \geq \inf_M u >0 $, (positivit\'e stricte):

$$ \sup_M u\times \inf_K u \geq c(\inf_K u) >0, $$

On a bien une in\'egalit\'e de Harnack \`a 2 parametres entre $ \sup_M u $ et $ \inf_K u $, en admetttant la condition de Brezis-Merle ou Brezis-Li-Shafrir: $\inf_M u > 0 $. Qui est possible en considerant le probleme sur $ (M, g) $ avec $ u >0 $ sur $ M $, puis on prend un sous-ensemble ouvert connexe, $ \tilde M $, alors on aura : $ \inf_{\tilde M} u >0 $. (Par exemple autour d'un point, comme le font Li-Zhang).

\smallskip

Comme le probleme se pose autour de chaque point $ x_0 \in M $, comme le font Li-Zhang. On prenant un ouvert relativement compact de $ M $, $ \tilde M \subset \subset M $. On a toujours $ \inf_{\tilde M} u >0 $. On a:

$$ \forall K \subset \subset \tilde M, \,\,\sup_{\tilde M} u\times \inf_K u \geq c(\inf_K u) >0, $$

1) Ici, comme l'operateur est coercif, par le principe du maximum, on a toujours: pour tout compact $ K $ de $ M $: $ \inf_K u >0 $. Il n'y a pas besoin de supposer la condition, de Brezis-Merle ou de Brezis-Li-Shafrir.

Par le principe du maximum, on a toujours: $ \inf_K u >0 $ (aussi, ici, on a directement $ \inf_K u >0 $, car $ u >0 $ est strict.positive sur $ M $ et continue sur $ \bar M $ et $ K $ est compact). Donc, on a l'inegalit\'e de Harnack sans la condition de Brezis-Merle ou de Brezis-Li-Shafrir.

\smallskip

2) Ce qui parait ambigu, est que la fonction $ m\to c(m)$ est croissante de $ m >0 $, on ne connait pas sa forme explicite. Par exemple, on ne peut pas avoir une fonction linaire, $ c(m) =c_0 m $, car l'exemple qu'on a donn\'e avant: $ x\to (\mu/\mu^2+|x|^2)^{(n-2)/2}, \mu \to +\infty  $, implique qu'on a : $ \sup_M u \geq c_0 \frac{\sup_K u}{\inf_K u} \to c_0 >0$, alors que $ \sup_M u \to 0, \mu \to +\infty $ (et aussi le cas $ \mu \to 0 $, on ne peut aps avoir $ c(m)=c_0 m^k, k >1 $, une puissance). Donc, parfois, on ne peut pas avoir $ c(m)$ lineaire. Ce qui veut dire aussi, que l'in\'egalit\'e qu'on a obtenue, est une vraie in\'egalit\'e de Harnack. Par exemple, on n'a pas une in\'egalit\'e triviale du type $ \sup_M u \geq c_0 >0 $.

\smallskip

3) on a:

$$ \forall K \subset \subset  M, \,\,\sup_M u\times \inf_K u \geq c(\inf_K u) >0, $$

ou bien:

$$ \forall K \subset \subset M, \,\,\sup_M u\geq \frac{c(\inf_K u)}{\inf_K u} >0, $$ 

avec, $ m\to c(m) >0 $ fonction croissante de $ m >0$.

Cette in\'egalit\'e \`a 2 parametres se deduit du th 2 de l'article du Bull.Sci.Math.2006. en prenant $ x_0 \in K $ et $ m=\inf_K u >0 $. Ce qui etait ambigu: il ne faut pas prouver l'in\'egalit\'e en prenant $ v >0$, $ \inf_K v >0$, mais en considerant un point $ x_0$ et on la prouve en prenant $ v(x_0) \geq m >0$. Si on suppose $ \inf_K v =m >0$, il n'y a rien a prouver. Mais en prenant $ v(x_0) \geq m >0$, on prouve l'in\'egalit\'e en prenant $ m=\inf_K u >0 $, puis on l'applique \`a $ v=u$, $ v $ particuliere, $ v=u$, comme $ x_0 \in K $, on a bien $ v(x_0)=u(x_0) \geq \inf_K u >0 $.

\smallskip

Le th 2 de l'article du Bull.Sci.Math.2006. ecrit, permet de deduire l'in\'egalit\'e \`a 2 parametres, en prenant, $ x_0 \in K $,$ m=\inf_K u >0 $ et $ u_i=u $.

\smallskip

Comme on l'a dit, \underbar{ la fonction $ c(m)$ ne peut pas etre triviale}, en prenant, les fonctions particlieres ci-dessus. $ x\to (\mu/\mu^2+|x|^2)^{(n-2)/2}, \mu \to +\infty  $ ou $ \mu \to 0 $.

\smallskip

//////////////////////////////////////////////////////////
\smallskip

Dans Einstein-Lichnerowicz, la seule maniere de combiner, l'eq. de la courbure scalaire prescrite, eq. de Schrodinger non lin\'eaire, et la relativit\'e g\'en\'erale, est de considerer l'eq. de Yamabe. Ici, \`a la limite, il est possible de de recuprer l'eq. de Yamabe. Pour ce qui nous concerne, on a l'eq. de la courbure scalaire prescrite.

\smallskip

Dans Einstein-Lichnerowicz, il y a eu separation de la variable temporelle et spaciale, puis il y a eu utilisation de metriques conformes. On a: \`a un instant $ t =0 $, par approximation, on obtient l'eq. de la courbure scalaire prescrite: car on part de $ -\Delta u+ (\frac{n-2}{4(n-1)}R_g-|\nabla \Psi|^2)u= \frac{1}{2} m \Psi^2 u^{N-1}, u >0, N=\frac{2n}{n-2}, n\geq 3 $, comme on consider\'e $ \Psi $ voisine d'une constante (impulsion \`a l'instant $ t=0 $, puis le potentiel se stabilise autour d'une constante), on peut supposer $ |\nabla \Psi | \approx 0 $, et on obtient l'eq. de la courbure scalaire prescrite, avec un potentiel $ V $ voisin d'une constante $ >0 $.(On en a parl\'e pour la dimension 4, d'o\`u les conditions de platitudes). L'eq, devient:$ -\Delta u+ \frac{n-2}{4(n-1)}R_g u= \frac{1}{2} m \Psi^2 u^{N-1}, u >0, n\geq 3, N=\frac{2n}{n-2} $, eq. de la courbure scalaire prescrite. On en a parl\'e pour la dimension 4. Ceci reste valable en toute dimension $ n\geq 3 $: $ \Psi $ voisine d'une $ constante >0  \Rightarrow |\nabla \Psi | \approx 0 \Rightarrow $ l'eq. de la courbure scalaire prescrite $ \Rightarrow $ eq. de l'invariance conforme $ \Rightarrow $ astronomie et g\'eom\'etrisation de l'astronomie.

\smallskip

On peut dire que l'eq. de la courbure scalaire prescrite est aussi l'eq. de la supersymetrie, comme l'eq. de Yamabe. L'eq. de Yamabe est un cas particulier de l'eq. de la courbure scalaire prescrite.

\smallskip

On peut dire qu'avec ce proc\'ed\'e d'approximation: que l'eq. de la courbure scalaire prescrite est une eq.de l'astronomie. Au d\'epart, on s'est focalis\'e  sur l'eq. de Yamabe (qui est obtenue aussi par approximation), on peut considrer de la meme maniere , l'eq. de la courbure scalaire prescrite. C'est une equation plus g\'en\'erale que l'eq. de Yamabe. C'est l'eq. de l'invariance conforme. (Laplacien conforme, courbure scalaire, tenseur de Weyl, eq. de la courbure scalaire prescrite: invariants conformes). Masse et Masse positive: courbure, scalaire, de Ricci, de Weyl.

\smallskip

Il y a le formalisme de Penrose: g\'eom\'etrie conforme+compactification+courbure+masse.

\smallskip

Si on veut etre plus precis et expliquer mieux: on commence par poser la question suivante: Trouver $ (u,V), u >0, V>0 $ solutions de:

$$ -\Delta u+\frac{(n-2)}{4(n-1)} R_g u= Vu^{N-1}, u >0, N=\frac{2n}{n-2}, n\geq 3, \quad (*) $$

tels que ($ K $ d\'esigne un compact quelconque de la vari\'et\'e Riemannienne $ (M,g) $ non necessairement compacte sans bord):

$$ \sup_K u\times \inf_M u \leq c,$$

et,

$$ \int_K u^{2n/(n-2)} dV_g \leq c, $$

Ici on est dans le cas positif. Ceci inclut le cas non loc.conf.plat.

\smallskip

////////////////////////////////////////////////

\smallskip

a) Dans le cas negatif on a:

$$ -\Delta u+\frac{(n-2)}{4(n-1)} R_g u= -u^{N-1}, u >0, \,\,\, {\rm ici}, \,\, V\equiv -1, $$

avec la compacit\'e locale:

$$ \sup_K u\leq c.$$

///////////////////////////////////////////////

\smallskip

b) Dans le cas nul, on a:

$$ -\Delta u+\frac{(n-2)}{4(n-1)} R_g u=0, u >0, \,\,\, {\rm ici},\,\, V\equiv 0.$$

avec l'in\'egalit\'e de Harnack usuelle:

$$ \sup_K u \leq c\inf_K u, $$

et,

$$ {\rm principe \,\,d'Harnack}.$$

////////////////////////////////////////////

Ici, on s'interesse au cas positif et l'equation $ (*) $ ci-dessus. Alors on commence par resoudre l'equation suivante:

$$ -\Delta v-\lambda v= n(n-2) v^{N-1}, v >0, $$

avec, $ 0< m \leq \lambda + \frac{(n-2)}{4(n-1)} R_g \leq 1/m $.

\smallskip

G\'eom\'etrisation du probleme et de la solution: on prend alors:

$$ V=n(n-2)+(\frac{(n-2)}{4(n-1)}R_g+\lambda )v^{2-N}, $$

et

$$ u=v, $$

Si on considere les metriques conformes: $ \tilde g=u^{4/(n-2)} g $, alors, la courbure scalaire de $ \tilde g $ est:

$$ R_{\tilde g}=[-\Delta u+\frac{(n-2)}{4(n-1)} R_g u]u^{1-N} =[-\Delta v+\frac{(n-2)}{4(n-1)} R_g v]v^{1-N}=n(n-2)+(\frac{(n-2)}{4(n-1)}R_g+\lambda )v^{2-N}=V>0, $$

Donc, $ V $ est courbure scalaire de $ \tilde g $.

\smallskip

Si l'operateur $ -\Delta -\lambda $ est coercif, on utilisant sa fonction de Green, on a majoration locale du volume = energie.

\smallskip

Donc le couple $ (u, V) $ est solution. Ceci inclut le cas non loc.conf.plat. On a:

\smallskip

L'in\'egalit\'e de Harnack optimale:

$$ \sup_K u\times \inf_M u \leq c,$$

et, majoration locale du volume = energie, si $ -\Delta -\lambda $ est coercif:

$$ \int_K u^{2n/(n-2)} dV_g \leq c, $$

De plus, on connait explicitment le potentiel $ V >0 $, qui doit etre $ C^{\infty} $, si on connait explicitement la fonction $ v >0 $, qui doit etre $ C^{\infty} $, pour ce qui concerne les consid\'erations numeriques.

\smallskip

//////////////////////////////////////////////////////////////////////

\smallskip

\section{ Consequence du resultat d'unicit\'e et de rigidit\'e: une in\'egalit\'e de Sobolev et une in\'egalit\'e d'interpolation}

\smallskip

Soit $ (M,g) $ une vari\'et\'e Riemannienne compacte sans bord de dimension $ n \geq 4 $ et de courbure scalaire $ S_g >0 $ partout sur $ M $ et orientable. Dans le livre d'Aubin, pour $\epsilon >0 $ assez petit $ \epsilon <S_g $, l'equation  $ 4\dfrac{n-1}{n-2} \Delta u+\epsilon u =u^{N-1}, u>0 $ a une solution et le resultat d'unicit\'e dit que cette solution est unique. (Voir dans le livre d'Aubin, la fonctionnelle associ\'ee \`a cette equation possede un infimum atteint par cette solution unique). On obtient(voir aussi l'article de L.Veron et J.R. Licois et l'article sur Hal de J.Dolbeault, Esteban, Loss), une in\'egalit\'e d'interpolation et une in\'egalit\'e de Sobolev:

$ \exists \, \epsilon_0=\epsilon_0(n,M,g)>0,$ tel que pour $ 0 < \epsilon <\epsilon_0 $ on ait:

$$ \inf_{u\in H^1(M)-\{0\}} \dfrac{4\dfrac{n-1}{n-2} \int_M |\nabla u|^2dV_g+\epsilon\int_M u^2dV_g}{(\int_M |u|^N dV_g)^{2/N}}=\epsilon (vol(M))^{2/n}, $$

\smallskip

On a alors, pour tout $ u \in H^1(M), $ on a:

$$ 4\dfrac{n-1}{n-2} \int_M |\nabla u|^2dV_g \geq \epsilon\left ( \left ( \int_M |u|^N dV_g \right )^{2/N} (vol(M))^{2/n}- \int_M u^2dV_g \right ), $$

ou encore,

$$ \forall u\in H^1(M),\,\, 4\dfrac{n-1}{n-2} ||\nabla u||_2^2\geq \epsilon \left ( (vol(M))^{2/n}||u||_{2^*}^2-||u||_2^2 \right )\geq 0. $$

avec $ 2^*=N=\dfrac{2n}{n-2} $ pour $ n\geq 4 $.

Finalement: il existe $ C=C(n,M,g) >0 $ telle que:

$$ \forall u\in H^1(M), \,\, ||\nabla u||_2^2\geq C \left ( (vol(M))^{2/n}||u||_{2^*}^2-||u||_2^2 \right )\geq 0, $$

avec,  $ n\geq 4 $ et $ 2^*=N=\dfrac{2n}{n-2}.$

\section{ Consequence du resultat de compacit\'e en dimension 2: application de: compacit\'e avec energie ou volume born\'es: in\'egalit\'es du type Moser-Trudinger:}

\smallskip

Compacit\'e avec contrainte: 

\smallskip

1) Topologie des espaces de domaines \`a bord: en plus du fait qu'on a convergence de domaines avec bord et metriques conformes, suites de domaines \`a bord. geometrie des metriques. Espaces d'espace ou espaces de domaines a bord muni de metriques Riemanniennes. Volume fix\'e: metrique blow-up, eclatement, ou effondrement d'espaces. ou compacit\'e de metriques, non effondrement d'espaces. Ceci est juste un apercu, car il faut faire attention aux hypotheses des th. il faut que les donn\'ees soient formul\'ees en termes de \underbar{distance geodesique} en general.

\smallskip

2) Resultat de compacit\'e de fonctions dans les EDP. et les espaces de fonctions. En topologie et en analyse fonctionnelle.

\smallskip

3) En chimie: emballement thermique (thermal runaway, jusqu'a l'eclatement, tend vers l'abime) ou non emballement (compacit\'e: objet ne s'abime pas)

\smallskip

4) Problemes variationnels:

Par exemple, regardons le cas regulier, sans singularit\'e:

On a vu dans le print,"Cas d'existence de solutions d'EDP": que la probleme variationnel:

$$ \mu = \inf \{ ||\nabla u||_2^2, u\in \dot H_1^2(\Omega), \int_{\Omega} Ve^u=1, \} $$
 a une solution positive $ u >0 $ avec la condition $  b|\Omega| < 1 $.

$$ \Delta u = \lambda V e^u, u=0, \,\, {\rm au\,\, bord},\,\, \lambda >0, $$

a) En supposant: $ 0 <a \leq V\leq b <+\infty $ avec la fonction propre, le coefficient, le multiplicateur de lagrange verifie $ \lambda \times a \leq C <+\infty$

b) En supposant $ V $ Lipschitz, on a la compacit\'e. et en fait $ \lambda \not \to 0 $ uniform\'ement. Donc, on a: $ \mu >0 $ uniform\'ement en $ u \in \dot H_1^2(\Omega) $ avec $ \int_{\Omega} Ve^u = 1 $.

Donc avec la condition $ b|\Omega| <1$:

$$ \forall u \in \dot H_1^2(\Omega), \int_{\Omega} Ve^u=1 \Rightarrow ||\nabla u||_2 \geq \mu >0.$$

Donc, la contrapos\'ee donne:

$$ ||\nabla u ||_2 < \mu \Rightarrow \int_{\Omega} Ve^u \not = 1. $$

On va  voir qu'on peut avoir une in\'egalit\'e du type Moser-Trudinger, et les constantes, et $ V $ determinent cette nouvelle in\'egalit\'e.

Soit: $ u\in \dot H_1^2(\Omega), u\not \equiv 0, $ et $ v=\frac{\mu}{2} \frac{u}{||\nabla u||_2} $. Alors, $ ||\nabla v||_2 =\mu/2<\mu $, donc, $ \int_{\Omega} V e^v \not = 1 $.

Supposons qu'il existe $ u_0 $ et $ u_1 $ non nuls tels que $ \int_{\Omega} Ve^{v_0} <1, \int_{\Omega} V e^{v_1} >1 $, avec $ v_0, v_1 $, les fonctions construites  a partir de $ u_0,  u_1$, comme $ u $ et $ v $. On a: $  v_0=\frac{\mu}{2} \frac{u_0}{||\nabla u_0||_2},  v_1=\frac{\mu}{2} \frac{u_1}{||\nabla u_1||_2} $.

Alors, on considere le chemin $ w_t=tv_0+(1-t) v_1$, alors, $ ||\nabla w_t||\leq \mu/2 <\mu $, d'ou, $ \int_{\Omega} Ve^{w_t} \not = 1, \forall \, t \in [0,1] $, or par le theoreme des valeurs intermdiaires en considerant la fonction continue de $ t, g(t)=\int_{\Omega} V e^{w_t} dx $, on aura un $ t_0 $ tel que $ g(t_0)=1=\int_{\Omega} Ve^{w_{t_0}} $, ce qui est constardictoire.

c) Donc, on a soit tout le temps $ \int_{\Omega} Ve^v<1$ ou tout le temps, $ \int_{\Omega} Ve^v > 1 $.

Or, en considerant les fonctions $ h\leq 0, h\not \equiv 0 $, on $ \int_{\Omega} Ve^h \leq b|\Omega|<1$. (Par exemple si $ u \in \dot H_1^2(\Omega),  u\not \equiv 0 $, $ |\nabla (|u|)|=|\nabla u| $, on peut prendre $ h=\frac{\mu}{2}\frac{-|u|}{||\nabla u||_2} $)

d) Donc, on a tout le temps:

$$ \int_{\Omega} Ve^v < 1. $$

e) Soit, 

$$ \nu=\sup \{ \int_{\Omega} V e^v, v=\frac{\mu}{2}\frac{u}{||\nabla u||_2}, u \in \dot H_1^2(\Omega), u\not \equiv 0 \}, $$

Par la compacit\'e de l'injection de Moser Trudinger, ce $ \sup $
 est atteint: $ ||\nabla v||_2\leq \liminf_i ||\nabla v_i||_2 =\mu/2 <\mu $ et la compacit\'e de l'injection de Moser-Tridinger: $\int_{\Omega} Ve^{v_i} \to \int_{\Omega} Ve^v =\nu \leq 1 $. Alors, $ v \in \dot H_1^2(\Omega), ||\nabla v||_2 \leq \mu/2 < \mu $, d'ou, $ \int_{\Omega} Ve^v \not =1$. Donc $ \nu <1 $.

f) Finalement:

$$ \forall \, a, b, A >0, \, b|\Omega|<1, \exists \,\mu >0, \exists \, 0 < \nu <1, \,\, \forall u \in \dot H_1^2(\Omega)-\{0\}, \,\, \int_{\Omega} V e^{\frac{\mu}{2} \frac{u}{||\nabla u||_2}} dx \leq \nu < 1.$$

C'est une in\'egalit\'e du type Moser-Trudinger, la consition, $ b|\Omega| <1 $ et les constantes $ \mu>0 , 0 <\nu <1 $ determinent cette in\'egalit\'e. Donc, pour $ a,b,A>0 $ avec la condition $ b|\Omega| < 1 $, on obtient:

$$ \exists \,\mu >0, \exists \, 0 < \nu <1, \,\, \forall u \in \dot H_1^2(\Omega)-\{0\}, \,\, \int_{\Omega} V e^{\frac{\mu}{2} \frac{u}{||\nabla u||_2}} dx \leq \nu < 1.$$

g) On peut faire la meme chose avec, le resultat de compacit\'e avec singualrit\'e au bord. On a une in\'egalit\'e du type Moser-Trudinger avec singularit\'e au bord, avec des constantes $ \mu >0, 0 < \nu <1$ et la condition sur $ b: b \times \int_{\Omega} \frac{1}{|x-x_0|^{2\alpha}} dx <1, \alpha \in ]0,1/2[ $, $ x_0 \in \partial \Omega $ et $ \Omega $ domaine analytique:

$$ \exists \,\mu >0, \exists \, 0 < \nu <1, \,\, \forall u \in \dot H_1^2(\Omega)-\{0\}, \,\, \int_{\Omega} \frac{V}{|x-x_0|^{2\alpha}} e^{\frac{\mu}{2} \frac{u}{||\nabla u||_2}} dx \leq \nu < 1.$$

h) On peut mettre la valeur absolue pour $ u $ et remplacer $ u $ par $ |u| $ dans ces in\'egalit\'es.

Ces fonctions $ u $ font apparaitre l'Eq. de la courbure scalaire prescrite en dimension 2. Donc: pour un $ u \in \dot H_1^2(\Omega) $, $ \int_{\Omega} e^{ku} dx \approx |\Omega|_{g= e^u \delta}^k $ correspond \`a un volume ou une surface, \`a un r\'eel positif, $ k >0 $ pr\'es (en utilisant l'in\'egalit\'e de Holder par exemple). Le terme $ ||\nabla u||_2^2 \approx |\partial \Omega|_{g=e^u\delta}^2 $, (en utilisant l'inegalit\'e de Cauchy-Schwarz par exemple), correspond au perimetre, en passant par les fonctions BV, voir l'article de. O. Druet dans Numdam, 2001-2002: in\'egalit\'es de Sobolev et in\'egalit\'es isop\'erimetriques.

\smallskip

Dans le cas avec singularit\'e au bord, on a des in\'egalit\'es entre surface et perimetre avec la metrique $ g=\frac{e^u}{|x-x_0|^{2\alpha}} \delta, \alpha \in ]0, 1/2[, x_0\in \partial \Omega $ et $ \Omega $ analytique.

\smallskip

Cela s'applique aussi a l'operateur: $ -div (e^{\epsilon |x|^2/2} \nabla) $. La fonctionnelle est: $ \int_{\Omega} e^{\epsilon |x|^2/2} |\nabla u|^2 dx $.

\smallskip

Donc, ces in\'egalit\'es en dimension 2 de Moser-Trudinger, mettent en relation la surface conforme (volume conforme) et le p\'erimetre (conforme): ce sont des in\'egalit\'es isoperimetriques particulieres.

\smallskip

Il y a aussi l'interpretation en physique, en termes d'energies. (Qui correspond aussi au cas ou l'operateur n'est pas n\'ecessairement le laplacien).

\smallskip

On peut prendre comme contrainte: $ \int_{\Omega} Ve^u = k $, $ k \geq 1 $ au lieu de $ 1 $. Dans ce cas, la condition sur $ b $ est: $ b|\Omega| < k $.

\smallskip

En mettant la valeur absolue dans $ u $, on obtient: in\'egalit\'es du type Moser-Trudinger:

\smallskip

Avec, $ a, b, A>0 $, $ V $ et $ b|\Omega|< 1$:

$$ \exists \,\mu >0, \exists \, 0 < \nu <1, \,\, \forall u \in \dot H_1^2(\Omega)-\{0\}, \,\, \int_{\Omega} V e^{\frac{\mu}{2} \frac{|u|}{||\nabla u||_2}} dx \leq \nu < 1.$$

et, avec $ a, b, A >0 $, $ \alpha \in ]0, 1/2[$, $ V $, $ b\int_{\Omega} \frac{1}{|x-x_0|^{2\alpha}} dx < 1$ et $ \Omega $ analytique:

$$ \exists \,\mu >0, \exists \, 0 < \nu <1, \,\, \forall u \in \dot H_1^2(\Omega)-\{0\}, \,\, \int_{\Omega} \frac{V}{|x-x_0|^{2\alpha}} e^{\frac{\mu}{2} \frac{|u|}{||\nabla u||_2}} dx \leq \nu < 1.$$

\bigskip

\section{Quelques remarques: I:}

\smallskip

1) 

-Th\'eorie de Yang-Mills, Equations de Yang-Mills. 

\smallskip

-Th\'eorie conforme des champs de Liouville. Th\'eorie de la gravitation quantique de Liouville. 

\smallskip

-Th\'eorie de la gravitation quantique; en relativit\'e g\'en\'erale($ n=3 $), de Kaluza-Klein ($ n=4$), dans la th\'eorie des cordes $ (n=5,6$) et des supercordes ($ n=9 $). Particules: Axion, graviton. Supersymetrie, symetrie quantique(symetrie conforme): Eq. de Yamabe. Eq. d'Einstein-Lichnerowicz (quantique relativiste) : particules de spin entier et de spin demi-entier. 

\smallskip

-Eq. de Schrodinger(quantique non relativiste).

\smallskip

-Modele des interactions des particules. 

\smallskip

-Modele cosmologique. Trous noirs.

\smallskip

2) 

\smallskip

-Eq. de Liouville. Eq.du type Liouville (Eq. de la courbure scalaire prescrite en dimension 2). 

\smallskip

-Eq. de Yamabe. Eq. de la courbure scalaire prescrite en dimension $ n\geq 3 $. Eq. du type courbure scalaire prescrite. Eq.d'Einstein-Lichnerowicz.

\smallskip

3) 

D-Branes: Eq. avec condition de Dirichlet au bord (vari\'et\'e \`a bord, bord regulier: $ n=2 $ Eq. de Liouville ou de courbure scalaire et $ n\geq 3 $ Eq. d'Einstein-Lichnerowicz, Eq. de Yamabe). 

\smallskip

Masse positive, expansion de l'univers. Operateurs coercif.

\smallskip

4) 

-Enroulement. Torsion. Distrotion. Noeuds.

\smallskip

-Stability. Coh\'erence. 

\smallskip

-Effondrement d'espaces: blow-up. Non-effondrement d'espaces: compacit\'e, regidit\'e.

\bigskip

{\bf Equation de la courbure scalaire prescrite (et du type Yamabe) en dimension 4 et champs de Yang-Mills}:

\smallskip

a) exemple: 1-si on note $ SW $ la vari\'et\'e hyperbolique sans bord et orientable de dimension 3, de Seifert-Weber, alors $ SW\times S_1$ est une vari\'et\'e loc.conf.plate de courbure scalaire constante stric.negative de dimension 4, elle est li\'e aux champs de Yang-Mills. 2- En dimension 4, aussi on considere, la vari\'et\'e de Davis ou le produit de deux surfaces de courbure sectionelle non oppos\'ees et de courbure scalaire constante strc.negative (donc, non-loc-conf.plate $ S_{-1} \times S_{-1} $ avec $ S_{-1} $ une surface de courbure scalaire constante $ <0 $, $ S_{-1} \times S_{-1} $ est d'Einstein), on obtient une vari\'et\'e Riemmannienne li\'e aux champs de Yang-Mills. On peut considerer sur ses vari\'et\'es l'Equation de la courbure scalaire prescrite (dimension 4).

\smallskip

Dans le cas positif, on peut prendre $ S_4, S_2\times S_2, S_3\times S_1 $.

\smallskip

b) En physique (description de la force nucl\'eaire responsable de la cohesion des protons-neutrons dans le noyau) le cas le plus important est port\'e aux vari\'et\'es de dimension 4 avec une metrique Riemannienne ou Lorentzienne.

\smallskip

c) Regarder l'article de, Andrzej Derdzinski, dont le titre est: "Riemannian manifolds with harmonic curvature", il y a la fonctionnelle de Yang-Mills et pour lui la deriv\'ee du tenseur de courbure donne les points critiques de la fonctionelle de Yang-Mills. En fait il s'agit de la vari\'et\'e, de la metrique et de la connection. C'est la courbure de Riemann qui est un champ de Yang-Mills.

\smallskip

d) Pour ce qui est des champs de Yang-Mills en dimension 4 et la th\'eorie de Seiberg-Witten, apparait la fonctionelle de Yamabe et l'invariant de Yamabe et la courbure scalaire et les metriques conformes, dans des articles de M.Gursky par exemple en dimension 4:

\smallskip

 Si on prend une vari\'et\'e Riemannienne compacte sans bord $ (M^4,g) $, de dimension 4 avec champ de Yang-Mills, cela veut dire que le tenseur de Riemann, $ Riem=F_{\nabla} $ est un champ de Yang-Mills. Comme la courbure scalaire $  S_g $ se d\'eduit du tenseur de Riemann (elle est li\'ee au tenseur de Riemann), alors, $ S_g $ est un objet de la physique et l'Equation de la courbure scalaire prescrite, en dimension 4, sur les vari\'et\'es de Yang-Mills, est li\'ee \`a la physique puisqu'elle contient $ S_g $. De meme pour la metrique conforme, qui est li\'ee \`a $ g $, la metrique Riemannienne $ g $ est appel\'e en dimension 4 un instanton.

\smallskip

La th\'eorie de Yang-Mills existe sur des vari\'et\'es non-compactes ou completes. (non compact, or complete Manifolds with harmonic curvature). Voir Taubes et les articles de Gabor Etesi et Yawei Chu. Exemples: vari\'et\'es d'Einstein (Einstein manifolds) $ n\geq 3 $ et les vari\'et\'es loc.conf.plates de courbure scalaire constante (pour $ n\geq 4 $, en particulier pour $ n=4 $) ou Ricci parallel, en particulier les vari\'et\'es de courbure sectionelle constante.

\smallskip

 Du point de vue de la physique, en dimension 4, on peut considerer les champs de Yang-Mills sur une vari\'et\'e compacte ou complete ou non compacte et localement on a une notion de champs de Yang-Mills. On peut considerer au depart la vari\'et\'e et puis on regarde ce qui se passe localement, on a encore des champs de Yang-Mills. Puis, on regarde ce qui se passe globalement ou localement (mesures, obsevation, \`a partir des donn\'ees et du champs de Yang-Mills).

\smallskip

Comme pour l'equation de Schrodinger: La donn\'ee est la courbure scalaire prescrite $ V $ (elle est prescrite, pulsion ou signal ou potentiel), la solution $ u $ est une "fonction d'onde" qui lie $ V $ \`a $ S_g $, $ S_g $ est un "champ scalaire de Yang-Mills". L'equation de la courbure scalaire prescrite en dimension 4 et sur une vari\'et\'e de Yang-Mills de dimension 4 est un objet de la physique. 

Aussi, c'est ecrit dans l'article de T.H. Parker (Gauge Theories on four dimensional Riemannian manifold, Comm.Math.Physics, 1982), si on considere un terme ("potentiel de Higgs") du type $ P(u)= au^4+mu^2 $ avec $ u >0$ la solution et $ a, m $ deux fonctions r\'eelles avec $ m \leq 0 $ possible, l'equation bosonique dans un champs de Yang-Mills devient du type Yamabe (ou de type courbure scalaire prescrite). ($ \Delta_g u + (\frac{s_g}{6}+2m) u=-4au^3 $, avec $ \Delta_g=-\nabla^i \nabla_i $, $ s_g $ la courbure scalaire). On conclut, voir l'article de T.H. Parker (dans un champ de Yang-Mills, l'equation bosonique, induit (avec le fait qu'on a un champ de Yang-Mills, les lagrangiens s'annulent et la d\'eriv\'ee par rapport \`a la connexion est nulle), que le r\'esidu est nul et constitue une solution de couplage des champs).(le lagrangien Yang-Mills bosonique, $ Y=B_1+B_2 $ avec $ B_1 $ la fonctionnelle de Yang-Mills et $ B_2 $ la lagrangien bosonique, la loi de la particule: $ \dot Y= dY=0= dB_1+dB_2=0 $, si on se place dans un champ de Yang-Mills, alors, $ dB_1=0 $, ce qui implique que $ dB_2=0 $, si de plus, $ u >0 $ est solution de l'equation bosonique (cela veut dire que $ \frac{\partial B_2}{\partial u}=0$), ceci implique que $ \frac{\partial B_2}{\partial \nabla}=0 $ et le r\'esidu est nul: $ \Sigma <\nabla_j u|\rho (e_i) u>=0 $ et donc, les deux equations du systeme de couplage sont v\'erifi\'ees. Donc l'equation bosonique (du type courbure prescrite) dans un champ de Yang-Mills suffit \`a determiner le systeme. Ces deux conditions, champs de Yang-Mills et equation bosonique constituent une solution particuliere du systeme.(En partant de la formulation de Yang-Mills bosonique, $ \dot Y=dY=0$. Donc, dans la formulation de Yang-Mills bosonique, l'Eq de la courbure scalaire prescrite ou du type courbure scalaire prescrite et le champ de Yang-Mills realisent un systeme physique.) On a aussi, l'equation de la courbure scalaire prescrite ou de type courbure scalaire prescrite y apparait (dans la formulation de Yang-Mills bosonique, $ \dot Y=dY=0 $), c'est une equation de la physique. )

(On a: $ Y=B_1+B_2 $, donc: $ dY=dB_1+dB_2=dB_1+dB_{21}+dB_{22} $, donc: Eq de la courbure scalaire prescrite ou du type courbure scalaire prescrite (Eq bosonique)+Yang-Mills: $ \Rightarrow dB_{22}=dB_1=0 \Rightarrow dB_{21}=dY $ et dans la formulation de Yang-Mills bosonique: $ \dot Y=dY=0 $, donc: $ dB_{21}=dY=0 $.)

On peut aussi considerer chaque lagrangien seul dans n'importe quelle formulation, par exemple ici, celui de Yang-Mills seul et celui bosonique seul (On peut considerer le lagrangien bosonique et choisir une connexion particuliere de Yang-Mills. dynamique du boson de Higgs dans un champs de Yang-Mills ou vari\'et\'e de Yang-Mills). Ici (dans l'article de Parker, T.H) on a le lagrangien de Yang-Mills bosonique, l'Eq de la courbure scalaire prescrite ou du type courbure scalaire prescrite y apparait, c'est l'Eq bosonique dans un systeme de Yang-Mills bosonique.

\smallskip

{\bf En dimension 2: equation de Liouville:} 

\smallskip

a) Equation de la courbure scalaire prescrite sur la sphere de dimension 2 (vortex equation) ou un ouvert de $ {\mathbb R}^2 $ (vortex equation) ou surface de dimension 2 (vortex equation).  Mean-Field-equation.

\smallskip

b)Voir les articles de Crandall-Rabinowitz, De Figueiredo-Lions-Nussbaum, Chen-Li, pour la provenence de ses equations: G\'eometrie (courbure de Gauss), Gazs, combustion, astronomie et astrophysique.

\smallskip

c) En dimension 2, aussi, c'est un cas particulier des champs de Yang-Mills, la theorie de Glashow-Weinberg-Salam, modelise les interactions (faibles) electromagnetiques des particules. C'est la theorie de jauge (theorie de champs) avec un groupe de jauge (groupe de symetries locales) $ U(1) \times SU(2) $ 

voir le livre de G. Tarantello Self-dual Gauge theories.

\smallskip

d) L'equation de Liouville ou de type Liouville apparait aussi dans le ph\'enom\`ene de "cordes cosmiques", un objet de l'univers "cosmic strings" et \`a ne pas confondre avec la theorie des cordes. Voir l'article de J. Spruck et Yisong Yang (cosmic strings).(On prend une metrique Lorentzienne $ ds^2=-dt^2 +dz^2+g_{ij}dx^idx^j $, la metrique $ g_{ij} $ est definie sur une surface, $ M $ de dimension 2 et apres on peut choisir $ g_{ij} $ conforme a une autre metrique, par exemple $ g_{ij}=e^u \delta_{ij} $ et on choisit le champ particulier). C'est l\`a, par exemple l'id\'ee g\'en\'erale.

\smallskip

e) En dimension 2: c'est aussi, la theorie des champs de Liouville en dimension 2. Qui a \'et\'e etudi\'ee pour des surfaces \`a courbure n\'egative et plus r\'ecemment sur la sphere de dimension 2. (LCFT, Liouville conformal field theory, on peut mettre la th\'eorie de Liouville avec de la gravit\'e 3D (2+1), la theorie de Liouville est une th\'eorie des champs ou des cordes particuli\`ere, l'action ou le lagrangien dans la theorie de Liouville donne l'Eq de Liouville ou de courbure scalaire prescrite en dimension 2). 

f) Il y a aussi la th\'eorie de Chern-Simons ou l'Eq de Liouville avec singularit\'es apparait, formulation (2+1) dans un espace de Minkowski. Voir l'article de C.C.Chen-C.S.Lin-G.Wang: Concentration pheneomena of two vortex solutions in Chern-Simons model, Ann. Scuola. Nor. Sup. Pisa.2004. 

g) Dans la th\'eorie de Liouville, le lagrangien, l'equation de mouvement de la corde ou des particules est une g\'en\'eralisation de l'equation de Liouville ou equation de Gauss ou de courbure de Gauss. En dimension 2, les courbures scalaire, de Ricci, de Gauss ou sectionnelle sont identiques, c'est la 'courbure'(valeur) moyenne (mediane), d'ou le terme, champs moyen.

\smallskip

De meme ici, du point de vue de la physique, de la chimie, astronomie ou la physique quantique, en dimension $ n=2 $, on peut considerer des sources, des champs sur un ouvert de $ {\mathbb R}^2 $ ou une surface compacte ou complete ou non compacte et localement on a une notion de champs local. On peut considerer au depart l'ouvert ou la surface et puis on regarde ce qui se passe localement, on a encore des champs. Puis, on regarde ce qui se passe globalement ou localement (mesures, obsevation, \`a partir des donn\'ees et du champs ou des parametres).

\smallskip

Ici aussi, comme pour l'equation de Schrodinger: La donn\'ee est la fonction $ V $, ($ V $ est la source ou pulsion ou le signal ou le potentiel), la solution $ u $ est une "fonction d'onde" ou emission, qui lie $ V $ \`a la donn\'ee de d\'epart, par exemple $ S_g $ (dans le cas d'un ouvert $ S_g=0 $, $ g=\delta $ metrique euclidienne), $ S_g $ est un "champ scalaire". Il se peut que ce ne soit pas $ S_g $ comme pour l'operateur $ \Delta + \epsilon (x_1\partial_1+x_2\partial_2) $, $ \Delta=\partial_{11}+\partial_{22} $.

\smallskip

{\bf En dimension $ n\geq 3 $: Relativit\'e g\'en\'erale ($ n=3 $) et Cosmologie quantique ($ n\geq 3 $):}

\smallskip

Cette id\'ee de partir d'un espace Lorentzien $ n+1 $ et de choisir une metrique conforme, sur la vari\'et\'e $ M $, en partant de l'espace temps $ M\times (0,\delta) $,  se g\'eneralise aux dimensions superieures. (Il y a la correspondance, ADS-CFT(Anti-de-Sitter, Conformal field theory) propos\'ee par Juan Maldacena ($\Lambda <0 $, constante cosmologique strictement n\'egative dans l'equation d'Einstein). Les equations d'Einstein classiques et quand la dimension $ n\geq 3 $ c'est la "quantum-cosmology qui correspond au cas $ \Lambda =0 $, equation d'einstein du d\'ebut, (1915)). On a $ ds^2=-dt^2+g_{ij} dx^idx^j $ avec $ g_{ij} =u^{4/(n-2)} g_0 $, une metrique conforme a celle de d\'epart $ g_0 $ sur $ M $. 

L'id\'ee g\'en\'erale est la suivante: (On part de la formulation Lorentzienne $ (n+1), M\times (0, \delta) $, d'Einstein en relativit\'e g\'enerale et on introduit sur la variable espace des metriques conformes dans la fonctionelle du champ (ici au lieu de la fonctionelle du champ comme dans la th\'eorie de Yang-Mills, on a l'equation d'Einstein), on a des champs particuliers li\'es au changement de metriques conformes, tout cela dans le but d'etudier les interactions des particules et des astres (gravitation d\^ue \`a la courbure de l'espace temps ou \`a la formulation en $ n+1 $ et les ondes qui sont les effets de la gravitation ou des champs, en particulier gravitationnel, un autre champ est le champ $ T $ d'energie-impulsion, peut etre vu comme une donn\'ee, l'interaction est ecrite dans l'equation d'Einstein ci-dessous, gravitation, $ Ricci^{\gamma}, R^{\gamma} $, et les donn\'ees ou un autre champ, est le tenseur $ T $). C'est l'adpatation \`a $ (n+1) $ de la th\'eorie d'Einstein, en relativit\'e g\'en\'erale, qui est en $ (3+1) $. L'espace Anti-de-Sitter est une solution particuliere des equations d'Einstein, comme l'espace de Minkowski ou la metrique de Schwazrchild, l'inconnu ici est un espace Lorentzien $(V, \gamma) $ de metrique Lorentzienne $ \gamma $, qui verifie l'equation d'Einstein $ Ricci^{\gamma}-\frac{1}{2} R^{\gamma} \gamma = T $, avec $ Ricci^{\gamma}, R^{\gamma}, T $ respectivement le tenseur de Ricci et la courbure scalaire et le tenseur d'Energie-impulsion. Dans certains cas, cela revient \`a resoudre ce qu'on appelle le probleme de Cauchy en relativit\'e.par exemple $ T=0 \Rightarrow Ricci^{\gamma}=0 $. Un autre cas particulier, on en parl\'e ci-dessus, est de trouver des solutions du type $ V=M\times (0,\delta) $ et $ \gamma=-dt^2+g_{ij} dx^idx^j $ avec $ g=u^{4/(n-2)} g_0 $.  Resoudre les equations des contraintes pour $ (V, \gamma) $ particuliers avec la m\'ethode conforme de Lichnerowicz. On prend $ V=M\times (0,\delta) $ avec $ \gamma = -dt^2+g_{ij}dx^i dx^j $, on transfome les equations d'Einstein en 2 Equations des contraintes (Probleme de Cauchy, par les relations Gauss-Codazzi), puis on prend $ g=u^{4/(n-2)} g_0 $ sur $ M $, on transforme les 2 equations des contraintes en un systeme de deux equations elliptiques. Einstein scalar field Lichnerowicz equations. Voir, Hebey, Pollack, Pacard, Choquet-Bruhat, Bartnik.)

\smallskip

{\bf Remarque:} Pour ce qui est de la th\'eorie de Yang-Mills, on s'interesse \`a une fonctionelle du champ, on peut se ramener \`a la dimension 2 et utiliser la methode conforme comme dans Spruck-Yisong (cosmic-strings) ou d'autres th\'eories avec fonctionnelle de champs sur une vari\'et\'e Lorentzinne de dimension 4. ($ ds^2=-dt^2+dz^2+g_{ij} dx^idx^j $) (exemple, Glashow-Weinberg-Salam, combine fonctionnelle de Yang-Mills et theorie conforme en 2 dimensions). Ici dans ADS-CFT ($ \Lambda <0 $ ou bien $ \Lambda=0 $ et on retrouve les Equations d'Einstein classiques avec $ n=3 $ pour la relativit\'e g\'enerale et $ n\geq 3 $ pour la cosmologie quantique ou "quantum cosmology"), on remplace la fonctionnlle par l'equation d'Einstein. De meme, ici, dans la relativit\'e g\'enerale ($ n=3 $) et pour la cosmologie quantique ou "quantum cosmology" ($ n\geq 3 $), on remplace la fonctionnlle par l'equation d'Einstein.

\smallskip

a) il y a une application et r\'ealisation en physique, en $(3+1) $ et $ (4+1) $ et $ (5+1) $ et $ (6+1)$ de la correspondance de ADS-CFT, de Juan.Maldacena (on voit bien le probleme du bon espace-temps). Et aussi $ (2+1) $ (LCFT, Liouville conformal field theory).

La correspondance ADS-CFT correspond a l'equation d'Einstein avec une constante cosmologique $ \Lambda <0 $ (Anti-de-Sitter): $ Ricci^{\gamma}-\frac{1}{2}R^{\gamma} g+\Lambda g= T=0 $. Si $\Lambda =0$, on retrouve l'equation d'Einstein (1915) ecrite plus haut, et les th\'eories en question sont: $ n=3 $, la relativit\'e g\'en\'erale et $ n\geq 3 $ la cosmologie quantique ou "quantum cosmology". 

\smallskip

Il y a des auteurs qui font r\'ef\'erence \`a la correspondance ADS-CFT, qui essaie de trouver le bon espace-temps avec $ \Lambda <0 $. Dans la th\'eorie classique d'Eisntein($ \Lambda=0 $), c'est la "quantum cosmology", trouver le bon espace-temps pour $ n\geq 3 $.

\smallskip

b) Pour ce qui est des equations d'Einstein et methode conforme: ($ n\geq 3 $ ("quantum-cosmology") et $ n=3 $ qui correspond \`a la th\'eorie de la relativit\'e g\'en\'erale:) Voir les articles de Choquet-Bruhat-Isenberg-Pollack, dans general relativity and quantum cosmology. Cosmologie quantique pour $ n\geq 3 $. (il y a des auteurs qui font r\'eference \` la correspondance ADS-CFT, qui essaie de trouver le bon espace-temps pour $ n \geq 3 $ avec $ \Lambda <0$. Ici, on a la th\'eorie classique d'Einstein ($ n=3 $ relativit\'e g\'en\'erale et "quantum cosmology" pour $ n\geq 3$ avec  $ \Lambda =0 $)).

\smallskip

(Equation d'Einstein classique: $ \Lambda =0 $, (1915) cosmologie). Voir le monograph de: E.Hebey, D. Pollack, F. Pacard, Y. Choquet-Bruhat, R. Bartnik: dans la th\'eorie des champs de Klein-Gordon massive le potentiel $ V(\Psi)=\frac{1}{2} m^2 \Psi^2 $. (les donn\'ees $ (\Psi, \sigma, \tau, \pi) $ peuvent etre choisies librement ainsi que la metrique $ g_0 $ et la vari\'et\'e). 

\smallskip

Voir aussi les articles de Choquet-Bruhat-Isenberg-Pollack, la vari\'et\'e Riemanninne $ (M, g_0) $ peut etre compacte sans bord ou non compacte. Voir aussi la Th\`ese de C.Valcu.

c) On part des equations de contraintes qui derivent des equations d'Einstein et qui sont obtenues par les equations de Gauss-Codazzi (on contracte le tenseur de Riemann dans les equations de Gauss-Codazzi pour avoir le tenseur de Ricci). On contracte le tenseur de Ricci pour faire apparaitre la courbure scalaire et on separe les termes en variables temporelles et spaciales. Puis on utilise une metrique conforme.

d) Inversement, d\`es qu'on a une solution des equations de contraintes avec champ scalaire, il existe une solution des equations d'Einstein, un espace-temps maximal unique avec une metrique Lorentzienne, d'apres Choquet-Bruhat et Choquet-Bruhat-Geroch. C'est ce qu'on appelle le formalisme de Choquet-Bruhat-Geroch-Lichnerowicz.(Il existe un espace temps $(L,\gamma)$ 'maximal' et unique \`a isometrie pr\'es avec $ \gamma $ Lorentzienne, et un plongement $ i:M\to L $, tel que $ i^*(\gamma)=g $, $(L,\gamma) $ 'prolonge' $ (M,g) $ et les autres parametres, comme le champ scalaire par exemple. A l'instant $ t=0 $ on a $ (M,g) $ et la fonction d'onde $ u>0 $, $ g=u^{4/(n-2)} g_0 $ et localement $ L=(-t_0,t_0) \times M $). L'espace $ (M,g) $ est la deformation de l'espace temps $ (L, \gamma) $ a l'instant $ t=0 $ ou bien dynamique des particules a l'instant $ t=0 $. Voir le print de B. Premoselli et l'article de A. Carlotto (The general relativistic constraints equations).

Pour chaque solution des equations des contraintes on a une solution des equation d'Einstein, on a une multitude d'univers ou plusieurs espace-temps. La multitude d'espace-temps, modelise aussi la multitude d'etats quantiques, en particulier c'est la combinaison entre la gravit\'e et les differents champs dont le champ electromagnetique et d'autres champs: Yang-Mills, Klein-Gordon, fluides, etc. Cela veut dire que du point de vue mathematique on a plusieurs espace-temps alors que, du point de vue de la physique, il y a plusieurs etats quantiques ou plusieurs fonctions d'ondes a l'instant $ t=0 $, dans le meme espace ou vari\'et\'e $ M $ avec des propri\'et\'es concernant les fonctions d'ondes ou etats quantiques (Plusieurs ondes au meme instant qui sont dans un espace commun, et, differents au meme temps). Ceci est un modele mathematique pour expliquer l'existence de plusieurs etats quantiques ou fonctions d'ondes, qui sont a la fois differents, dans plusieurs espaces differents, et dans un meme espace au meme temps: $ (M,g_0, \psi, \sigma, \tau, \pi) $, modele de base commun a tous les espaces, qui est lui meme un 'univers' de 'parametres'. Ph\'enomene physique $ \leftrightarrow $ Propri\'et\'e mathematique: i) Etat quantique $ \leftrightarrow $ fonction d'onde $ u >0 $, et, ii) extra-dimensions ou dimensions supplementaires $ \leftrightarrow $ propri\'et\'e des fonctions d'ondes: $ \sup u = f(\inf u) $, (notions d'enroulement, de torsion et les valeures $ (\sup,\inf)$). 

\smallskip

Dans un article de Choquet-Bruhat (1968), on ne peut pas avoir existence globale d'un espace temps car d'apres Penrose et Hawking, les Eq de la relativit\'e d'Einstein developpent des singularit\'es. Mais on a existence locale d'un espace temps $ (-t_0, t_0)\times M, t_0 >0 $ et unicit\'e locale et globale.

\smallskip

Th\'eorie de la graviation quantique: du type Kaluza-Klein. La th\'eorie de Kaluza-Klein c'est en $ (4+1) $ comme espace-temps avec une dimension cach\'ee. Les th\'eories du 'type' Kaluza-Klein, c'est en $ (5+1), (6+1), (n+1), n=9,10 $ ou $ (n+1), n\geq 3 $ avec plusieurs dimensions cach\'ees.

\smallskip

Pour ce qui est des solutions des Equations d'Einstein par la m\'ethode conforme de Lichnerowicz-Choquet-Bruhat-York: on a:

i) Si on prend $ \Psi =0 $, $ \tau \equiv constante \not =0 $, $ \sigma =\pi = 0 $, alors $ W=0 $ est une solution et $ u >0 $ est solution de l'equation de Yamabe (dans "le cas n\'egatif", la vari\'et\'e n'est pas n\'ecessairement, compacte sans bord).

\smallskip

ii) Si on prend $ \Psi\equiv constante \not =0 $, $ \tau =0 $, $ \sigma = \pi = 0 $, alors $ W=0 $ est une solution et $ u >0 $ est solution de l'equation de Yamabe (dans "le cas positif", la vari\'et\'e n'est pas n\'ecessairement, compacte sans bord).

\smallskip

iii) Si on prend $ \Psi \not \equiv constante $, $ \sigma = \pi =0 $ et $ \tau =0 $ alors $ W=0 $ est une solution et $ u >0 $ est solution de l'equation du type courbure scalaire prescrite avec $ h=S_{g_0}-|\nabla \Psi|_{g_0}^2 \leq S_{g_0} $ (dans "le cas positif", la vari\'et\'e n'est pas n\'ecessairement, compacte sans bord).

\smallskip

De meme ici, du point de vue de la physique, de l'astronomie ou de la cosmologie ou de la cosmologie quantique, en dimension $ n=3 $ ou $ n\geq 3 $, on peut considerer les champs sur une vari\'et\'e compacte ou complete ou non compacte et localement on a une notion de champs local. On peut considerer au depart la vari\'et\'e et puis on regarde ce qui se passe localement, on a encore des champs. Puis, on regarde ce qui se passe globalement ou localement (mesures, obsevation, \`a partir des donn\'ees et du champs ou des parametres, comme $ (\Psi, \sigma, \tau, \pi) $ et $ (M, g_0) $).

\smallskip

Ici aussi, comme pour l'equation de Schrodinger: La donn\'ee est la fonction $ V $ (du type courbure scalaire prescrite) qui peut etre constante (Yamabe), ($ V $ est prescrite (ou Yamabe, $ V\equiv constante $)), on peut prendre $ V=\frac{1}{2} m^2 \Psi^2 $, $ V $ est la pulsion ou le signal ou le potentiel), la solution $ u $ est une "fonction d'onde" qui lie $ V $ \`a $ S_{g_0} $ ou $ S_{g_0}-|\nabla \Psi|_{g_0}^2 $, $ S_{g_0} $ est un "champ scalaire".

\smallskip
On choisit $ W $ champ de vecteur de Killing conforme ($ D W=0 $, conformal Killing vector field), pour avoir un champ conforme et etre compatible avec la th\'eorie des champs conformes. Ici, on a pris $ W=0 $ (pour dire qu'il y a au moins toujours une solution), mais on peut prendre $ DW=0$. on a une th\'eorie des champs, gravitation+Klein-Gordon, ou Yang-Mills, ou Electromagnetique, et une th\'eorie conforme des champs (on avait deja une metrique conforme introduite par la methode de Lichnerowicz-York-Choquet-Brhuat-Isenberg-Pollack, c'est le contexte conforme. Ici conforme pour dire que les objets sont invariants par transformations conformes, 'Scale invariants' ou 'Zoom invariants').

\bigskip

{\bf En dimension $ n\geq 3$: r\'esultat d'unicit\'e et de rigidit\'e: en Biologie Math\'ematique:}

\smallskip

Interpretation en Biologie Math\'ematique: la Chimiotaxie (Chemotaxis en anglais), est l'etude du mouvement de certaines cellules (les amibes, amoebae en anglais) dans leur environement, mouvement d\^u  \`a la liberation de substances chimiques par les amibes. les amibes se deplacent vers les endroits ou la concentration de cette substance chimique est la plus grande et s'agregent, forment des agregats. L'interpretation biologique du resultat d'unicit\'e est que les amibes ne peuvent pas former d'agregats (ne se concentrent pas, un ou plusieurs points de concentrations) quand on a une solution unique et constante au systeme de Keller-Segel. Voir l'article de Lin-Ni-Takagi. Le systeme de Keller-Segel peut se reduire \`a l'etude d'une seule equation, voir l'article de Lin-Ni-Takagi, Journal of Diff. Equations, 1988. Cette formulation existe sur les vari\'et\'es compactes sans bord de dimension $ n\geq 3 $, voir T. Hillen et K.Painter.

\smallskip

{\bf Remarque sur les in\'egalit\'es du type Harnack $ \sup \,\inf $:} La th\'eorie des cordes (string theory) unifie les th\'eories de Jauge (des champs), des champs conformes, (th\'eorie des champs de Liouville, th\'eorie des cordes particuli\`ere, il est possible de melanger avec de la gravit\'e, 3D, (2+1), dans la th\'eorie de Liouville on considere l'action de Liouville ou lagrangien donnant l'Eq de Liouville), de Yang-Mills, et la gravitation quantique (Kaluza-Klein, $ n+1 , n=3, 4, 5, 6, 9, 10 $ ou plus g\'en\'erale $ n\geq 3 $). Les in\'egalit\'es du type Harnack apparaissent dans ces trois th\'eories, on peut considerer alors que cette notion est une notion de la th\'eorie des cordes. Corde qui est un 'fil' ou un ensemble  de particules.

Formulation d'Einstein + champs exterieurs (tenseur $ T $):$ Ricci^{\gamma}-\frac{R^{\gamma} }{2} \gamma = T $, avec un choix particulier de $ \gamma = -dt^2+g_{ij} dx^idx^j $ et $ g =e^u g_0 $ (dimension 2, Equation de Liouville) ou $ g=u^{4/(n-2)} g_0 $ (dimensions $ n\geq 3 $, relativit\'e generale et Kaluza-Klein), conforme. C'est ecrit dans le print de Choquet-Bruhat-Isenberg-York avec differents type de champs exterieurs $ T $, Yang-Mills en dimension 3 ou Klein-Gordon en dimension $ n \geq 3 $.

\smallskip

{\bf Remarque 2:} Il y a aussi les D-branes en dimension 2, cordes ouvertes avec condition de Dirichlet.(fait reference a l'equation du type Liouville (equation de mouvement des particules ou de la corde) avec condition de Dirichlet). Pour la theorie conforme des champs de Liouville.

\smallskip

{\bf Remarque 3:} En ce qui concerne l'existence de dimensions supplementaires (extra-dimensions) dans les th\'eories du type Kaluza-Klein: ceci est li\'e a l'existence  d'une nouvelle particule qui est l'axion qui se deplace dans ces dimensions cach\'ees, on peut les detecter ou savoir comment les eliminer avec des outils, des telescopes, voir les articles ou preprint de Horvat, Krcmar, Lakic. Il y a le telescope du CERN, et aussi d'autres appareils, pour essayer de detecter ou eliminer ces particules que sont les axions. C'est dit dans le print de Horvat, Krcmar, Lakic, les dimensions supplementaires $ (5+1), (6+1) $ sont li\'es a l'astrophysique, elles ont des r\'ealisations. Quant \`a $ (4+1)$, c'est Kaluza-Klein.

Voir le print de Collion-Vaugon sur l'id\'ee du proc\'ed\'e de r\'eduction dimensionnelle et qui est un exemple.(Dans leur exemple ils prennent le cercle $ S_1 $ qui est explicite, donc, la vari\'et\'e peut etre explicite aussi, $ S^3 $ ou d'autres vari\'et\'es sur les quelles le cercle peut agir. Dans le cas general d'une vari\'et\'e non explicite, un fibr\'e principal naturel est le fibr\'e des reperes).

\section{Remarques: II:}

\smallskip

Sur le premier article, 2003, "Differentes estimations de..."

\smallskip

Dans le cas n\'egatif, lorsque le bord de la vari\'et\'e $ \partial M $ n'est pas regulier, on est oblig\'e de supposer les solutions regulieres.

\smallskip

On suppose les solutions $ u \in C^2 $.

\smallskip

Dans ce cas on plusieurs points de vues:

\smallskip

Le point de vue d'analyse fonctionnelle et EDP, et le point de vue groupes topologiques.

\smallskip

On suppose donc que les solutions $ u\in C^2 $: le but est de prouver des estimations a priori: dans $ L^{\infty}_{loc}, C^0_{loc}, C^1_{loc}, C^2_{loc} $.

\smallskip

a) Le point de vue groupe topologique: l'estimation dans $ C^0_{loc} $ suffit, pour avoir la notion de groupe topologique localement compact: c'est une bonne estimation.

\smallskip

b) Par contre, pour le point de vue analyse fonctionnelle, EDP, ou de fonctions, comme on a suppos\'e les solutions $ u \in C^2 $, il faut que l'estimation a priori cherch\'ee soit de la taille de l'hypothese, c'est \`a dire $ C^2_{loc}$. Or pour avoir $ C^2_{loc} $, il faut utiliser les estimations de Schauder, c'est a dire qu'il faut supposer $ C^{2, \alpha}, \alpha >0 $, or cela necessite que le potentiel $ V $ soit $ C^{\alpha} $ uniform\'ement.

\smallskip

c) Tout ceci pour dire que selon le point de vue, la regularit\'e forte uniforme du potentiel est n\'ecessaire, le fait de supposer le potentiel $ V $ entre deux constantes unifor. n'est pas justifiable de maniere directe, en general, c'est pour appliquer le th. d'Ascoli, mais dans le cas n\'egatif on peut comprendre, la raison pour laquelle, on fait cette hypothese:

\smallskip

Par exemple pour $ -\infty < a \leq V \leq b <0$, si on suppose par exmple qu'on peut avoir $ 0 > b \geq V_i \geq a_i \to -\infty $, il se peut qu'on ait uniform. sur une boule ouverte $ 0 >a_i/2 \geq V_i \geq a_i \to -\infty$, et en multipliant l'eq. par une fonction test et en integrant par parties, (il se peut aussi que les solutions $ u_i$ restent born\'ees), et on aurait, un terme qui tend vers $ -\infty $, et est egal a un terme born\'e uniform. Donc, pour qu'on soit sur qu'il n'y ait pas de contradiction, il faut supposer $ -\infty < a \leq V_i \leq b <0 $ uniform.

Dans tous les cas $ V $ doit etre regulier $ C^{\alpha}, \alpha >0 $. Pour le point de vue groupe topologique, ce n'est pas n\'ecessaire de supposer la regularit\'e $ C^{\alpha} $ uniforme. Alors que pour le point de vue analyse fonctionnelle, EDP, et des fonctions, il est necessaire de supposer la regularit\'e $ C^{\alpha}, \alpha >0$, uniforme.

\smallskip

d) Quand on considerent les solutions $ C^2$ et on veut obtenir, des estimations $ C^2$, il faut considrer une eq et non une inequation. Ce n'est plus possible de considerer des inequations dans ce cas.

\smallskip

e) Quand on considere la regularit\'e $ C^2 $, une inequation peut s'ecrire comme une equation. En posant $ V= [(\Delta u+Ru)/(u^{q-1})] \in C^0$ ou $ V=[(\Delta u +R)e^{-u}] \in C^0, \Delta =-\nabla^i\nabla_i$.

\smallskip

En ce qui concerne le theoreme 1 de l'article de 2003. Dans le raisonnement par l'absurde: pour obtenir une contradiction: On utilise l'in\'egalit\'e de Holder au lieu de l'in\'egalit\'e de Jensen, car les fonctions peuvent toucher $ 0 $. On ne peut utiliser l'in\'egalit\'e de Jensen, car une fonction convexe, doit l'etre dans un intervalle ouvert, or le proc\'ed\'e diagonal, implique l'existence d'une fonction $ v \geq 0 $ sur tout $ {\mathbb R}^n $ solution de $ \Delta v= W(x_0) v^{q-1}, 2< q \leq N, \Delta = -\nabla^i\nabla_i $, apres on utilise les integrales superficielles, sur les spheres. On ne peut pas appliquer la formule de Jensen, car la fonction $ v $ peut parfois s'annuler, elle n'est pas strict.positive, et ses valeurs sont dans $ [0, +\infty[$ qui n'est pas ouvert, et la convexit\'e necessite des intervalles ouverts.

\smallskip

////////////////////////////////////////////////////////////////////////////////

\smallskip

 En ce qui concerne le "pinching", la metric geometry, les groupes topologiques, la symetrie et le cadre de la cosmologie quantique: on en a parl\'e un peu dans le print: "Quelques remarques sur les vari\'et\'es, fonctions de Green et formule de Stokes": 

\smallskip

 1) Pour etre sur d'etre dans le cas negatif, il suffit de prendre la courbure scalaire $ R_g\equiv -1 $. Dans ce cas l'action de symetrie $ f $, on a $ R_gof =-1$. C'est leq de Yamabe dans le cas negatif. De plus on a $ \nabla (R_gof)=0 $, on voit bien que l'hypothese supplementaire $ C^{\alpha} $ est verif\'ee. On a, a la fois le "pinching" et tout ce qui concerne les groupes topologiques et la symetrie, et l'eq.de Yamabe, et la metric geometry, dans le cas negatif. Cela correspond au cas de ce qu'on a dit sur la cosmologie quantique avec $ \tau = constante \not = 0, -\tau^2 <0 $, l'eq. d'Einstein-Lichnerowicz dans le cas negatif. Pour etre sur d'etre dans ce cas, on prend $ \tau =constante \not = 0$. C'est l'eq. de Yamabe dans le cas negatif.le th 1 de l'article de 2003 est suffisant (on a suppos\'e les courbure holderiennes, ce qui est un peut genant en ce qui concerne les groupes topologiques, mais pour ce qui concerne le cadre, $ \tau =constante \not = 0, -\tau^2 <0 $, cette contrainte supplementaire disparait, car pour etre sur d'etre dans le cas negatif, il faut prendre cette condition sur $ \tau $).

\smallskip

On a aussi, le cas de la dimension 2, critique, dans le th.2. de l'article de 2003. En dimension 2. on a le "pinching", la metric geometry (surfaces hyperboliques), les groupes topologiques.

\smallskip

2) Dans le cas positif, dans l'article Bull.Sci.math. 2006. on a aussi le "pinching", dans un cas plus general incluant l'eq. de Yamabe, L'eq. de la courbure scalaire prescrite, la metric geometry.

\smallskip

3) Dans certains articles, Journ.fun.Anal. 2007. On avait le "pinching" pour l'eq. de Yamabe. Car le potential $ V\equiv constante=1 $. Ici aussi, on a la metric geometry.

Dans les articles avec potentiel $ V $ Lipschitzien, on a aussi la metric geometry. Car dans le cadre de la metric geometry, il faut que les fonctions  soit lipschitziennes: $ |V(x)-V(y)|\leq A  d_g(x,y), A>0, \forall x, y \in M $, on a des variations "lin\'eaires" des fonctions, difference entre $ V(x) $ et $ V(y) $ en termes de distance geodesique, de la metrique Riemannienne $ g $.

\smallskip

Donc, Lipschitzien est inclut dans la metric geometry et le "pinching".

\smallskip

Le fait d'avoir le potentiel $ 0 < a \leq V \leq b <+\infty $ est tres bien, on a bien le "pinching" et la metric geometry.

La condition de Lipschitz pour \underbar{la distance geodesique} est aussi tres interessant pour tout ce qui "pinching" et metric geometry. Avec la condition de Lipschitz on reste dans ces cadres, "pinching" et metric geometry.

\smallskip

/////////////////////////////////////////////////////////////////

\smallskip

En ce qui concerne la cosmologie quantique: 

\smallskip

a) pour etre sur d'etre dans le cas negatif, il faut prendre $ \tau =constante \not =0, -\tau^2 <0, \Psi=0$. Eq. Yamabe cas negatif.

\smallskip

b) pour etre sur d'etre dans le cas positif, avec symetrie, supersymetrie, il faut prendre $ \tau =0, \Psi=constante >0$. Eq; de Yamabe cas positif.

\smallskip

////////////////////////////////////////////////////////////////////

\smallskip

On a aussi:

\smallskip

c)l'article de 2003, met en avant le "pinching" et la metric geometry dans le cas negatif. Yamabe cas negatif.

\smallskip

c) L'article au Bull.Sci.math. 2006,  met en lumiere le "pinching" et la metrique geometry, dans le cas le plus general possible, dans le cas positif. Ici, la vari\'et\'e a un bord regulier et l'operateur $ \Delta+h $ est coercif.

\smallskip

d) La condition de Lipschitz pour la distance geodesique, est compatible avec la metric geometry. Pour ces articles, la vari\'et\'e peut avoir un bord non regulier et l'operateur principal est le laplacien, l'operateur $ \Delta +h $, n'est pas necessairement coercif.

\smallskip

e) Sur une vari\'et\'e Riemannienne, $ (M,g) $ de dimension $ n\geq 3 $, la condition (de la metric geometry): $ u\geq m >0 $, peut etre remplac\'e ou est equivalente \`a la condition suivante sur le volume conforme, ou en geometrie conforme et la "conformal geometry":

$$ \forall x\in M,\exists \rho_x >0,\,\, \forall \rho, \,\, 0 <\rho <\rho_x, \,\, \int_{B_{\rho}(x)} u^{2n/(n-2)} dV_g \geq \tilde m |B_{\rho}(x)|, \,\, \tilde m >0, $$

avec $ |B_{\rho}(x)|$ le volume de la boule centr\'ee en $ x\in M $ et de rayon $ \rho >0 $, pour la mesure Riemannienne relative \`a la metrique Riemannienne $ g $. La fonction $ u >0$ est continue sur $ M $. Par exemple l'article en dimension $ n=4 $, sur les vari\'et\'es Riemannienne de dimension 4. Ou la partie de l'article sur l'equation de Yamabe, en dimension 6. Cette condition est une condition sur le volume conforme ou de geometrie conforme. En termes d'hypotheses de geometrie conforme.

\smallskip

////////////////////////////////////////////////////////////////////////////////

\smallskip

Pour revenir au laplacien sur une boule de $ {\mathbb R}^n $, on a dit dans le print "Quelques remarques sur les vari\'et\'es, fonction de Green et formule de Stokes" qu'on obtenait directement son expression en coordonn\'ees polaires.

\smallskip

a) Explication 1: passage de coordonn\'ees cart\'esiennes en coordonn\'ees polaires: changement de coordonn\'ees. C'est ecrit dans le Dautray-Lions, le chapitre sur l'operateur de Laplace.

\smallskip

b) Explication 2: On a parl\'e de metrique euclidienne: $ dx^2=dr^2+r^2d\theta^2 $ et formule usuelle dans ne carte polaire, $ (r,\theta)$, comme c'est expliqu\'e dans le livre de Hebey.

Consid\'erons une boule $ B_r(0) $ et une fonction $ u \in C^2$ de $ {\mathbb R}^n $, alors la distance $ d $ de $ {\mathbb R}^n $ et la distance geodesique $ d_{\delta} $ pour la metrique Riemannienne euclidienne $ \delta $, sont egales: $ d=d_{\delta} $. La fonction $ u $ est $ C^2 $ sur $ B_r(0) $ en tant que vari\'et\'e Riemannienne munie de la metrique Riemannienne euclidienne. On ecrit alors:

$$ \Delta u= \Delta_{\delta} u, {\rm \,\, en\, tant \, que \, operateur \, Riemannien} $$

on a le "comme si ", on considerait une structure de vari\'et\'e Riemmanienne sur la boule $ B_r(0) $, alors le changement de variable s'interprete, comme, chercher l'expression du laplacien dans une "carte polaire", localement, en coordonn\'ees, ici, en coordonn\'ees polaires. D'o\`u la fait  de considerer, la metrique euclidienne: $ dx^2=dr^2+r^2d\theta^2 $. Apres on utilise le resultat du livre de Hebey: le laplacien en coordonn\'ees, ici en polaires.

\smallskip

Une fonction est $ C^2$ au sens usuel $=$ une fonction $ C^2$ au sens des cartes et de la metrique Riemannienne, de la connexion Riemannienne Euclidienne (connexion de Levi-Civita pour la metrique euclidienne, gradient usuel). Ici, il y a deja une carte, $ (B_r(0), id) $. Tout se confond, vari\'et\'e et ouvert de carte.

\smallskip

C'est la o\`u il faut supposer que l'ouvert est muni de la metrique Riemannienne Euclidienne, $ \delta $.

\smallskip

///////////////////////////////////////////////////////////////////////////////

\smallskip

{\bf Sur l'article donnant l'unicit\'e au.Bull.Sci.Math.de 2009:}

\smallskip

Remarque sur le th. d'Olivier Druet donnant la compacit\'e et le th. de l'article de 2009 au Bull.Sci.Math. donnant la compacit\'e pour $ \epsilon \to 0 $.

1) La preuve de Druet est vraie quand $ 1\geq \epsilon \geq \tilde \epsilon_0 >0, \forall \tilde \epsilon_0 >0 $. 

2) La preuve de l'article au Bull.Sci.Math. est vraie lorsque $ 0 < \epsilon \leq \epsilon_0 $: $ \exists  m_0 >0, \exists \epsilon_0 >0, $ tels que pour $ 0 < \epsilon \leq \epsilon_0 $, on ait: $ \max_M u_{\epsilon} \leq m_0 $.

3) En suite on a l'unicit\'e: $ \exists \epsilon_1 >0, 0 <\epsilon_1 \leq \epsilon_0 $ tel que pour $ 0 <\epsilon \leq \epsilon_1 $, on ait: $ u_{\epsilon} \equiv \epsilon^{(n-2)/4} $.

\smallskip

La remarque est la suivante: on ne peut pas dire que le resultat de Druet est vrai, par rapport au resultat du Bull.Sci.Math. 2009, pour $ 0< \epsilon_1 \leq \epsilon\leq \epsilon_0 $. Car des qu'on parle de $ \epsilon_0 $, cela nous amene au point 2) precedent qui dit qu'on a la compacit\'e pour $ 0 <\epsilon \leq \epsilon_0 $, or la preuve(la methode) de Druet n'est pas vraie pour $ 0 < \epsilon \leq \epsilon_0 $, donc on ne peut pas parler de la preuve de Druet dans ce cas, car elle n'est pas vraie. On ne peut pas distinguer les cas $ 0 <\epsilon \leq \epsilon_1 $ et $ 0 < \epsilon_1 \leq \epsilon \leq \epsilon_0 $. Car des qu'on parle de $ \epsilon_0$ cela fait reference au point 2) precedent qui dit qu'on a la compacit\'e dans ce cas et que la preuve de Druet n'est pas vrai dans ce cas.

\smallskip

On ne peut pas parler de $\epsilon_0$ et $\epsilon_1$ ou de $ \epsilon_0$, et du resultat de Druet, au meme temps, ce n'est pas compatible.

\smallskip

(Par exemple, considerer $ \epsilon_0 >0$ du point 2), c'est prendre en compte les fonctions $ 
u_{\epsilon} $ avec $ \epsilon >0 $ voisin de $ 0 $, ce qui n'est pas le cas du resultat de Druet).
\smallskip

Tout ceci pour dire qu'on a bien un resultat de compacit\'e dans l'article du Bull.Sci.Math. de 2009, sans faire reference au resultat de Druet et sans faire reference \`a l'unicit\'e. L'unicit\'e vient apres.

\smallskip

On peut mettre un $ \epsilon_2 >0, \epsilon_2 \geq \epsilon_0 \geq \epsilon_1 $, pour l'in\'egalit\'e $ \sup \times \inf $ sur une vari\'et\'e compacte sans bord.

\smallskip

Dans cet article, on a l'in\'egalit\'e optimale locale (vari\'et\'e Riemannienne non necessairement sans bord), $ \sup \times \inf $ pour une equation de Schrodinger: quantique non relativiste. C'est reellement quantique non relativiste, car, cette equation, ne provient pas, par exemple, d'Einstein-Lichnerowicz. C'est une equation de Schrodinger. On a alors un exemple d'eq. quantique non relativiste. Eq. de Schrodinger, non relativiste, apparait en optique. 

\smallskip

Mais quand on met une condition sur la courbure scalaire, pour avoir l'inegalit\'e optimale $ \sup \times \inf $, cette equation devient relativiste, car, la condition implique la dependance en la courbure scalaire.

\smallskip

{\bf Sur la 2\`eme note aux Comptes Rendus Math. Acad. Sci. Paris. 2006:}

\smallskip

On a l'in\'egalit\'e de Brezis-Gallouet, qui prend en compte les deriv\'ees d'ordre 2, $ H^2_2 $ et la norme $ \sup $. On peut voir dans le dernier th de la note de 2006 (on utilise la borne uniforme des solutions au voisinage du bord, de l'equation ($ \Delta^2 u= u^{p-\epsilon}, p=\frac{(n+4)}{(n-4)}, n\geq 5 $ dans $ \Omega $ et avec les conditions au bord $ u=\Delta u=0$ sur $ \partial \Omega $, sur un ouvert regulier strictement convexe, par exemple une boule ou un ellipsoide)):

$$ \sup_{\Omega} u \times \inf_K u \geq c(K,\Omega,n) \int_{\Omega} (\Delta u)^2 dx \geq \tilde c(K, \Omega, n) >0, $$

C'est une in\'egalit\'e liant (et melange) les normes essentielles et la norme $ H^2_2(\Omega) $. Operateur bilaplacien. Condition au bord.

\smallskip


\smallskip

\begin{thebibliography}{99} 

\bibitem{1}{T. Aubin. Some Nonlinear Problems in Riemannian Geometry. Springer-Verlag, 1998.}

\bibitem{2}{C. Bandle. Isoperimetric Inequalities and Applications. Pitman, 1980.}



\bibitem{3}{H. Brezis, YY. Li. Some nonlinear elliptic equations have only constant solutions. J. Partial Diff. Eqs. Vol 19, no 3, (2006), 208–217.}

\bibitem{4}{H. Brezis, F. Merle. Uniform estimates and Blow-up behavior
for solutions of $ -\Delta u=V(x) e^u $ in two dimensions. Commun. in
Partial Differential Equations, 16 (8 and 9), 1223-1253(1991).}

\bibitem{5}{H. Brezis, L. Nirenberg. Positive solutions of nonlinear elliptic equations involving critical Sobolev exponents. Comm. Pure and Appl. Math. Volume 36, 4, pp 437-477. (1983).}

\bibitem{6}{Crandall.M.G, Rabinowitz.P.H. Some continuation and variational methods for positive solutions of nonlinear
elliptic eigenvalue problems. Arch.rational Mech.Anal, 58, 207-218, 1975.}

\bibitem{5}{C.C.Chen. C.S.Lin. Estimates of the conformal scalar curvature equation via the method of moving planes. Comm.
Pure Appl. Math. 50, (1997) 0971-1017.}

\bibitem{5}{C.C.Chen. C.S.Lin. Estimates of the conformal scalar curvature equation via the method of moving-planes.II. J. Differential Geom. 49(1): 115-178 (1998).}

\bibitem{6}{Davis.M.W. A hyperbolic 4-manifold Proc. Amer. Math. Soc.93 (1985), pp.325-328.}


\bibitem{7}{Druet.O. Generalized scalar curvature type equations on compact Riemannian manifolds. C. R. Acad. Sci. Paris, 327, serie I, p 377-381. 1998.}

\bibitem{8}{ Dupaigne. Stable Solutions of Elliptic Partial Differential Equations (Monographs and Surveys in Pure and Applied Mathematics). 2011.}

\bibitem{9}{J. Escobar. R. Schoen. Conformal metrics with prescribed scalar curvature. Inventiones mathematicae (1986) 86, page 243-254.}

\bibitem{10}{ Gallot. Hulin. Lafontaine. Riemannian Geometry.}

\bibitem{11}{E. Hebey. Introduction a l'analyse non-lineaire sur les varietes. Diderot Editions.}

\bibitem{12}{E. Hebey, M. Vaugon. Existence and multiplicity of nodal solutions for nonlinear elliptic equations with critical Sobolev growth. Journal of Functional Analysis, 119, 298-318, 1994}

\bibitem{13}{D. Holcman.  Solutions nodales sur les varietes riemanniennes non localement conformement plates a bord. Comment. Math. Helv. 76 (2001), no. 3, 373-387.}

\bibitem{14}{D.Huybrechts. Lecture on K3 surfaces.}

\bibitem{15}{D. Joyce. Constant scalar curvature metrics on connected sums. International Journal of Mathematics and Mathematical Sciences.(2003), Issue 7, pp 405-450.}

\bibitem{16}{Karp. L, Pinsky. M . The first eigenvalue of a small geodesic ball in a Riemannian manifold. Bull. Sci. Math. 
(2)111 (1987), no. 2,229-239.}

\bibitem{17}{Li.YY, Zhu.M. Yamabe type equations on three dimensional Riemannian manifolds. Communications in Contemporary Math. 1 (1999), 1-50.}

\bibitem{18}{J. Lafontaine.Introduction a la geometrie Differentielle.2010.}


\bibitem{19}{J.G. Ratcliffe.S.T. Tschantz. On the Davis hyperbolic 4-manifold. Topology and its Applications
Volume 111, Issue 3. pp 327-342}

\bibitem{20}{Troyanov. M. Prescribing curvature on compact surfaces with conical singularities. Trans. Amer. Math. Soc. 324 (1991), no. 2, 793-821.}

\bibitem{21}{Troyanov. M. Un principe de concentration-compacit\'e pour les suites de surfaces Riemanniennes. Annales de l'IHP, Analyse non lin\'eaire. Vol 8, no 5, (1991), pp 419-441.}

\end{thebibliography}

\begin{thebibliography}{99} 

\bibitem{1}{A.D. Alexandrov. Uniqueness thoerems for surfaces in the large. V. Vestnik Leningrad Univ. Mat. Astronom 13, 5-8 (1958); Amer. Math. Soc. Transl. 21, 412-416 (1962).}

\bibitem{2}{T. Aubin. Some Nonlinear Problems in Riemannian Geometry. Springer-Verlag, 1998.}

\bibitem{3}{ Ambrosio. L, Fusco. N, Pallara, D. Functions of Bounded variations and Free discontinuity Problems, Oxford Press. 2000.}
 
\bibitem{4}{ S.S.Bahoura. Majorations du type $ \sup u \times \inf u \leq c $ pour l'\'equation de la courbure scalaire sur un ouvert de $ {\mathbb R}^n, n\geq 3 $. J. Math. Pures. Appl.(9) 83 2004 no, 9, 1109-1150.}

\bibitem{5}{Bahoura.S.S. Lower bounds for sup+inf and sup * inf and an extension of Chen-Lin result in dimension 3.  Acta Math. Sci. Ser. B Engl. Ed.  28  (2008),  no. 4, 749-758}

\bibitem{6}{Bahoura.S.S. About Brezis Merle problem with Holderian condition. Mathematica Aeterna, Vol 4, no 1, 13-25, 2014.}

\bibitem{7}{Bahoura.S.S. A uniform boundedness result for solutions to the Liouville type equation with boundary singularity. J. Math. Sci. Univ. Tokyo, 23, no 2, 487-497. 2016.}

\bibitem{8}{Bahoura.S.S. A compactness result for an equation with Holderian condition. Commun. Math. Anal. Vol 21, no 1, 23-34, 2018.}

\bibitem{9}{Bahoura.S.S. A local uniform boundedness result for an elliptic equation. hal-02495444. 2020.}

\bibitem{10}{C. Bandle. Isoperimetric Inequalities and Applications. Pitman, 1980.}

\bibitem{11}{Berestycki. H. Nirenebrg. L. Varadhan. The principal eigenvalue and maximum principle for second-order elliptic operators in general domains. Comm. Pure Appl. Math. 47 (1994), no. 1, 47-92}

\bibitem{12}{L. Boccardo, T. Gallouet. Nonlinear elliptic and parabolic
equations involving measure data. J. Funct. Anal. 87 no 1, (1989),
149-169.} 

\bibitem{13}{H. Brezis. Analyse Fonctionelle et Applications. 1983.}

\bibitem{14}{Brezis. H, Marcus. M, Ponce. A. C. Nonlinear elliptic equations with measures revisited. Mathematical aspects of nonlinear dispersive equations, 55-109, Ann. of Math. Stud., 163, Princeton Univ. Press, Princeton, NJ, 2007.} 

\bibitem{15}{H. Brezis, F. Merle. Uniform estimates and Blow-up behavior
for solutions of $ -\Delta u=V(x) e^u $ in two dimensions. Commun. in
Partial Differential Equations, 16 (8 and 9), 1223-1253(1991).}

\bibitem{16}{H. Brezis, W. A. Strauss. Semi-linear second-order elliptic equations in L1. J. Math. Soc. Japan 25 (1973), 565-590.}

\bibitem{17}{H. Brezis, YY. Li , I. Shafrir. A sup+inf inequality for some
nonlinear elliptic equations involving exponential
nonlinearities. J.Funct.Anal.115 (1993) 344-358.
}

\bibitem{18}{L. Caffarelli, B. Gidas, J. Spruck. Asymptotic symmetry and local
behavior of semilinear elliptic equations with critical Sobolev
growth. Comm. Pure Appl. Math. 37 (1984) 369-402.
}

\bibitem{19}{Chen. W,  Li, C. Classification of solutions of some nonlinear elliptic equations. Duke Math. J. 63 (1991), no. 3, 615-622.}
\bibitem{20}{W. Chen, C. Li. A priori Estimates for solutions to Nonlinear
Elliptic Equations. Arch. Rational. Mech. Anal. 122 (1993) 145-157.}

\bibitem{21}{C-C. Chen, C-S. Lin. Local behavior of singular positive solutions of semilinear elliptic equations with Sobolev exponent. Duke Math. J. 78 (1995), no. 2, 315-334.}

\bibitem{22}{C-C.Chen, C-S. Lin. Estimates of the conformal scalar curvature
equation via the method of moving planes. Comm. Pure
Appl. Math. L(1997) 0971-1017.}

\bibitem{23}{C-C. Chen, C-S. Lin. A sharp sup+inf inequality for a
nonlinear elliptic equation in ${\mathbb R}^2$. Commun. Anal. Geom.
6, No.1, 1-19 (1998).}

\bibitem{24}{de Figueiredo, D. G,  Lions, P.-L, Nussbaum, R. D. A priori estimates and existence of positive solutions of semilinear elliptic equations. J. Math. Pures Appl. (9) 61 (1982), no. 1, 41-63.}

\bibitem{25}{ O. Druet, E. Hebey, F.Robert, Blow-up theory in Riemannian Geometry, Princeton University Press 2004.}

\bibitem{26}{ B. Gidas, W-M. Ni, L. Nirenberg. Symmetry and Related Properties via the Maximum Principle. Commun. Math. Phys. 68, 209-243 (1979).}

\bibitem{27}{ D. Gilbarg, N.S. Trudinger. Elliptic Partial Differential Equations of Second order, Berlin Springer-Verlag, Second edition, Grundlehern Math. Wiss.,224, 1983.}

\bibitem{28}{Han Z-C. Asymptotic approach to singular solutions for nonlinear
elliptic equations involving critical Sobolev
exponent. Ann. Inst. Henri Poincar{\'e}. Analyse non lineaire 8 (1991) 159-174.}

\bibitem{29}{E. Hebey, Analyse non lineaire sur les Vari\'et\'es, Editions Diderot.}
\bibitem{30}{O. Kavian. Introduction \`a la Th\'eorie des points critiques. Springer-Verlag.1993.}

\bibitem{31}{N. Korevaar, F. Pacard, R. Mazzeo, R. Schoen. Refined asymptotics for constant scalar curvature metrics with isolated singularities.  Invent. Math.  135  (1999),  no. 2, 233--272.}
\bibitem{32}{YY. Li. Harnack Type Inequality: the Method of Moving Planes. Commun. Math. Phys. 200,421-444 (1999).}

\bibitem{33}{YY.Li, L.Zhang. Liouville type theorems and Harnack type inequalities for semilinear elliptic equations, Journal d'Analyse Mathematique, 90 (2003), 27-87.}


\bibitem{34}{YY. Li, L. Zhang. A Harnack type inequality for the Yamabe equation in low dimensions.  Calc. Var. Partial Differential Equations  20  (2004),  no. 2, 133--151.}

\bibitem{35}{L. Ma, J-C. Wei. Convergence for a Liouville equation.
Comment. Math. Helv. 76 (2001) 506-514.}

\bibitem{36}{Nazaraov, S, A.  Sweers, G. A hinged plate equation and iterated Dirichlet Laplace operator on domains with concave corners. Journal of Differential Equations. 233, 1, 2007, 151-180.}

\bibitem{37}{J. Serrin, A symmetry problem in potential theory, Arch. Rational Mech. Anal. 43 (1971), 304-318.}
\end{thebibliography}

\smallskip

\begin{thebibliography}{99} 

\bibitem{1}{Bahoura.S.S. Estimations du type $ \sup \times \inf $ sur une vari\'et\'e compacte. Bull.Sci.math. 130 (7), 2006, pp 624-636.}

\bibitem{2}{Bahoura.S.S. Estimations uniformes pour l'equation de Yamabe en dimensions 5 et 6. J. Funct. Anal.  242  (2007),  no. 2, 550-562.}

\bibitem{3}{Bahoura.S.S. Harnack inequalities for Yamabe type equations.  Bull. Sci. Math.  133  (2009),  no. 8, 875-892}

\bibitem{4}{Bahoura.S.S. $ \sup \times \inf $ inequality on manifold of dimension 3. Math. Aeterna, 1 (01) (2011), pp. 13-26}

\bibitem{5}{Bahoura.S.S. A uniform estimate for scalar curvature equation on manifolds of dimension 4. J. Math. Anal. Appl.
388 (2012), no. 1, 386-392}

\bibitem{6}{Bahoura.S.S. An estimate on Riemannian manifolds of dimension 4. Analysis in Theory and Applications. No 32, 3, (2016) pp 272-282.}

\bibitem{7}{Bahoura.S.S. Some uniform estimates for scalar curvature type equations. Pacific. J. Math, vol 301, no 1, 2019, pp 55-65.}

\bibitem{8}{Druet.O. Compactness for Yamabe metrics in low dimensions. Int. Math. Res. Not. 23, 1143-1191 (2004)}

\bibitem{9}{Li.YY, Zhang.L. Compactness of solutions to the Yamabe problem. II, Calculus of Variations and PDEs 24 (2005), 185-237.}

\bibitem{10}{F.C.Marques. A priori estimates for the Yamabe problem in the non-locally conformally flat case. J. Diff. Geom. 71(2), pp 315-346, 2005.}

\bibitem{11}{G.Tarantello. A Harnack Inequality for Liouville-type Equations with Singular Sources. Indiana University Mathematics Journal. 54, no 2 (2005) 599-615.}


\end{thebibliography}

\begin{thebibliography}{99} 

\bibitem{1}{Carlotto.A. The general relativistic constraints equations. Living reviews in Relativity. 24, 2, 2021. Springer.}

\bibitem{2}{Chen.C.C. Lin.C.S. Wang.G. Concentration Phenomena of two Vortex Solutions in a Chern-Simons Model. Ann. Scuola.Norm.Sup.Pisa.Cl.Sci.(5) Vol 3, (2004), pp 367-397.}

\bibitem{3}{Choquet-Bruhat.Y. Th\'eor\`eme global d'unicit\'e pour les solutions des equations d'Einstein. Bulletin de la S.M.F. tome 96, 1968, 181-192.}

\bibitem{4}{Choquet-Bruhat.Y, Isenberg.J, Pollack.D. The Einstein-scalar field constraints on asymptotically Euclidean manifolds. arXiv:gr-qc/0506101.}

\bibitem{5}{Choquet-Bruhat.Y, Isenberg.J, Pollack.D. The constraint equations for the Einstein-scalar field system on compact manifolds. arXiv:gr-qc/0610045.}

\bibitem{6}{Choquet-Bruhat.Yvonne, Isenberg.James, York.James. W. Einstein Constraints on Asymptotically Euclidean Manifolds. arXiv:gr-qc/9906095.}

\bibitem{7}{Choquet-Bruhat.Y, Geroch.R. Global aspects of the Cauchy problem in general relativity, Comm. Math. Phys. 14 (1969), no. 4, 329-33.}

\bibitem{8}{Collion.S, Vaugon.M. A new approach to Kaluza-Klein Theory. arXiv:1709.04172v3.}

\bibitem{9}{Four\`es-Bruhat.Y. Th\'eor\`eme d'existence pour certains syst\`emes d'\'equations aux d\'eriv\'ees
partielles non lin\'eaires, Acta Math. 88 (1952), 141-225.}

\bibitem{10}{Horvat.R, Krcmar.M, Lakic.B. Recent searches for solar axions and large extra dimensions.
arXiv:hep-ph/0112224.}

\bibitem{11}{Horvat.R, Krcmar.M, Lakic.B. CERN \,\, Axion \,\, Solar \,\, Telescope \,\, as a probe of large extra dimensions. arXiv:astro-ph/0312030.}

\bibitem{12}{Hillen.T, Painter.K. Global existence for a parabolic chemotaxis model with prevention of overcrowding. Adv. Appl. Math. 26, 280-301, 2001.}

\bibitem{13}{Keller. E. F, Segel, L.A. Initiation of slime mold aggregation viewed as an instability.J. Theoret. Biol., 26 (1970), pp. 399-415.}

\bibitem{14}{Lin. C.S, Ni. W.M, Takagi.I Large amplitude stationary solutions to a chemotaxis system. Journ. Diff. Equations, volume 72, no 1, 1988, pp 1-27.}

\bibitem{15}{Parker. T.H. Gauge theories on four dimensional Riemannian manifolds. Comm. Math. Phys. 85, 563-602, 1982.}

\bibitem{16}{Premoselli.B. Stability and instability of the Einstein-Lichnerowicz constraint system. arXiv:1502.04233v1 }

\end{thebibliography}

\begin{thebibliography}{99} 

\bibitem{1}{S.S.Bahoura. Diff\'erentes estimations du  $ \sup u \times \inf u $ pour l'\'equation de la courbure scalaire prescrite en dimension $ n\geq 3 $. J. Math. Pures Appl. (9), 82 (1) (2003), pp. 43-66 }

\bibitem{2}{ S.S.Bahoura. Majorations du type $ \sup u \times \inf u \leq c $ pour l'\'equation de la courbure scalaire sur un ouvert de $ {\mathbb R}^n, n\geq 3 $. J. Math. Pures. Appl.(9) 83 2004 no, 9, 1109-1150.}

\bibitem{3}{Bahoura.S.S. Harnack inequalities for Yamabe type equations.  Bull. Sci. Math.  133  (2009),  no. 8, 875-892}

\bibitem{4}{S.S.Bahoura. In\'egalit\'es de Harnack pour les op\'erateurs elliptiques d'ordre 2 et 4 et ph\'enom\`ene de concentration. C.R.Acad.Sci.Paris, Ser.I, 342, (2006), 755-758. }

\bibitem{5}{S.S.Bahoura. In\'egalit\'es de Harnack et ph\'enom\`ene de concentration. J. Math. Soc. Japan. Vol 59, No 4, (2007) pp, 1011-1030.}

\bibitem{6}{Druet.O. Compactness for Yamabe metrics in low dimensions. Int. Math. Res. Not. 23, 1143-1191 (2004)}

\end{thebibliography}
\end{document}